\documentclass[12pt,final]{amsart}
\usepackage[utf8]{inputenc}
\usepackage[T1]{fontenc}
\usepackage{lmodern}
\usepackage{amsmath,amssymb,amsthm,mathtools}
\usepackage[margin=1in]{geometry}
\usepackage{mathrsfs,microtype}
\usepackage[shortlabels]{enumitem}
\usepackage[ocgcolorlinks,linkcolor=blue,citecolor=blue,urlcolor=blue]{hyperref}
\hypersetup{pdftitle={Sharp entanglement and decay thresholds for fractional elliptic operators},pdfauthor={Yi-Hsuan Lin and Hongyu Liu}}
\numberwithin{equation}{section}
\allowdisplaybreaks[2]
\mathtoolsset{showonlyrefs}
\newtheorem{theorem}{Theorem}[section]
\newtheorem{lemma}[theorem]{Lemma}
\newtheorem{proposition}[theorem]{Proposition}
\newtheorem{corollary}[theorem]{Corollary}
\theoremstyle{definition}
\newtheorem{definition}[theorem]{Definition}

\newcommand{\R}{\mathbb R}
\newcommand{\C}{\mathbb C}
\newcommand{\N}{\mathbb N}
\newcommand{\Z}{\mathbb Z}
\newcommand{\Nzero}{\mathbb N_0}
\newcommand{\eps}{\varepsilon}
\newcommand{\cO}{\mathcal O}
\newcommand{\cE}{\mathscr E'}
\newcommand{\cS}{\mathscr S}
\newcommand{\HL}{\mathcal H_L}
\newcommand{\norm}[1]{\lVert #1\rVert}
\newcommand{\ip}[2]{\langle #1,#2\rangle}
\newcommand{\para}[1]{\vspace{3mm}\noindent\textbf{#1.}}
\DeclareMathOperator{\supp}{supp}
\DeclareMathOperator{\dist}{dist}
\DeclareMathOperator{\Rea}{Re}
\DeclareMathOperator{\Ima}{Im}
\DeclareMathOperator{\Div}{div}
\DeclareMathOperator{\Dom}{Dom}

\title[Entanglement and decay thresholds]{Sharp entanglement and decay thresholds for fractional elliptic operators}
\author[Y.-H. Lin]{Yi-Hsuan Lin}
\address{Department of Applied Mathematics, National Yang Ming Chiao Tung University, Hsinchu, Taiwan}
\email{yihsuanlin3@gmail.com}
\author[H. Liu]{Hongyu Liu}
\address{Department of Mathematics, City University of Hong Kong, Kowloon, Hong Kong, China}
\email{hongyu.liuip@gmail.com, hongyliu@cityu.edu.hk}
\date{}

\begin{document}
	
	\begin{abstract}
		We establish sharp entanglement principles for fractional powers of uniformly elliptic divergence-form operators $L$ on $\mathbb{R}^n$, $n\ge 2$, with smooth real symmetric coefficients. Let $u_1,\ldots,u_N$ and the sum $\sum_{k=1}^N L^{s_k}u_k$ vanish on a common nonempty open set. Assuming that $L$ coincides with the Laplacian outside a compact set, we prove that root-exponential decay of these functions forces every $u_k$ to vanish globally, under sharp dimension-dependent conditions on the fractional orders $s_k$. The exponent $1/2$ is optimal within the stretched-exponential scale, and the sharp logarithmic endpoint is $e^{-a\sqrt{r}/\log r}$ in every dimension. Our approach relies on a novel resolvent characterization of entanglement, which we use both to separate the fractional orders at zero energy and to construct exact cancellations from outgoing waves. For perturbations extending to infinity, this wave construction yields two-dimensional counterexamples with smooth compactly supported inputs and coefficient fields converging to the identity in every fixed unweighted $C^k$ norm. Finally, for a known underlying operator, a single exterior measurement determines the number, orders, and coefficients of the fractional terms, together with a continuous potential.
	\end{abstract}

	\maketitle
	\tableofcontents
	
	\section{Introduction}\label{sec:introduction}
	
	The entanglement principle asks whether a linear relation among different fractional powers of an elliptic operator, acting on functions that vanish on an open set, forces every function to vanish globally. The principle was introduced and proved on closed Riemannian manifolds in \cite{FKU24}. On $\R^n$, the position of the observation set, the decay of the functions, and the coefficients at infinity lead to additional obstructions.
	
	For the fractional Laplacian on $\R^n$, \cite{FL24} established entanglement under exponential decay, with integer differences of orders excluded in even dimensions and half-integer differences excluded in odd dimensions. The optimal decay assumption and the necessity of the additional odd-dimensional restriction were left open. Related entanglement principles for nonlocal parabolic operators and on hyperbolic spaces were studied in \cite{LLY25,Lin26}.
	
	We distinguish exterior observation from observation in a bounded component of the complement of the supports. For compactly supported distributions, an observation set meeting the unbounded component gives entanglement for general smooth uniformly elliptic coefficients whenever the orders have pairwise noninteger differences. Observation in a bounded component permits additional half-integer cancellations in odd dimensions. For compact perturbations of the Laplacian and arbitrary observation, we prove the sharp order restrictions and determine the exact threshold $\gamma=1/2$ in the scale $e^{-a|x|^\gamma}$. In particular, Schwartz decay alone does not suffice.
	
	The coefficients at infinity give another obstruction. In dimension two, we construct compactly supported counterexamples for coefficients converging to the identity in every fixed unweighted $C^k$ norm. The fractional orders satisfy the restrictions of the positive theorem for compact perturbations. Smallness in these norms cannot replace the assumption on the coefficients at infinity.
	
	The proofs use a novel resolvent characterization of entanglement that applies to general smooth uniformly elliptic coefficients. On an open set where the inputs vanish, Proposition~\ref{prop:resolvent} identifies a single fractional relation with the family of resolvent identities \eqref{eq:intro-resolvent}, in which each order $s$ enters only through the scalar weight $\sin(\pi s)\lambda^s$. Applied in the converse direction to an outgoing wave that stays outside the observation ball, the same characterization gives the exact cancellation in \eqref{eq:wave-cancellation}. It also yields unique continuation for finite sums of distinct positive noninteger powers acting on one $L^2$ function, without additional decay assumptions or restrictions on differences of orders. The inverse problems considered below use these results for general coefficients.
	
	Throughout the paper, let
	\begin{equation}\label{eq:L}
		L=-\Div(A(x)\nabla)\quad\text{on }\R^n,\quad n\ge2,
	\end{equation}
	where $A\in C^\infty(\R^n;\R^{n\times n})$ is real symmetric and uniformly elliptic. More precisely, we assume that there exist constants $0<a_0<a_1<\infty$ such that
	\begin{equation}\label{eq:ellipticity}
		a_0|\xi|^2\le A(x)\xi\cdot\xi\le a_1|\xi|^2\quad\text{for all }x,\xi\in\R^n.
	\end{equation}
	We use the nonnegative self-adjoint realization of $L$ on $L^2(\R^n)$ and define its fractional powers by the spectral theorem. All powers in an entanglement equation are powers of the same operator.
	
	Let $N\in\N$, let $s_1,\ldots,s_N$ be positive nonintegers, and let $b_1,\ldots,b_N\in\C\setminus\{0\}$. Given a nonempty open set $\cO\subset\R^n$, we consider
	\begin{equation}\label{eq:entanglement}
		u_1|_{\cO}=\cdots=u_N|_{\cO}=0,\quad\sum_{k=1}^N b_kL^{s_k}u_k\Big|_{\cO}=0.
	\end{equation}
	The equations are understood in the sense of distributions. The distribution spaces and the corresponding fractional powers are specified in Section~\ref{sec:preliminaries}.
	
	\subsection{Entanglement principles}
	
	We denote by $\cE(\R^n)$ the space of compactly supported distributions on $\R^n$, and we begin with exterior observation.
	
	\begin{theorem}[Exterior observation]\label{thm:exterior}
		Let $A$ satisfy \eqref{eq:ellipticity}, and let $u_1,\ldots,u_N\in\cE(\R^n)$. Set $K=\bigcup_{k=1}^N\supp u_k$, and let $E$ be the unbounded connected component of $\R^n\setminus K$. Let $\cO\subset\R^n$ be a nonempty open set satisfying $\cO\cap E\ne\varnothing$, and assume that
		\begin{equation}\label{eq:exterior-orders}
			s_k-s_\ell\notin\Z\quad\text{whenever }k\ne \ell.
		\end{equation}
		If \eqref{eq:entanglement} holds, then $u_k=0$ on $\R^n$ for every $k$.
	\end{theorem}
	
	The compact set $K$ may have holes and need not have a regular boundary; the open set $\cO$ need only intersect the unbounded connected component of $\R^n\setminus K$. The condition \eqref{eq:exterior-orders} is necessary. Indeed, if $t=s+m$ with $m\in\N$, a nonzero $v\in C_c^\infty(\R^n)$ and $u=-L^mv$ satisfy $L^su+L^tv=0$ globally. Equal orders also allow cancellation by taking $u=-v$.
	
	For an arbitrary nonempty open set $\cO$, we assume in addition that
	\begin{equation}\label{eq:compact-perturbation}
		A-I\in C_c^\infty(\R^n;\R^{n\times n}).
	\end{equation}
	We impose no smallness assumption on $A-I$. We use the following distributional decay condition, as in \cite{FL24}.
	
	\begin{definition}\label{def:decay}
		For $\gamma>0$, let $\mathscr D_\gamma$ consist of those $u\in H^{-r}(\R^n)$, for some $r\ge0$, for which there are constants $a,C>0$ such that
		\begin{equation}\label{eq:distribution-decay}
			|\ip{u}{\phi}|\le Ce^{-aR^\gamma}\norm{\phi}_{H^r(\R^n)}
		\end{equation}
		for all $R\ge1$ and all $\phi\in C_c^\infty(\R^n\setminus\overline{B_R})$.
	\end{definition}
	
	Here and below, $B_R$ is the open ball of radius $R$ centered at the origin. Every compactly supported distribution belongs to $\mathscr D_\gamma$ for every $\gamma>0$. For $L^2$ functions, the weighted condition $e^{a|x|^\gamma}u\in L^2(\R^n)$ implies $u\in\mathscr D_\gamma$. The choice $\gamma=1/2$ gives root-exponential decay.
	
	\begin{theorem}[Arbitrary observation]\label{thm:arbitrary}
		Let $A$ satisfy \eqref{eq:ellipticity} and \eqref{eq:compact-perturbation}, let $\cO\subset\R^n$ be nonempty and open, and let $u_k\in\mathscr D_{1/2}$ for every $k$. Assume one of the following:
		\begin{enumerate}[\rm(i)]
			\item\label{item_1_thm_arb} $s_k-s_\ell\notin\Z$ for $k\ne \ell$ when $n$ is even;
			\item\label{item_2_thm_arb} $s_k-s_\ell\notin\frac12\Z$ for $k\ne \ell$ when $n$ is odd.
		\end{enumerate}
		If \eqref{eq:entanglement} holds, then $u_k=0$ on $\R^n$ for every $k=1,\ldots,N$.
	\end{theorem}
	
	Theorem~\ref{thm:arbitrary} relaxes the exponential decay assumption in \cite{FL24} to root-exponential decay and allows arbitrary smooth uniformly elliptic changes of the coefficients in a compact region. Theorem~\ref{thm:sharpness} proves that both the order restrictions and the exponent $1/2$ are sharp.
	
	Even for $A=I$ and compactly supported sources, the position of $\cO$ changes the admissible orders in odd dimensions. Differences in $\frac12+\Z$ are allowed under the exterior observation hypothesis of Theorem~\ref{thm:exterior}, but they admit compactly supported counterexamples when $\cO$ lies in a bounded component of the complement of the supports. These examples are given in Theorem~\ref{thm:sharpness}\ref{item_1_thm_sharp}.
	
	For a finer decay condition, put
	\begin{equation}\label{eq:def_w_p_r}
		w_p(r)=\frac{\sqrt r}{\log^p(e+r)}, \quad \text{for }p\ge 0 \text{ and }r\geq 0.
	\end{equation}
	We denote by $\mathscr D_{\log,p}$ the class in Definition~\ref{def:decay} with $e^{-aR^\gamma}$ replaced by $e^{-aw_p(R)}$. In particular, $\mathscr D_{\log,0}=\mathscr D_{1/2}$, and $e^{aw_p(|x|)}u\in L^2$ implies $u\in\mathscr D_{\log,p}$.
	
	\begin{proposition}\label{prop:odd-logarithmic}
		Let $n\ge3$ be odd, and let $A$ satisfy \eqref{eq:ellipticity} and \eqref{eq:compact-perturbation}. Suppose that $u_k\in\mathscr D_{\log,p}$ for every $k$, where $0\le p\le1$, and that $s_k-s_\ell\notin\frac12\Z$ for $k\ne \ell$. If \eqref{eq:entanglement} holds on a nonempty open set $\cO$, then $u_k=0$ on $\R^n$ for every $k$.
	\end{proposition}

	\begin{proposition}\label{prop:even-logarithmic}
		Let $n\ge2$ be even, and let $A$ satisfy \eqref{eq:ellipticity} and \eqref{eq:compact-perturbation}. Suppose that $u_k\in\mathscr D_{\log,p}$ for every $k=1,\ldots,N$, where $0\le p\le1$, and that $s_k-s_\ell\notin\Z$ for $k\ne\ell$. If \eqref{eq:entanglement} holds on a nonempty open set $\cO$, then $u_k=0$ on $\R^n$ for every $k$.
	\end{proposition}
	
	Both propositions include the endpoint weight $e^{-a\sqrt r/\log(e+r)}$.

	\subsection{Counterexamples and sharpness}
	
	The next result proves that the additional odd-dimensional order restriction is necessary, even for smooth compactly supported functions. It also shows in every dimension that the decay exponent $1/2$ in Theorem~\ref{thm:arbitrary} cannot be decreased.
	
	\begin{theorem}\label{thm:sharpness}
		Let $n\ge2$, and let $A$ satisfy \eqref{eq:ellipticity} and \eqref{eq:compact-perturbation}. Fix positive nonintegers $s<t$ and $R>0$ such that $\supp(A-I)\Subset B_R$.
		\begin{enumerate}[\rm(i)]
			\item\label{item_1_thm_sharp} If $n$ is odd and $t-s\in\frac12\N$, there are nonzero real functions $u,v\in C_c^\infty(\R^n\setminus\overline{B_R})$ satisfying
			\begin{equation}\label{eq:counterexample}
				(L^su+L^tv)|_{B_R}=0.
			\end{equation}
			\item\label{item_2_thm_sharp} For every $0<\gamma<1/2$, there are nonzero real functions $u,v\in\cS(\R^n)$, both zero on $B_R$, satisfying \eqref{eq:counterexample} and
			\begin{equation}\label{eq:counterexample-decay}
				|\partial^\alpha u(x)|+|\partial^\alpha v(x)|
				\le C_\alpha e^{-c|x|^\gamma}\quad\text{for all }x\in\R^n
			\end{equation}
			for every multi-index $\alpha$, with a constant $c>0$ independent of $\alpha$.
		\end{enumerate}
	\end{theorem}
	
	The compactly supported examples in part~\ref{item_1_thm_sharp} show that the additional half-integer restriction in odd dimensions is necessary, answering the question raised in \cite[Remark 1.3]{FL24}. Together with the algebraic cancellation for integer differences, these examples prove that the order conditions in Theorem~\ref{thm:arbitrary} are sharp.
	
	Part~\ref{item_2_thm_sharp} gives counterexamples for every $0<\gamma<1/2$ and every pair of positive noninteger orders $s<t$, including pairs satisfying those order conditions. Combined with Theorem~\ref{thm:arbitrary}, this proves that $\gamma=1/2$ is the optimal exponent within the scale $e^{-a|x|^\gamma}$ in every dimension. In particular, the assumption $u_k\in\cS(\R^n)$ alone does not suffice for entanglement.
	
	The logarithmic decay conditions in Propositions~\ref{prop:odd-logarithmic} and~\ref{prop:even-logarithmic} are also sharp.
	
	\begin{proposition}\label{prop:logarithmic-sharpness}
		Let $n\ge2$, let $A$ satisfy \eqref{eq:ellipticity} and \eqref{eq:compact-perturbation}, and let $R>0$ satisfy $\supp(A-I)\Subset B_R$. For every $p>1$ and every pair of positive nonintegers $s<t$, there are nonzero real functions $u,v\in\cS(\R^n)$, both zero on $B_R$, satisfying \eqref{eq:counterexample} and
		\begin{equation}\label{eq:logarithmic-counterexample-decay}
			|\partial^\alpha u(x)|+|\partial^\alpha v(x)|\le C_\alpha e^{-c w_p(|x|)},\quad x\in\R^n,
		\end{equation}
		for every multi-index $\alpha$, with $c>0$ independent of $\alpha$.
	\end{proposition}
	
	The functions in this proposition belong to $\mathscr D_{\log,p}$. Choosing the orders to satisfy the restrictions in Theorem~\ref{thm:arbitrary} shows that the range $0\le p\le1$ in Propositions~\ref{prop:odd-logarithmic} and~\ref{prop:even-logarithmic} is optimal in every dimension.
	
	The next result concerns the coefficients at infinity. We construct a sequence of coefficient fields converging to the identity in every fixed unweighted $C^k$ norm for which entanglement fails, even for smooth compactly supported inputs and orders allowed by the even-dimensional theorem. We use the notation $\Nzero=\N\cup\{0\}$.
	
	\begin{theorem}\label{thm:general-counterexample}
		There are smooth real symmetric matrices $A_\ell$ on $\R^2$, uniformly elliptic with constants independent of $\ell$, such that $A_\ell=I$ near the origin and $\norm{A_\ell-I}_{C^k(\R^2)}\to0$ for every fixed $k\in\Nzero$, with the following property. For $L_\ell=-\Div(A_\ell\nabla)$, there are nonzero real functions $u_\ell,v_\ell\in C_c^\infty(\R^2\setminus\overline{B_1})$ satisfying
		\begin{equation}\label{eq:general-counterexample}
			(L_\ell^{1/4}u_\ell+L_\ell^{3/4}v_\ell)|_{B_1}=0.
		\end{equation}
	\end{theorem}
	
	The matrices in this theorem are not equal to $I$ at infinity. The construction uses an angular mode for which the radial wave equation admits an outgoing solution. Its support moves away from $B_1$ for all positive times, and its Laplace transform gives the fractional relation. Outgoing radial waves also give the compactly supported examples in Theorem~\ref{thm:sharpness}\ref{item_1_thm_sharp}.
	
	\subsection{Unique continuation and the resolvent characterization}
	
	For a single $L^2$ function, finite sums of fractional powers satisfy unique continuation for general smooth uniformly elliptic coefficients. The orders may differ by integers, and the nonzero coefficients may be complex.
	
	\begin{theorem}\label{thm:single}
		Let $A$ satisfy \eqref{eq:ellipticity}, let $s_1,\ldots,s_N$ be distinct positive nonintegers, and let $b_k\in\C\setminus\{0\}$. If $u\in L^2(\R^n)$ satisfies
		\begin{equation}\label{eq:single}
			u|_{\cO}=0, \quad \sum_{k=1}^N b_kL^{s_k}u\Big|_{\cO}=0
		\end{equation}
		on a nonempty open set $\cO$, then $u=0$ on $\R^n$.
	\end{theorem}
	
	The theorem also holds in the negative spectral Sobolev spaces associated with $L$; see Corollary~\ref{cor:single-distributions}. A polynomial in $L$ may be added to the sum, since it vanishes locally wherever $u$ does. Related unique continuation results for fractional polyharmonic equations, with additional assumptions on the function in the interior or exterior of a bounded domain, were obtained in \cite{LLL25}.
	
	The proofs are based on a characterization in terms of resolvents. For $u_1,\ldots,u_N\in L^2(\R^n)$ vanishing on $\cO$, Proposition~\ref{prop:resolvent} shows that the fractional relation in \eqref{eq:entanglement} is equivalent to
	\begin{equation}\label{eq:intro-resolvent}
		\sum_{k=1}^N b_k\sin(\pi s_k)\lambda^{s_k}(L+\lambda)^{-1}u_k\Big|_{\cO}=0\quad\text{for every }\lambda>0.
	\end{equation}
	Applying integer powers of $L$ to the fractional relation and using a local resolvent formula gives a scalar function of $\lambda$ whose integer moments all vanish. An off-diagonal resolvent estimate makes this moment problem determinate and gives \eqref{eq:intro-resolvent}. The proof uses no explicit heat kernel formula and applies under the general coefficient assumptions in \eqref{eq:ellipticity}. When all the inputs coincide, the sum factors into a nonzero scalar function of $\lambda$ and a single resolvent. This gives Theorem~\ref{thm:single}.
	
	Exterior observation permits a reduction to a bounded domain. Local unique continuation extends the vanishing of the resolvent combination in \eqref{eq:intro-resolvent} to the unbounded component of the complement of the supports. This combination is then compactly supported and can be represented using a Dirichlet resolvent on a large ball. The Dirichlet resolvent is holomorphic in $\lambda$ at zero, so noninteger differences between the orders separate the terms. The residues at its poles recover the spectral projections of each input, and unique continuation for the corresponding eigenfunctions proves Theorem~\ref{thm:exterior}. This reduction explains why no assumption on the coefficients at infinity is needed for exterior observation.
	
	For arbitrary observation, we use the resolvent near zero on the whole space. Under \eqref{eq:compact-perturbation}, an equation for the gradient on a ball containing $\supp(A-I)$ separates the coefficient perturbation from the free resolvent. Uniform ellipticity gives an energy estimate that makes the operator in this equation invertible at zero. In odd dimensions, the localized free resolvent is holomorphic in the parameter $\kappa=\sqrt\lambda$ near zero. In even dimensions, its expansion contains integer powers of $\lambda$ and logarithmic terms. These different expansions account for the two order restrictions in Theorem~\ref{thm:arbitrary}.
	
	The order restrictions force all coefficients in the corresponding expansions of the resolvent pairings to vanish. For inputs without compact support, the passage from vanishing coefficients to vanishing pairings requires estimates on the remainders. Root-exponential decay gives Gevrey bounds of order two in odd dimensions and quantitative remainders for finite power--logarithm expansions in even dimensions. The resulting boundary estimates force the holomorphic pairings to vanish. The logarithmic improvements use finer moment estimates and quasianalyticity for holomorphic functions with bounds on all derivatives. In even dimensions, differentiated remainders give these bounds after the power--logarithm coefficients vanish. The argument requires only finite expansions for the noncompact input.
	
	The converse implication in \eqref{eq:intro-resolvent} provides a way to construct counterexamples. In dimension three, the radial resolvent on a ball disjoint from the support is a fixed spatial factor times a one-dimensional Laplace transform. A causal fractional derivative supplies the power of $\kappa$ needed to cancel two terms. Integer derivatives preserve compact support and give the half-integer examples. For noninteger derivatives, vanishing moments remove the algebraic tails that would otherwise prevent the required decay. The logarithmic examples use zero-moment profiles obtained from the boundary values of a holomorphic function on $\C\setminus[0,\infty)$. Lemmas~\ref{lem:fractional-tail} and~\ref{lem:logarithmic-fractional-tail} give the required decay after causal fractional differentiation. Abel inversion gives the corresponding construction in dimension two. The radial differential identity in Lemma~\ref{lem:dimension} gives higher dimensions, and Lemma~\ref{lem:transfer} transfers the cancellations to compact coefficient perturbations. The counterexamples in Theorem~\ref{thm:general-counterexample} arise instead from outgoing waves for coefficients that are not constant at infinity.

	\para{Comparison with other settings}
	The dimension-independent order conditions in \cite{FKU24,Lin26,LLY25} can be understood through the large time behavior of the semigroup profiles used in their proofs. On a closed manifold, applying the Laplacian removes the constant mode, and the spectral gap gives exponential decay of the resulting heat semigroup profiles; see \cite[Remark~1.9 and Eq.~(2.5)]{FKU24}. Together with small time off-diagonal decay, this yields entire Mellin transforms. The separation argument then uses Gamma factors whose pole sequences have integer spacing and are disjoint when the orders differ by nonintegers; see \cite[Proposition~3.1 and Lemma~3.5]{FKU24}.
	
	Hyperbolic space provides a noncompact example with the same dimension-independent order condition. The continuous spectrum of $-\Delta_{\mathbb H^n}$ starts at $(n-1)^2/4>0$ for $n\ge2$, and this positive lower bound gives the factor $e^{-(n-1)^2t/4}$ in the heat-semigroup estimate; see \cite[Section~2.1.1 and Eq.~(3.13)]{Lin26}. Under the decay hypotheses of \cite[Theorem~1.1]{Lin26}, that estimate supplies the bounds at both ends required by the separation criterion in \cite[Proposition~3.3]{Lin26}. This example shows that noncompactness and continuous spectrum alone do not account for the additional odd-dimensional restriction in the Euclidean setting: the hyperbolic spectrum is separated from zero.
	
	For the parabolic operator $H_g=\partial_t-\Delta_g$, the identity $(e^{-\tau H_g}u)(x,t)=(e^{\tau\Delta_g}[u(\cdot,t-\tau)])(x)$ shows that large values of the auxiliary parameter $\tau$ sample the distant past of the input; see \cite[Eq.~(2.3)]{LLY25}. The backward-time decay assumption in \cite[Eq.~(1.4)]{LLY25} gives exponential decay of this profile as $\tau\to\infty$ for fixed $t$. Under the zero history condition in \cite[Eq.~(4.1)]{LLY25}, the profile vanishes for sufficiently large $\tau$. These temporal assumptions supply the large time control needed for the same separation argument, without a distinction between odd and even spatial dimensions; see also \cite[Remark~3.2]{LLY25}.
	
	In the present elliptic setting on $\R^n$, compact support of the coefficient perturbation does not produce a spectral gap: zero remains a continuous spectrum threshold. The dimension-dependent low energy structure of the free resolvent, described by \eqref{eq:odd-kernel} and \eqref{eq:even-kernel}, enters the separation of fractional orders in \eqref{eq:intro-resolvent}. The compactly supported counterexamples in Theorem~\ref{thm:sharpness}\ref{item_1_thm_sharp} show that the additional half-integer restriction in odd dimensions is necessary for arbitrary observation. Exterior observation permits the bounded-domain reduction described above, where the Dirichlet resolvent is holomorphic in $\lambda$ at zero, yielding the dimension-independent order condition in Theorem~\ref{thm:exterior}.

	\subsection{Applications to inverse problems}\label{subsec:inverse-introduction}
	
	Theorems~\ref{thm:exterior} and~\ref{thm:single} give inverse results for a known operator $L$ satisfying \eqref{eq:ellipticity}, with no assumption on its coefficients at infinity. We first recover a potential and a local anisotropic principal term, extending the applications in \cite{FL24}. For operators consisting of noninteger powers of $L$, we then allow the orders and their constant coefficients to be unknown. Under the solvability assumptions below, one nonzero exterior datum and its measurement determine the number, orders, and coefficients of these terms, together with a continuous bounded potential. For this simultaneous recovery, the measurement is taken on a nonempty open subset of the unbounded connected component of the complement of the union of the closure of the domain and the support of the datum.
	
	All potentials in our inverse problems are real-valued. Let $\Omega\subset\R^n$ be a bounded Lipschitz domain, and write $\Omega_e:=\R^n\setminus\overline\Omega$. We first consider
	\begin{equation}\label{eq:inverse-P}
		P=cL^m+\sum_{k=1}^N b_kL^{s_k},
		\quad c\in\R,\quad m\in\N,\quad b_k\in\R\setminus\{0\},
	\end{equation}
	where $0<s_1<\cdots<s_N$ are nonintegers. The local term is omitted when $c=0$. For $q\in L^\infty(\Omega)$, let $u_f$ solve
	\begin{equation}\label{eq:inverse-problem}
		\begin{cases}
			(P+q)u_f=0&\text{in }\Omega,\\
			u_f=f&\text{in }\Omega_e.
		\end{cases}
	\end{equation}
	We use the weak exterior Dirichlet formulation in Section~\ref{sec:inverse-problems} and assume that its homogeneous problem has only the zero solution. When $c\ge0$ and all $b_k$ are positive, this holds in particular if $q\ge0$. For $f\in C_c^\infty(\Omega_e)$, define the exterior Dirichlet-to-Neumann (DN) map by
	\begin{equation}\label{eq:inverse-map}
		\Lambda_{P,q}f=(Pu_f)|_{\Omega_e},
	\end{equation}
	with values in distributions on $\Omega_e$. Maps associated with different fractional orders can be compared by evaluating them on $C_c^\infty(\Omega_e)$.
	
	\begin{theorem}[Recovery of a potential]\label{thm:inverse-potential}
		Let $P$ be known and given by \eqref{eq:inverse-P}. Let $q_1,q_2\in L^\infty(\Omega)$ satisfy the preceding solvability assumption, and let $W_1,W_2\subset\Omega_e$ be nonempty open sets. If
		\begin{equation}\label{eq:inverse-potential-data}
			\left.\Lambda_{P,q_1}f\right|_{W_2}
			=\left.\Lambda_{P,q_2}f\right|_{W_2}
			\quad\text{for every }f\in C_c^\infty(W_1),
		\end{equation}
		then $q_1=q_2$ almost everywhere in $\Omega$.
	\end{theorem}
	
	No restriction on differences of the orders is needed. The proof uses unique continuation for the whole sum from Theorem~\ref{thm:single} and the Runge approximation in Lemma~\ref{lem:inverse-Runge}; see \cite{GSU20,GLX,FL24} for the corresponding fractional inverse problems.
	
	For the anisotropic application, let $\mathcal N=\sum_{k=1}^N b_kL^{s_k}$, where $0<s_1<\cdots<s_N$ are nonintegers and $b_k\in\R\setminus\{0\}$. We consider
	\begin{equation}\label{eq:inverse-anisotropic-problem}
		\begin{cases}
			-\Div(C\nabla u_f)+\mathcal Nu_f+qu_f=0&\text{in }\Omega,\\
			u_f=f&\text{in }\Omega_e,
		\end{cases}
	\end{equation}
	where the local matrix $C$ is unknown. The matrix $A$ defining the fractional powers of $L$ remains known. We take $C\in C^\infty(\overline\Omega;\R^{n\times n})$ real symmetric and uniformly positive definite, and $q\in L^\infty(\Omega)$. Assuming again that the homogeneous weak Dirichlet problem has only the zero solution, define $\Lambda_{C,q}^{\mathcal N}f=(\mathcal Nu_f)|_{\Omega_e}$.
	
	\begin{theorem}[Recovery of a local anisotropic operator]\label{thm:inverse-anisotropic}
		Let $\mathcal N$ be known as above, and let $C_j,q_j$, $j=1,2$, satisfy the preceding assumptions. For any nonempty open sets $W_1,W_2\subset\Omega_e$, the equality
		\begin{equation}\label{eq:inverse-anisotropic-data}
			\Lambda_{C_1,q_1}^{\mathcal N}f\big|_{W_2}=\Lambda_{C_2,q_2}^{\mathcal N}f\big|_{W_2}
			\quad\text{for every }f\in C_c^\infty(W_1)
		\end{equation}
		implies $C_1=C_2$ in $\Omega$ and $q_1=q_2$ almost everywhere in $\Omega$.
	\end{theorem}
	
	The known nonlocal term allows recovery of the local matrix $C$ in the fixed coordinates of $\Omega$. We next allow the fractional orders and their constant coefficients to be unknown. For $j=1,2$, write
	\begin{equation}\label{eq:inverse-unknown-P}
		P_j=\sum_{k=1}^{N_j}b_{j,k}L^{s_{j,k}},
		\quad N_j\in\N,\quad 0<s_{j,1}<\cdots<s_{j,N_j},\quad
		s_{j,k}\notin\N,\quad b_{j,k}\in\R\setminus\{0\}.
	\end{equation}
	Let $q_j\in L^\infty(\Omega)$. For each $j=1,2$, assume that the homogeneous weak exterior Dirichlet problem associated with $P_j+q_j$ has only the zero solution, as in \eqref{eq:inverse-nonresonance}. Integer differences between fractional orders are allowed, both within each operator and between the two operators. For a fixed nonzero $f\in C_c^\infty(\Omega_e)$, let $F_f$ denote the unbounded connected component of $\R^n\setminus(\overline\Omega\cup\supp f)$.
	
	\begin{theorem}[Recovery of fractional orders and coefficients]\label{thm:inverse-orders}
		Let $P_j,q_j$ be as above, and let $W\subset F_f$ be a nonempty open set. If
		\begin{equation}\label{eq:inverse-orders-data}
			\left.\Lambda_{P_1,q_1}f\right|_W
			=\left.\Lambda_{P_2,q_2}f\right|_W,
		\end{equation}
		then $N_1=N_2$ and, in increasing order,
		\[
		s_{1,k}=s_{2,k},\quad b_{1,k}=b_{2,k}
		\quad\text{for every }k=1,\ldots,N_1.
		\]
		The corresponding solutions also coincide on $\R^n$.
	\end{theorem}
	
	This conclusion uses one exterior datum, with no restriction on integer differences between the orders. In the proof, orders with the same fractional part are first grouped into polynomials in $L$. Theorem~\ref{thm:exterior} separates the distinct fractional parts. The datum $f$ determines each polynomial, because a nonzero polynomial in $L$ cannot annihilate a nonzero compactly supported smooth function. The condition on $W$ holds, for example, if $W$ lies outside a ball containing $\overline\Omega\cup\supp f$.
	
	\begin{corollary}[Simultaneous recovery from one exterior datum]\label{cor:inverse-simultaneous}
		Under the assumptions of Theorem~\ref{thm:inverse-orders}, suppose in addition that $q_j\in C(\Omega)\cap L^\infty(\Omega)$ for $j=1,2$. Then the equality \eqref{eq:inverse-orders-data} also implies $q_1=q_2$ in $\Omega$.
	\end{corollary}
	
	For known $P$ in \eqref{eq:inverse-P}, a continuous bounded potential is likewise determined by one nonzero $f\in C_c^\infty(\Omega_e)$ and its image on any nonempty open $W\subset\Omega_e$. Single-measurement recovery for the fractional Schr\"odinger equation was studied in \cite{GRSU20}. Here continuity of the potential allows the proof to use only unique continuation from open sets.

	\subsection{Earlier literature}\label{subsec:earlier-literature}
	
	The extension formulation for the fractional Laplacian in \cite{CS07} and its generalization to nonnegative self-adjoint operators in \cite{ST10} realize powers of order $0<s<1$ through a boundary value problem in one additional variable. Extension methods have been used to prove strong unique continuation for fractional elliptic equations with potentials, variable coefficients, and higher orders; see \cite{FF14,Rul15,Yu17,GFR19}. For $0<s<1$, local approximation by compactly supported functions that are $s$-harmonic on a ball was established in \cite{DSV17}. The unique continuation and Runge approximation results in \cite{GSU20,GLX,CMR20} are central to exterior inverse problems.
	
	The monograph \cite{LL25} gives a systematic treatment of inverse problems for integro-differential operators. For a known fractional operator, exterior measurements determine drift terms and potentials \cite{CLR20}, embedded obstacles and surrounding potentials \cite{CLL19}, and linear local perturbations of higher-order fractional equations \cite{CMRU22}. Recovery of a potential from a single measurement was proved in \cite{GRSU20}; see also \cite{HL26} for recovery from measurements at a single exterior point. Quantitative approximation, stability, exponential instability, and monotonicity-based reconstruction were studied in \cite{RS20,RS18,HL19,HL20}. Unique continuation for fractional polyharmonic operators under additional assumptions was considered in \cite{LLL25}. The entanglement principle in \cite{FL24} permits recovery of local anisotropic principal terms in fractional polyharmonic equations.
	
	Inverse problems for the operator defining the fractional power have also been studied. Local source-to-solution data determine closed Riemannian manifolds up to isometry under the assumptions in \cite{Fei24,FGKU25}, and an extension approach was developed in \cite{Rul25}. Simultaneous recovery of a closed Riemannian manifold and a potential from local fractional Cauchy data was established in \cite{FKU26}. For divergence form operators in dimensions at least three, partial fractional exterior data determine the corresponding local boundary data under the assumptions in \cite{CGRU26}. For smooth metrics equal to the Euclidean metric outside a compact set, fractional exterior data determine the metric up to a diffeomorphism fixing the exterior \cite{FGKRSU25}. Potential recovery for space-time fractional parabolic equations was proved in \cite{LLR20,BS24}, and leading-coefficient inverse problems were reduced to local parabolic problems in \cite{LLU22}.

	\para{Organization of the paper}
	We begin with the resolvent characterization and the regularization of distributional inputs in Section~\ref{sec:preliminaries}, followed by the proof of exterior entanglement in Section~\ref{sec:exterior}. For arbitrary observation, we need estimates for the resolvent near zero. These are developed in Section~\ref{sec:low-energy} and used in Section~\ref{sec:arbitrary} to prove the root-exponential and logarithmic decay results. The counterexamples follow in Section~\ref{sec:sharpness}. We then examine the role of the coefficients at infinity in Section~\ref{sec:general-counterexample}, before turning to inverse problems in Section~\ref{sec:inverse-problems}.

	\section{Fractional powers and resolvent identities}\label{sec:preliminaries}
	
	\subsection{Functional calculus}
	
	In the estimates below, $a_1$ denotes the upper ellipticity bound in \eqref{eq:ellipticity}.
	
	All inner products and duality pairings are linear in the first argument. The sesquilinear form 
	\begin{equation}\label{eq:def_of_mathfrak_a}
		\mathfrak a(u,v)=\int_{\R^n}A\nabla u\cdot\overline{\nabla v}\,dx, \quad \text{for any }u,v\in H^1(\R^n)
	\end{equation}
	defines the nonnegative self-adjoint operator $L$. In particular, $\Dom(L^{1/2})=H^1(\R^n)$, with equivalent norms after adding the $L^2$ norm. For $r\ge0$, we set $\HL^r=\Dom((I+L)^{r/2})$ and equip it with the norm $\norm{(I+L)^{r/2}u}_2$. We denote by $\HL^{-r}$ the space of continuous conjugate linear functionals on $\HL^r$, identifying each $f\in L^2(\R^n)$ with the functional $v\mapsto(f,v)_{L^2}$. The spectral calculus gives a bounded map $L^s:\HL^r\to\HL^{r-2s}$ for every $r\in\R$ and $s\ge0$.
	
	Every $\phi\in C_c^\infty(\R^n)$ belongs to $\Dom(L^m)$ for every $m\in\Nzero$, since $L^m\phi$ is again smooth and compactly supported. For $u\in L^2(\R^n)$, the distribution $L^su$ is defined by $\ip{L^su}{\phi}=(u,L^s\phi)_{L^2}$. This agrees with the spectral definition whenever $u\in\Dom(L^s)$. For every compact set $K$ and every $k\in\Nzero$, iterated interior elliptic estimates and Sobolev embedding give
	\begin{equation}\label{eq:local-spectral-regularity}
		\norm{v}_{C^k(K)}\le C_{K,k,m}\norm{(I+L)^mv}_2,\quad v\in\Dom(L^m),\quad 2m>k+n/2.
	\end{equation}
	These estimates use only bounds on the coefficients and their derivatives near $K$. Every compactly supported distribution belongs to $\HL^{-2m}$ for sufficiently large $m$. The density of $C_c^\infty(\R^n)$ in the positive spectral spaces, proved in Lemma~\ref{lem:smoothing}, allows us to regard elements of the corresponding negative spectral spaces as distributions. For $u\in\cE(\R^n)\cap\HL^{-2m}$, we define $L^su\in\HL^{-2m-2s}$ by the spectral calculus.
	
	Under \eqref{eq:compact-perturbation}, the coefficients and all their derivatives are bounded. Global elliptic estimates, interpolation, and duality give
	\begin{equation}\label{eq:spectral-sobolev-equivalence}
		\HL^r=H^r(\R^n),\quad r\in\R,
	\end{equation}
	with equivalent norms. In particular, the fractional powers are defined on each class $\mathscr D_\gamma$ and $\mathscr D_{\log,p}$. For general coefficients satisfying \eqref{eq:ellipticity}, we work in the spectral spaces and use the local estimate \eqref{eq:local-spectral-regularity}.
	
	We use unique continuation for the (local) equations $(L+\lambda)v=0$ and $(L-\mu)v=0$ on connected open sets; see \cite{Aronszajn57}. We also use the fractional unique continuation theorem \cite[Theorem 1.2]{GLX}: if $0<\alpha<1$, $v\in H^\alpha(\R^n)$, and $v=L^\alpha v=0$ on a nonempty open set, then $v=0$ on $\R^n$.
	
	\subsection{The local resolvent formula}
	
	Throughout this subsection, we write $T_t=e^{-tL}$ for $t>0$ and $R_z=(L+z)^{-1}$ for $\Rea z>0$, in particular $R_\lambda=(L+\lambda)^{-1}$ for $\lambda>0$. The spectral theorem gives $\norm{R_z}_{L^2\to L^2}\le(\Rea z)^{-1}$ and holomorphic dependence on $z$ in this half-plane. The notation $(\cdot,\cdot)_{L^2}$ and $\norm{\cdot}_2$ refers to the inner product and norm on $L^2(\R^n)$.
	
	All distances are Euclidean, and $d>0$ denotes a separation parameter satisfying $\dist(\supp u,\supp\phi)\ge d$. For the local formula, when $u\ne0$ vanishes on an open set $\cO$ and $0\ne\phi\in C_c^\infty(\cO)$, we choose $d=\dist(\supp\phi,\R^n\setminus\cO)>0$.
	
	For a nonintegral $s>0$, we write $s=m+\alpha$, where $m=\lfloor s\rfloor\in\Nzero$ and $0<\alpha<1$.
	
	\begin{lemma}\label{lem:off-diagonal}
		Suppose that $u,\phi\in L^2(\R^n)$ and $\dist(\supp u,\supp\phi)\ge d>0$. There are constants $c,C>0$, depending only on the ellipticity constants, such that, for every $\lambda>0$,
		\begin{equation}\label{eq:off-diagonal}
			|(R_\lambda u,\phi)_{L^2}|\le C\lambda^{-1}e^{-cd\sqrt\lambda}\norm{u}_2\norm{\phi}_2.
		\end{equation}
	\end{lemma}
	
	\begin{proof}
		We may assume that $u\ne0$ and $\phi\ne0$. Choose $\rho(x)=\min\{\dist(x,\supp u),d\}$. Then $\rho$ is bounded and Lipschitz, $|\nabla\rho|\le1$ almost everywhere, and $\rho=0$ on $\supp u$ and $\rho=d$ on $\supp\phi$. Fix a constant $\eta>0$, and set $v(t)=T_tu$ and $w(t)=e^{\eta\rho}v(t)$. For $t>0$, we have $v(t)\in\Dom(L)$ and $e^{2\eta\rho}v(t)\in H^1(\R^n)$. By \eqref{eq:def_of_mathfrak_a} and the real symmetry of $A$,
		\begin{align*}
			\frac12\frac{d}{dt}\norm{w(t)}_2^2
			&=-\Rea\mathfrak a(v(t),e^{2\eta\rho}v(t))\\
			&=-\mathfrak a(w(t),w(t))+\eta^2\int_{\R^n}(A\nabla\rho\cdot\nabla\rho)|w(t)|^2\,dx\\
			&\le a_1\eta^2\norm{w(t)}_2^2.
		\end{align*}
		Here the cross terms have zero real part. Gronwall's inequality and strong continuity at $t=0$ give $\norm{w(t)}_2^2\le e^{2a_1\eta^2t}\norm{u}_2^2$, since $e^{\eta\rho}u=u$. Taking square roots and applying the Cauchy--Schwarz inequality, we obtain
		\begin{align*}
			|(T_tu,\phi)_{L^2}|=|(w(t),e^{-\eta\rho}\phi)_{L^2}|\le e^{a_1\eta^2t-\eta d}\norm{u}_2\norm{\phi}_2.
		\end{align*}
		For each fixed $t>0$, choosing $\eta=d/(2a_1t)$ gives the Davies--Gaffney estimate
		\begin{equation}\label{eq:DG}
			|(T_tu,\phi)_{L^2}|\le e^{-d^2/(4a_1t)}\norm{u}_2\norm{\phi}_2,\quad t>0.
		\end{equation}
		Its abstract formulation appears in \cite[Section~3, equation~(3.2)]{CS08}, with Euclidean distance scaled by $a_1^{-1/2}$.
		
		The spectral theorem gives $R_\lambda u=\int_0^\infty e^{-\lambda t}T_tu\,dt$, with absolute convergence in $L^2(\R^n)$ because $\norm{T_tu}_2\le\norm{u}_2$. Using \eqref{eq:DG} and the arithmetic--geometric mean inequality $\lambda t/2+d^2/(4a_1t)\ge d\sqrt{\lambda/(2a_1)}$, we obtain
		\begin{align*}
			|(R_\lambda u,\phi)_{L^2}|\le\norm{u}_2\norm{\phi}_2\int_0^\infty e^{-\lambda t-d^2/(4a_1t)}\,dt\le e^{-d\sqrt{\lambda/(2a_1)}}\norm{u}_2\norm{\phi}_2\int_0^\infty e^{-\lambda t/2}\,dt,
		\end{align*}
		which proves \eqref{eq:off-diagonal} with $C=2$ and $c=(2a_1)^{-1/2}$.
	\end{proof}
	
	\begin{lemma}\label{lem:local-power}
		Let $u\in L^2(\R^n)$ vanish on an open set $\cO$, and let $s>0$ be nonintegral. For every $\phi\in C_c^\infty(\cO)$,
		\begin{equation}\label{eq:local-power}
			\ip{L^su}{\phi}=-\frac{\sin(\pi s)}{\pi}\int_0^\infty\lambda^s(R_\lambda u,\phi)_{L^2}\,d\lambda,
		\end{equation}
		where the integral converges absolutely.
	\end{lemma}
	
	\begin{proof}
		Fix $\phi\in C_c^\infty(\cO)$, and assume that $u\ne0$ and $\phi\ne0$. With $s=m+\alpha$ as above, the Balakrishnan formula \cite[p.~420, formula~(2.1)]{Bal60}, applied with the operator there equal to $-L$, gives
		\begin{equation}\label{eq:Bala_formula_Dom}
			L^\alpha f=\frac{\sin(\pi\alpha)}{\pi}\int_0^\infty\lambda^{\alpha-1}R_\lambda Lf\,d\lambda,\quad f\in\Dom(L).
		\end{equation}
		The spectral multiplier $\mu/(\mu+\lambda)$ is bounded by $1$ for $\mu\ge0$, and $\norm{R_\lambda}_{L^2\to L^2}\le\lambda^{-1}$. These bounds give $\norm{R_\lambda Lf}_2\le\norm{f}_2$ and $\norm{R_\lambda Lf}_2\le\lambda^{-1}\norm{Lf}_2$. The integrand is bounded by $\lambda^{\alpha-1}\norm{f}_2$ near zero and by $\lambda^{\alpha-2}\norm{Lf}_2$ at infinity. Both bounds are integrable for $0<\alpha<1$, so the integral converges absolutely in $L^2(\R^n)$.
		
		Locality and smoothness of $L$ give $L^k\phi\in C_c^\infty(\cO)$ and $\supp(L^k\phi)\subset\supp\phi$ for every $k\in\Nzero$. In particular, $\phi\in\Dom(L^{m+1})$, so the formula applies to $f=L^m\phi$, with $L^\alpha L^m\phi=L^s\phi$. Iterating $R_\lambda L\chi=\chi-\lambda R_\lambda\chi$ for $\chi\in\Dom(L)$ gives
		\begin{equation}\label{eq:resolvent-polynomial}
			R_\lambda L^{m+1}\phi=\sum_{\ell=0}^m(-\lambda)^\ell L^{m-\ell}\phi+(-\lambda)^{m+1}R_\lambda\phi.
		\end{equation}
		Every term in the finite sum pairs to zero with $u$, since $u=0$ on $\cO$. Pairing \eqref{eq:resolvent-polynomial} first and using the self-adjointness of $R_\lambda$, we obtain, for $0<\eps<M<\infty$,
		\begin{equation}\label{eq:truncated-local-power}
			\int_\eps^M\lambda^{\alpha-1}(u,R_\lambda L^{m+1}\phi)_{L^2}\,d\lambda=(-1)^{m+1}\int_\eps^M\lambda^s(R_\lambda u,\phi)_{L^2}\,d\lambda.
		\end{equation}
		
		Since $\supp u\subset\R^n\setminus\cO$, our choice of $d$ gives $\dist(\supp u,\supp\phi)\ge d>0$. Multiplying \eqref{eq:off-diagonal} by $\lambda^s$ yields
		\begin{equation}\label{eq:local-power-integrand-bound}
			\lambda^s|(R_\lambda u,\phi)_{L^2}|\le C\lambda^{s-1}e^{-cd\sqrt\lambda}\norm{u}_2\norm{\phi}_2,\quad \lambda>0.
		\end{equation}
		This bound is integrable near zero because $s>0$. At infinity, the substitution $r=\sqrt\lambda$ reduces its integral to a constant multiple of $\int_1^\infty r^{2s-1}e^{-cdr}\,dr<\infty$.
		
		Multiply \eqref{eq:truncated-local-power} by $\sin(\pi\alpha)/\pi$. The Balakrishnan formula \eqref{eq:Bala_formula_Dom} applied to $f=L^m\phi$ gives
		\[
		\frac{\sin(\pi\alpha)}{\pi}\int_\eps^M\lambda^{\alpha-1}R_\lambda L^{m+1}\phi\,d\lambda\longrightarrow L^\alpha L^m\phi=L^s\phi
		\]
		in $L^2(\R^n)$ as $\eps\to0$ and $M\to\infty$. Pairing with the fixed function $u\in L^2(\R^n)$ shows that the left-hand side of the multiplied equality converges to $(u,L^s\phi)_{L^2}=\ip{L^su}{\phi}$. On the right-hand side, the integrable bound \eqref{eq:local-power-integrand-bound} permits passage from $\int_\eps^M$ to $\int_0^\infty$. The coefficient is $(-1)^{m+1}\sin(\pi\alpha)/\pi=-\sin(\pi s)/\pi$, which proves \eqref{eq:local-power}.
	\end{proof}
	
	\begin{lemma}\label{lem:moment}
		Let $H:(0,\infty)\to\C$ be measurable with $H\in L^1(0,1)$. Suppose that, for some constants $C,c>0$ and $q\ge0$,
		\begin{equation}\label{eq:assum_H}
			|H(\lambda)|\le C(1+\lambda)^q e^{-c\sqrt\lambda},\quad \lambda\ge1,
		\end{equation}
		and that
		\begin{equation}\label{eq:all-moments}
			\int_0^\infty\lambda^kH(\lambda)\,d\lambda=0\quad\text{for every }k\in\Nzero.
		\end{equation}
		Then $H=0$ almost everywhere.
	\end{lemma}
	
	\begin{proof}
		Set $g(r)=2rH(r^2)$ for $r>0$ and $c_0=c/2$. The substitution $\lambda=r^2$ and the assumption \eqref{eq:assum_H} on $H$ give
		\[
		\int_0^\infty e^{c_0r}|g(r)|\,dr=\int_0^\infty e^{c_0\sqrt\lambda}|H(\lambda)|\,d\lambda<\infty.
		\]
		Define $G(z)=\int_0^\infty g(r)\cos(zr)\,dr$ for $|\Ima z|<c_0$. For every $0<b<c_0$ and $k\in\Nzero$, the function $r^ke^{br}|g(r)|$ is integrable. This justifies differentiation under the integral on each smaller strip $|\Ima z|\le b$, so $G$ is holomorphic in $|\Ima z|<c_0$.
		
		For every $k\in\Nzero$, the substitution $\lambda=r^2$ and \eqref{eq:all-moments} give
		\[
		G^{(2k)}(0)=(-1)^k\int_0^\infty r^{2k}g(r)\,dr=(-1)^k\int_0^\infty\lambda^kH(\lambda)\,d\lambda=0.
		\]
		All odd derivatives at zero vanish because $G$ is even. Its Taylor series at zero is identically zero, and the identity theorem gives $G=0$ throughout the strip.
		
		Let $\widetilde g\in L^1(\R)$ be the even extension of $g$. With the Fourier transform convention $\widehat{\widetilde g}(\xi)=\int_\R e^{-\mathsf{i}\xi r}\widetilde g(r)\,dr$, we have $\widehat{\widetilde g}(\xi)=2G(\xi)=0$ for every $\xi\in\R$. Injectivity of the Fourier transform on $L^1(\R)$ gives $g=0$ almost everywhere. The substitution $\lambda=r^2$ then gives $H=0$ almost everywhere.
	\end{proof}
	
	\begin{proposition}[Resolvent characterization]\label{prop:resolvent}
		Let $u_k\in L^2(\R^n)$ vanish on a nonempty open set $\cO\subset\R^n$. For positive nonintegral orders $s_k$ and constants $b_k\in\C$, the fractional relation in \eqref{eq:entanglement} holds if and only if \eqref{eq:intro-resolvent} holds for every $\lambda>0$.
	\end{proposition}
	
	\begin{proof}
		The assertion is immediate if $\cO=\R^n$. Otherwise, fix $0\ne\phi\in C_c^\infty(\cO)$ and set
		\[
		H_\phi(\lambda)=\sum_{k=1}^N b_k\sin(\pi s_k)\lambda^{s_k}(R_\lambda u_k,\phi)_{L^2},\quad \lambda>0.
		\]
		Set $d=\dist(\supp\phi,\R^n\setminus\cO)>0$. Since $\supp u_k\subset\R^n\setminus\cO$, estimate~\eqref{eq:off-diagonal} applies to every nonzero $u_k$ and gives
		\[
		|H_\phi(\lambda)|\le C e^{-cd\sqrt\lambda}\norm{\phi}_2\sum_{k=1}^N |b_k\sin(\pi s_k)|\lambda^{s_k-1}\norm{u_k}_2,\quad \lambda>0.
		\]
		This bound is integrable near zero because every $s_k>0$, and it has the decay required in Lemma~\ref{lem:moment} at infinity.
		
		Suppose that \eqref{eq:entanglement} holds. For every $\ell\in\Nzero$, locality gives $L^{\ell}\phi\in C_c^\infty(\cO)$, and the spectral calculus gives $\ip{L^{s_k}u_k}{L^{\ell}\phi}=\ip{L^{s_k+\ell}u_k}{\phi}$. Testing the fractional relation against $L^{\ell}\phi$ and applying Lemma~\ref{lem:local-power}, we obtain
		\begin{align*}
			0=\sum_{k=1}^N b_k\ip{L^{s_k}u_k}{L^{\ell}\phi}=\sum_{k=1}^N b_k\ip{L^{s_k+\ell}u_k}{\phi}=-\frac{(-1)^{\ell}}{\pi}\int_0^\infty\lambda^{\ell}H_\phi(\lambda)\,d\lambda,
		\end{align*}
		where we used $\sin(\pi(s_k+\ell))=(-1)^{\ell}\sin(\pi s_k)$. Every integral here converges absolutely by Lemma~\ref{lem:local-power}. All integer moments of $H_\phi$ vanish, so Lemma~\ref{lem:moment} gives $H_\phi=0$ almost everywhere.
		
		The resolvent identity gives $\norm{R_\lambda-R_\mu}_{L^2\to L^2}\le|\lambda-\mu|/(\lambda\mu)$ for $\lambda,\mu>0$. In particular, $H_\phi$ is continuous, and its almost everywhere vanishing implies $H_\phi(\lambda)=0$ for every $\lambda>0$. Since $\phi$ was arbitrary, this proves \eqref{eq:intro-resolvent}.
		
		Conversely, suppose that \eqref{eq:intro-resolvent} holds. Pairing with $\phi$ gives $H_\phi(\lambda)=0$ for every $\lambda>0$. The absolute convergence in Lemma~\ref{lem:local-power} allows us to integrate each term, giving
		\[
		\sum_{k=1}^N b_k\ip{L^{s_k}u_k}{\phi}=-\frac1\pi\int_0^\infty H_\phi(\lambda)\,d\lambda=0,
		\]
		which proves the fractional relation on $\cO$.
	\end{proof}
	
	\begin{lemma}\label{lem:individual}
		Let $\cO\subset\R^n$ be a nonempty open set. If $u\in L^2(\R^n)$ and $(R_\lambda u)|_{\cO}=0$ for every $\lambda>0$, then $u=0$ on $\R^n$.
	\end{lemma}
	
	\begin{proof}
		Set $v=R_1u\in\Dom(L)\subset H^1(\R^n)$. By hypothesis, $v=0$ on $\cO$. The resolvent identity gives
		\[
		R_\lambda v=\frac{R_\lambda u-R_1u}{1-\lambda},\quad \lambda>0,\ \lambda\ne1.
		\]
		Both terms in the numerator vanish on $\cO$, so $(R_\lambda v)|_{\cO}=0$ for $\lambda\ne1$. Since $R_\lambda v\to R_1v$ in $L^2(\R^n)$ as $\lambda\to1$, we also have $R_1v=0$ on $\cO$. Lemma~\ref{lem:local-power}, applied with $s=1/2$, now gives
		\[
		\ip{L^{1/2}v}{\phi}=-\frac1\pi\int_0^\infty\lambda^{1/2}(R_\lambda v,\phi)_{L^2}\,d\lambda=0
		\]
		for every $\phi\in C_c^\infty(\cO)$.
		
		Since $v\in H^1(\R^n)\subset H^{1/2}(\R^n)$ and $v=L^{1/2}v=0$ on $\cO$, the fractional unique continuation theorem \cite[Theorem~1.2]{GLX} gives $v=0$ on $\R^n$. Finally, $u=(L+1)v=0$.
	\end{proof}
	
	\begin{proof}[Proof of Theorem~\ref{thm:single}]
		Relabel the orders so that $s_1<\cdots<s_N$, and set $p(\lambda)=\sum_{k=1}^N b_k\sin(\pi s_k)\lambda^{s_k}$. Proposition~\ref{prop:resolvent}, applied with $u_k=u$, gives
		\[
		p(\lambda)(R_\lambda u)|_{\cO}=0,\quad \lambda>0.
		\]
		Since
		\[
		\lim_{\lambda\to0}\lambda^{-s_1}p(\lambda)=b_1\sin(\pi s_1)\ne0,
		\]
		there exists $\lambda_0>0$ such that $p(\lambda)\ne0$ for $0<\lambda<\lambda_0$. It follows that $(R_\lambda u)|_{\cO}=0$ on this interval.
		
		For each $\phi\in C_c^\infty(\cO)$, the function $z\mapsto(R_z u,\phi)_{L^2}$ is holomorphic in $\Rea z>0$ and vanishes for $0<z<\lambda_0$. The identity theorem gives $(R_z u,\phi)_{L^2}=0$ for every $\Rea z>0$, in particular for every real $z=\lambda>0$. Since $\phi$ was arbitrary, $(R_\lambda u)|_{\cO}=0$ for every $\lambda>0$. Lemma~\ref{lem:individual} gives $u=0$.
	\end{proof}
	
	\subsection{Regularization with finite propagation}\label{sec:distributions}
	
	To pass from $L^2$ functions to distributions, we need a regularization that commutes with fractional powers of $L$ and preserves vanishing on a smaller open set. Throughout this subsection, $\eps>0$ denotes the regularization parameter, and we set $c_L=\sqrt{a_1}$. We write $\delta_y$ for the point mass at $y$ and $\partial^\beta\delta_y$ for its distributional derivatives. Choose a real even function $\zeta\in C_c^\infty((-1,1))$ with $\int_\R\zeta(t)\,dt=1$, and set $m(\sigma)=\int_\R\zeta(t)\cos(t\sigma)\,dt$ for $\sigma\in\R$. Define
	\begin{equation}\label{eq:spectral-smoothing}
		S_\eps=m(\eps\sqrt L)=\int_\R\zeta(t)\cos(\eps t\sqrt L)\,dt.
	\end{equation}
	For each $f\in L^2(\R^n)$, the integral in \eqref{eq:spectral-smoothing}, applied to $f$, converges absolutely in $L^2(\R^n)$ because $\norm{\cos(\eps t\sqrt L)f}_2\le\norm{f}_2$.
	
	\begin{lemma}\label{lem:smoothing}
		For $\eps>0$ and $r,M\in\R$, the operator $S_\eps$ maps $\HL^r$ continuously into $\HL^M$, commutes with $L^s$ for every $s>0$, and satisfies $S_\eps u\to u$ in $\HL^r$ as $\eps\to 0$. The operator $S_\eps$ has a kernel $k_\eps\in C^\infty(\R^n\times\R^n)$ such that $S_\eps f(x)=\int_{\R^n}k_\eps(x,y)f(y)\,dy$ for every $f\in C_c^\infty(\R^n)$. The kernel satisfies $k_\eps(x,y)=k_\eps(y,x)$ and $k_\eps(x,y)=0$ whenever $|x-y|>c_L\eps$. Moreover, $C_c^\infty(\R^n)$ is dense in $\HL^r$ for every $r\ge0$.
	\end{lemma}
	
	\begin{proof}
		Since $\zeta$ is even, $m(\sigma)=\int_\R\zeta(t)\cos(t\sigma)\,dt=\int_\R\zeta(t)e^{-\mathsf{i}t\sigma}\,dt$. For $k,N\in\Nzero$, differentiation in $\sigma$ and integration by parts in $t$ give $|\sigma|^N|m^{(k)}(\sigma)|\le\norm{\partial_t^N((-\mathsf{i}t)^k\zeta(t))}_{L^1(\R)}$, with no boundary terms because $\zeta$ is compactly supported. Every derivative of $m$ is bounded and decays faster than any inverse power, so $m$ is Schwartz. Moreover, $m(0)=\int_\R\zeta(t)\,dt=1$. 
		
		Fix $\eps>0$ and $r,M\in\R$, and choose an integer $q>\max\{M-r,0\}$. Since $m$ is Schwartz, $|m(\eps\sqrt\lambda)|\le C_{\eps,q}\lambda^{-q/2}$ for $\lambda\ge1$. The function $(1+\lambda)^{(M-r)/2}m(\eps\sqrt\lambda)$ is bounded at infinity because $q>M-r$, and it is bounded on $[0,1]$ by continuity. Set $C_{\eps,r,M}=\sup_{\lambda\ge0}(1+\lambda)^{(M-r)/2}|m(\eps\sqrt\lambda)|<\infty$. Let $u\in\HL^r$ and set $v=(1+L)^{r/2}u\in L^2(\R^n)$. The spectral theorem gives
		\begin{equation}\label{eq:spectral-smoothing-bound}
			\norm{S_\eps u}_{\HL^M}=\norm{(1+L)^{(M-r)/2}m(\eps\sqrt L)v}_2\le C_{\eps,r,M}\norm{v}_2=C_{\eps,r,M}\norm{u}_{\HL^r}.
		\end{equation}
		
		For $s>0$, the operator $L^s$ maps $\HL^r$ continuously into $\HL^{r-2s}$, and the two compositions $S_\eps L^s$ and $L^sS_\eps$ have the same spectral multiplier $\lambda^sm(\eps\sqrt\lambda)$. This proves commutation. Since $m(\eps\sqrt\lambda)\to1$ as $\eps\to0$ for every $\lambda\ge0$, and $|m(\eps\sqrt\lambda)|\le\norm{\zeta}_{L^1(\R)}$, dominated convergence in the spectral representation gives $S_\eps u\to u$ in $\HL^r$ as $\eps\to0$. By \eqref{eq:local-spectral-regularity}, the space $\HL^{2q}$ embeds continuously into $C^k(K)$ for every compact set $K\subset\R^n$ when $2q>k+n/2$. Combining this embedding with \eqref{eq:spectral-smoothing-bound} gives $S_\eps u\in C^\infty(\R^n)$ for every $u\in\HL^r$.
		
		By \eqref{eq:DG} and \cite[Theorem~3.4 and (3.5)]{CS08}, applied with the distance $|x-y|/c_L$, we have $(\cos(t\sqrt L)f,g)_{L^2}=0$ for $f,g\in L^2(\R^n)$ whenever $\dist(\supp f,\supp g)>c_L|t|$. Since $\supp\zeta\subset(-1,1)$, integrating this identity in \eqref{eq:spectral-smoothing} gives $(S_\eps f,g)_{L^2}=0$ whenever $\dist(\supp f,\supp g)>c_L\eps$. The support satisfies $\supp(S_\eps f)\subset\{x: \, \dist(x,\supp f)\le c_L\eps\}$ for every $f\in L^2(\R^n)$.
		
		We now prove density in the positive spectral spaces. Fix $r\ge0$ and $f\in\HL^r$, and choose $f_\ell\in C_c^\infty(\R^n)$ with $f_\ell\to f$ in $L^2(\R^n)$. For each fixed $\eps>0$, the smoothing and support properties proved above give $S_\eps f_\ell\in C_c^\infty(\R^n)$, and
		\[
		\norm{S_\eps f_\ell-f}_{\HL^r}\le C_{\eps,0,r}\norm{f_\ell-f}_2+\norm{S_\eps f-f}_{\HL^r}.
		\]
		Given $\eta>0$, first choose $\eps$ so that the second term is smaller than $\eta/2$, and then choose $\ell$ so that the first term is smaller than $\eta/2$. This proves the density of $C_c^\infty(\R^n)$ in $\HL^r$.
		
		To view elements of $\HL^{-r}$ as distributions, fix a compact set $K$ and an integer $q$ with $2q\ge r$. Smoothness and locality of $L$ give $\norm{\phi}_{\HL^r}\le\norm{(1+L)^q\phi}_2\le C_{K,q}\norm{\phi}_{H^{2q}}$ for $\phi\in C_c^\infty(\R^n)$ supported in $K$. By duality, every element of $\HL^{-r}$ acts continuously on test functions, and the density just proved makes this identification injective.
		
		For each integer $\ell\ge0$, choose $q$ with $2q>\ell+2+n/2$. Estimate~\eqref{eq:local-spectral-regularity} and Taylor's formula with integral remainder show that $y\mapsto\delta_y$ is locally $C^\ell$ with values in $\HL^{-2q}$. Estimates~\eqref{eq:spectral-smoothing-bound} and \eqref{eq:local-spectral-regularity} give a continuous map $S_\eps:\HL^{-2q}\to C^k(K)$ for every compact set $K\subset\R^n$ and every integer $k\ge0$. The kernel $k_\eps(x,y)=(S_\eps\delta_y)(x)$ is smooth in both variables, with $\partial_x^\alpha\partial_y^\beta k_\eps(x,y)=(-1)^{|\beta|}\partial_x^\alpha(S_\eps\partial^\beta\delta_y)(x)$.
		
		The multiplier $m$ is real, and $L$ has real coefficients. Self-adjointness of $S_\eps$ and its preservation of real-valued functions give $k_\eps(x,y)=k_\eps(y,x)\in\R$ and $S_\eps f(x)=\int_{\R^n}k_\eps(x,y)f(y)\,dy$ for $f\in C_c^\infty(\R^n)$. If $|x-y|>c_L\eps$, choose small neighborhoods of $x$ and $y$ whose distance exceeds $c_L\eps$. The vanishing of $(S_\eps f,g)_{L^2}$ for test functions supported in these neighborhoods forces the smooth kernel to vanish there. This proves the asserted support.
		
		The kernel defines $S_\eps$ on distributions by duality, since finite propagation gives $S_\eps\phi\in C_c^\infty(\R^n)$ for every $\phi\in C_c^\infty(\R^n)$. Its action is given by $\ip{S_\eps u}{\phi}=\ip{u}{S_\eps\phi}$ and, pointwise, by $S_\eps u(x)=\ip{u}{k_\eps(x,\cdot)}$. On the negative spectral spaces, this agrees with the spectral definition by self-adjointness.
	\end{proof}
	
	Suppose that $u=0$ on $\cO$. If $\dist(x,\R^n\setminus\cO)>c_L\eps$, the support property of $k_\eps$ gives $k_\eps(x,\cdot)\in C_c^\infty(\cO)$ and $S_\eps u(x)=\ip{u}{k_\eps(x,\cdot)}=0$. Now assume that \eqref{eq:entanglement} holds, and set $F=\sum_{k=1}^N b_kL^{s_k}u_k$. Let $\cO_0\Subset\cO$ be a nonempty open set and choose $\eps>0$ such that $c_L\eps<\dist(\overline{\cO_0},\R^n\setminus\cO)$. Applying the preceding argument to $u_1,\ldots,u_N$ and $F$, and using the commutation of $S_\eps$ with each $L^{s_k}$, gives
	\begin{equation}\label{eq:smoothed-entanglement}
		S_\eps u_1=\cdots=S_\eps u_N=\sum_{k=1}^N b_kL^{s_k}S_\eps u_k=0\quad\text{on }\cO_0.
	\end{equation}
	
	\begin{lemma}\label{lem:preserve-decay}
		Assume \eqref{eq:compact-perturbation}. If $u\in\mathscr D_\gamma$, then, for every fixed $\eps>0$, there is $a'>0$ such that $e^{a'|x|^\gamma}\partial^\alpha S_\eps u\in L^2(\R^n)$ for every multi-index $\alpha$.
	\end{lemma}
	
	\begin{proof}
		Let $u\in H^{-r}(\R^n)$ satisfy \eqref{eq:distribution-decay}, with $r\ge0$ and constants $a,C>0$. By \eqref{eq:spectral-sobolev-equivalence}, estimate~\eqref{eq:spectral-smoothing-bound} holds in the usual Sobolev spaces.
		
		Fix $\eps>0$ and a multi-index $\alpha$, and choose $M>n/2+|\alpha|$. The Fourier transform of $\partial^\alpha\delta_x$ gives $\norm{\partial^\alpha\delta_x}_{H^{-M}}^2=C_n\int_{\R^n}(1+|\xi|^2)^{-M}|\xi^\alpha|^2\,d\xi<\infty$, independently of $x$. Symmetry of the kernel gives $\partial_x^\alpha k_\eps(x,\cdot)=(-1)^{|\alpha|}S_\eps(\partial^\alpha\delta_x)$. Applying \eqref{eq:spectral-smoothing-bound}, transferred to the usual Sobolev spaces by \eqref{eq:spectral-sobolev-equivalence}, from $H^{-M}$ to $H^r$, we obtain
		\begin{equation}\label{eq:kernel-row}
			\sup_{x\in\R^n}\norm{\partial_x^\alpha k_\eps(x,\cdot)}_{H^r}\le C_{\alpha,\eps}.
		\end{equation}
		
		For sufficiently large $|x|$, set $R_x=|x|-c_L\eps-1\ge1$. The support of $\partial_x^\alpha k_\eps(x,\cdot)$ is contained in $\overline{B_{c_L\eps}(x)}$, where $|y|\ge R_x+1$. This function is a smooth compactly supported test function outside $\overline{B_{R_x}}$. Applying \eqref{eq:distribution-decay} and \eqref{eq:kernel-row} gives
		\[
		|\partial^\alpha S_\eps u(x)|=|\ip{u}{\partial_x^\alpha k_\eps(x,\cdot)}|\le C_{\alpha,\eps}e^{-aR_x^\gamma}.
		\]
		For $|x|\ge2(c_L\eps+1)$, we have $R_x\ge|x|/2$. Choose $0<a'<2^{-\gamma}a$. The preceding estimate gives $e^{a'|x|^\gamma}|\partial^\alpha S_\eps u(x)|\le C_{\alpha,\eps}e^{-(2^{-\gamma}a-a')|x|^\gamma}$ for sufficiently large $|x|$. The right-hand side is square integrable because $\gamma>0$. Smoothness gives integrability on bounded sets. The same choice of $a'$ works for every $\alpha$.
	\end{proof}
	
	\begin{corollary}\label{cor:single-distributions}
		Theorem~\ref{thm:single} holds for $u\in\HL^{-r}$ with any $r\ge0$. Under \eqref{eq:compact-perturbation}, it holds for every $u\in H^{-r}(\R^n)$ without a decay assumption.
	\end{corollary}
	
	\begin{proof}
		Choose a nonempty open ball $B\Subset\cO$ and set $d_B=\dist(\overline B,\R^n\setminus\cO)>0$. For $0<\eps<d_B/c_L$, Lemma~\ref{lem:smoothing} gives $S_\eps u\in L^2(\R^n)$, and \eqref{eq:smoothed-entanglement} gives $S_\eps u=\sum_{k=1}^N b_kL^{s_k}S_\eps u=0$ on $B$. Theorem~\ref{thm:single}, applied on $B$, gives $S_\eps u=0$ on $\R^n$. Since $S_\eps u\to u$ in $\HL^{-r}$ as $\eps\to0$, we conclude that $u=0$. Under \eqref{eq:compact-perturbation}, equation~\eqref{eq:spectral-sobolev-equivalence} gives the assertion for $u\in H^{-r}(\R^n)$.
	\end{proof}

	\section{Exterior observation}\label{sec:exterior}
	
	\begin{proof}[Proof of Theorem~\ref{thm:exterior}]
		Let $K=\bigcup_{k=1}^N\supp u_k$, and let $E$ be the unbounded connected component of $\R^n\setminus K$. We first assume that $u_k\in L^2(\R^n)$ for $k=1,\ldots,N$. Set $c_k=b_k\sin(\pi s_k)\ne0$ and $W_\lambda=\sum_{k=1}^N c_k\lambda^{s_k}(L+\lambda)^{-1}u_k$ for $\lambda>0$. By Proposition~\ref{prop:resolvent}, condition~\eqref{eq:entanglement} gives $W_\lambda=0$ on $\cO$. Moreover,
		\begin{equation}\label{eq:W-equation}
			(L+\lambda)W_\lambda=\sum_{k=1}^N c_k\lambda^{s_k}u_k.
		\end{equation}
		The right-hand side vanishes on $E$, so $(L+\lambda)W_\lambda=0$ there. Interior elliptic regularity gives $W_\lambda\in C^\infty(E)$. Since $E$ is connected and $\cO\cap E\ne\varnothing$, local unique continuation gives $W_\lambda=0$ on $E$ for every $\lambda>0$.
		
		Choose $R_0>0$ such that $K\subset B_{R_0}$. Since $n\ge2$, the set $\R^n\setminus\overline{B_{R_0}}$ is connected and contained in $E$. The vanishing of $W_\lambda$ on $E$ gives $\supp W_\lambda\subset\overline{B_{R_0}}$ for every $\lambda>0$. Choose $R>R_0$ such that $\cO\cap E\cap B_R\ne\varnothing$, set $D=B_R$, and fix a ball $\cO_0\Subset\cO\cap E\cap D$.
		
		Let $L_D$ denote the operator $-\Div(A\nabla)$ on $L^2(D)$ with domain $H^2(D)\cap H_0^1(D)$. For $h\in L^2(D)$ and $\lambda>0$, the function $v=(L_D+\lambda)^{-1}h$ is the unique weak solution of $(L+\lambda)v=h$ in $D$ with $v=0$ on $\partial D$. Since $W_\lambda\in H^1(\R^n)$ and is supported in $\overline{B_{R_0}}\Subset D$, it vanishes near $\partial D$ and $W_\lambda|_D\in H_0^1(D)$. Uniqueness applied to \eqref{eq:W-equation} gives
		\begin{equation}\label{eq:Dirichlet-combination}
			W_\lambda|_D=\sum_{k=1}^N c_k\lambda^{s_k}(L_D+\lambda)^{-1}u_k.
		\end{equation}
		
		Fix $\phi\in C_c^\infty(\cO_0)$ and set $f_k(z)=((L_D+z)^{-1}u_k,\phi)_{L^2(D)}$. The operator $L_D$ is nonnegative and self-adjoint with compact resolvent. Let $\mu_1>0$ denote its first eigenvalue. The spectral theorem gives $\norm{L_D^{-1}}_{L^2(D)\to L^2(D)}=\mu_1^{-1}$, so the Neumann series yields
		\[
		(L_D+z)^{-1}=L_D^{-1}(I+zL_D^{-1})^{-1}=\sum_{m=0}^{\infty}(-z)^mL_D^{-(m+1)},\quad |z|<\mu_1.
		\]
		The series converges in operator norm, locally uniformly in this disk. Applying it to $u_k$ and taking the inner product with $\phi$ shows that every $f_k$ is holomorphic there. Since $W_\lambda=0$ on $\cO_0$, equation~\eqref{eq:Dirichlet-combination} gives
		\begin{equation}\label{eq:Dirichlet-separation}
			\sum_{k=1}^N c_k\lambda^{s_k}f_k(\lambda)=0,\quad 0<\lambda<\mu_1.
		\end{equation}
		
		Suppose that some $f_k$ does not vanish identically near zero. For each such $k$, write $f_k(z)=a_kz^{m_k}+O(z^{m_k+1})$, where $a_k\ne0$ and $m_k\in\Nzero$. The leading exponents $s_k+m_k$ are distinct by \eqref{eq:exterior-orders}: equality for $k\ne\ell$ would give $s_k-s_\ell=m_\ell-m_k\in\Z$. Let $k_0$ be the unique index with the smallest leading exponent. Dividing \eqref{eq:Dirichlet-separation} by $\lambda^{s_{k_0}+m_{k_0}}$ and letting $\lambda\to0$, we obtain
		\[
		0=\lim_{\lambda\to0}\lambda^{-(s_{k_0}+m_{k_0})}\sum_{k=1}^N c_k\lambda^{s_k}f_k(\lambda)=c_{k_0}a_{k_0}.
		\]
		This contradicts $c_{k_0}a_{k_0}\ne0$. Every $f_k$ vanishes identically near zero.
		
		Let $0<\mu_1<\mu_2<\cdots$ be the distinct eigenvalues of $L_D$, and let $\Pi_{\mu_\ell}$ denote the orthogonal projection onto the eigenspace corresponding to $\mu_\ell$. The spectral theorem gives
		\[
		(L_D+z)^{-1}=\sum_{\ell=1}^{\infty}\frac{\Pi_{\mu_\ell}}{\mu_\ell+z}.
		\]
		Since $\mu_\ell\to\infty$, this series converges in operator norm, locally uniformly away from the points $-\mu_\ell$. The resolvent is meromorphic, with poles at these points. Each $f_k$ is holomorphic on the connected domain $\C\setminus\{-\mu_\ell:\ell\in\N\}$ and vanishes near zero. The identity theorem gives $f_k(z)=0$ throughout this domain.
		
		Fix an eigenvalue $\mu$ and set $\delta=\dist(\mu,\sigma(L_D)\setminus\{\mu\})>0$. For $0<|z+\mu|<\delta/2$, the spectral theorem gives
		\begin{align*}
			\norm{(z+\mu)(L_D+z)^{-1}-\Pi_\mu}_{L^2(D)\to L^2(D)}=\sup_{\substack{\lambda\in\sigma(L_D)\\\lambda\ne\mu}}\frac{|z+\mu|}{|\lambda+z|}\le\frac{2|z+\mu|}{\delta}.
		\end{align*}
		We obtain $(z+\mu)(L_D+z)^{-1}\to\Pi_\mu$ in operator norm as $z\to-\mu$ and
		\[
		(\Pi_\mu u_k,\phi)_{L^2(D)}=\lim_{z\to-\mu}(z+\mu)f_k(z)=0.
		\]
		Since $\phi\in C_c^\infty(\cO_0)$ was arbitrary, $\Pi_\mu u_k=0$ on $\cO_0$. It satisfies $(L-\mu)\Pi_\mu u_k=0$ in $D$, so local unique continuation gives $\Pi_\mu u_k=0$ throughout $D$. This holds for every eigenvalue $\mu$, and completeness of the eigenspaces gives $u_k=0$ in $L^2(D)$. Since $\supp u_k\subset K\subset D$, we conclude that $u_k=0$ on $\R^n$ for $k=1,\ldots,N$.
		
		We now consider $u_k\in\cE(\R^n)$. Choose $r\ge0$ such that $u_k\in\HL^{-r}$ for $k=1,\ldots,N$, and set $u_{k,\eps}=S_\eps u_k$. Lemma~\ref{lem:smoothing} gives $u_{k,\eps}\in C^\infty(\R^n)\cap L^2(\R^n)$ and $\supp u_{k,\eps}\subset K_\eps$, where $K_\eps=\{x:\dist(x,K)\le c_L\eps\}$ and $c_L=\sqrt{a_1}$.
		
		Choose $x_0\in\cO\cap E$ and a ball $\cO_0$ centered at $x_0$ with $\overline{\cO_0}\subset\cO\cap E$. Since $E$ is open and connected, there is a compact path $\Gamma\subset E$ joining $x_0$ to a point $y$ with $|y|>R_0$. The compact sets $\Gamma$ and $\overline{\cO_0}$ are disjoint from $K$, so both have positive distance from $K$. The ray $\{ty:t\ge1\}$ also has positive distance from $K$, since $\dist(ty,K)\ge|y|-R_0>0$ for $t\ge1$.
		
		For all sufficiently small $\eps>0$, the set $K_\eps$ is disjoint from $\overline{\cO_0}$ and $\Gamma\cup\{ty:t\ge1\}$, and $c_L\eps<\dist(\overline{\cO_0},\R^n\setminus\cO)$. Put $\widetilde K_\eps=\bigcup_{k=1}^N\supp u_{k,\eps}\subset K_\eps$. The curve $\Gamma\cup\{ty:t\ge1\}$ joins $x_0$ to infinity without meeting $\widetilde K_\eps$. Since $\cO_0$ is connected, contains $x_0$, and is disjoint from $\widetilde K_\eps$, it lies in the unbounded connected component of $\R^n\setminus\widetilde K_\eps$.
		
		By \eqref{eq:smoothed-entanglement}, the functions $u_{1,\eps},\ldots,u_{N,\eps}$ satisfy \eqref{eq:entanglement} on $\cO_0$. The $L^2$ case proved above, applied with observation set $\cO_0$, gives $u_{k,\eps}=0$ on $\R^n$ for $k=1,\ldots,N$. By Lemma~\ref{lem:smoothing}, $S_\eps u_k\to u_k$ in $\HL^{-r}$ as $\eps\to0$, which proves $u_k=0$ on $\R^n$ for every $k=1,\ldots,N$.
	\end{proof}

	\section{The resolvent near zero}\label{sec:low-energy}

	For arbitrary observation, the proof requires estimates for the resolvent near zero on the whole space. We obtain them from an equation for its gradient on a ball containing the coefficient perturbation. An energy estimate gives invertibility of the operator in this equation at zero using uniform ellipticity, without a smallness assumption on the perturbation.
	
	In this section, we assume \eqref{eq:compact-perturbation} and write $B=A-I$. Fix a ball $D$ with $\supp B\Subset D$, and set $\mathcal V=L^2(D;\C^n)$. We extend elements of $\mathcal V$ by zero outside $D$ whenever they are used on $\R^n$. For $\Rea\kappa>0$, write $\mathcal R_A(\kappa)=(L+\kappa^2)^{-1}$ and $\mathcal R_0(\kappa)=(-\Delta+\kappa^2)^{-1}$. Both resolvents are holomorphic in $\kappa$ as bounded operators on $L^2(\R^n)$.
	
	\subsection{An equation for the gradient}
	
	Regard $B$ as a multiplication operator on $\mathcal V$. For $\Rea\kappa>0$, define
	\begin{equation}\label{eq:K-T}
		K(\kappa)f=\left.\nabla\mathcal R_0(\kappa)\Div f\right|_D,
		\quad T(\kappa)=I-K(\kappa)B.
	\end{equation}
	The whole space operator $\nabla\mathcal R_0(\kappa)\Div$ has Fourier multiplier $-\xi\otimes\xi/(|\xi|^2+\kappa^2)$, so $K(\kappa)$ and $T(\kappa)$ are bounded and holomorphic on this half-plane.
	
	For $\kappa>0$ and $u\in L^2(\R^n)$, set 
	\begin{equation}\label{eq:def_X_u_g_u}
		X_u(\kappa)=\nabla\mathcal R_A(\kappa)u|_D \quad \text{and}\quad g_u(\kappa)=\nabla\mathcal R_0(\kappa)u|_D.
	\end{equation}
	Since $L=-\Delta-\Div(B\nabla)$, the resolvent equation gives
	\begin{equation}\label{eq:resolvent-gradient}
		\mathcal R_A(\kappa)u=\mathcal R_0(\kappa)u+\mathcal R_0(\kappa)\Div(BX_u(\kappa)).
	\end{equation}
	Taking the gradient and restricting to $D$, we obtain
	\begin{equation}\label{eq:gradient-equation}
		T(\kappa)X_u(\kappa)=g_u(\kappa),
	\end{equation}
	where $T(\kappa)$ is given by \eqref{eq:K-T}.
	
	At $\kappa=0$, define $K(0)f$ by applying the Fourier multiplier $-\xi\otimes\xi/|\xi|^2$, for $\xi\ne0$, to the zero extension of $f$ and restricting the result to $D$. Plancherel's theorem gives $\norm{K(0)}_{\mathcal V\to\mathcal V}\le1$ in every dimension $n\ge2$. Using \eqref{eq:K-T} again, we set $T(0)=I-K(0)B$.
	
	\begin{lemma}\label{lem:zero-inverse}
		The operator $T(0)$ is invertible on $\mathcal V$, and
		\begin{equation}\label{eq:zero-inverse-bound}
			\norm{T(0)^{-1}}_{\mathcal V\to\mathcal V}\le1+a_0^{-1}\norm{B}_{L^\infty}.
		\end{equation}
	\end{lemma}
	
	\begin{proof}
		Let $\dot H^1$ denote the completion of $C_c^\infty(\R^n)$ under the gradient norm, with functions understood modulo constants. Given $f\in\mathcal V$, the Lax--Milgram theorem gives a unique $v\in\dot H^1$ satisfying
		\[
		\int_{\R^n}A\nabla v\cdot\overline{\nabla w}\,dx=-\int_{\R^n}Bf\cdot\overline{\nabla w}\,dx,\quad w\in\dot H^1.
		\]
		$Lv=\Div(Bf)$ on $\R^n$. Taking $w=v$ and using ellipticity \eqref{eq:ellipticity} gives $\norm{\nabla v}_2\le a_0^{-1}\norm{Bf}_2$.
		
		Set $X=f+\nabla v|_D$. Since $B$ is supported in $D$, the equation for $v$ becomes $-\Delta v=\Div(BX)$ on $\R^n$. Taking Fourier transforms gives $\nabla v|_D=K(0)BX$, which is equivalent to $T(0)X=f$. Moreover, $\norm{X}_{\mathcal V}\le(1+a_0^{-1}\norm{B}_{L^\infty})\norm{f}_{\mathcal V}$. This proves surjectivity.
		
		For injectivity, suppose $T(0)X=0$, and let $v\in\dot H^1$ solve $-\Delta v=\Div(BX)$ on $\R^n$. Then $\nabla v|_D=K(0)BX=X$, so $Lv=-\Delta v-\Div(B\nabla v)=0$. Testing against $v$ gives $\int_{\R^n}A\nabla v\cdot\overline{\nabla v}\,dx=0$. Ellipticity gives $\nabla v=0$ and $X=0$. This proves injectivity, and the estimate for $X$ above gives \eqref{eq:zero-inverse-bound}.
	\end{proof}
	
	The proof uses only gradients in $L^2(\R^n)$, and it also applies in dimension two, where the scalar free resolvent has a logarithmic singularity at zero.
	
	For $F,G\in\mathcal V$, let $[F,G]_D=\int_D F\cdot G\,dx$ denote the complex bilinear pairing. Let $\kappa>0$, $u\in L^2(\R^n)$, and $\phi\in C_c^\infty(\R^n)$. If $T(\kappa)$ is invertible, equation~\eqref{eq:gradient-equation} gives $X_u(\kappa)=T(\kappa)^{-1}g_u(\kappa)$. Integration by parts in \eqref{eq:resolvent-gradient} then yields
	\begin{equation}\label{eq:scalar-reconstruction}
		\ip{\mathcal R_A(\kappa)u}{\phi}=\ip{\mathcal R_0(\kappa)u}{\phi}-[BT(\kappa)^{-1}g_u(\kappa),g_{\overline\phi}(\kappa)]_D.
	\end{equation}
	For positive $\kappa$, the free resolvent commutes with complex conjugation, so $g_{\overline\phi}(\kappa)=\overline{g_\phi(\kappa)}$. For complex $\kappa$ with $\Rea\kappa>0$, both $g_u(\kappa)$ and $g_{\overline\phi}(\kappa)$ are holomorphic. Since $[\cdot,\cdot]_D$ is bilinear, both sides of \eqref{eq:scalar-reconstruction} are holomorphic on any right half-disk $U_\eps=\{\kappa:\Rea\kappa>0,\ |\kappa|<\eps\}$, where $T(\kappa)$ is invertible. The identity theorem extends \eqref{eq:scalar-reconstruction} to $U_\eps$. In the following, the localized free resolvent formulas will show that such an $\eps>0$ exists.
	
	\subsection{Localized free resolvents}
	
	For $\Rea\kappa>0$, the integral kernel of $\mathcal R_0(\kappa)$ is $G_n(\kappa,|x-y|)$, where
	\begin{equation}\label{eq:free-kernel}
		G_n(\kappa,r)=(2\pi)^{-n/2}\left(\frac{\kappa}{r}\right)^{n/2-1}K_{n/2-1}(\kappa r),\quad r>0.
	\end{equation}
	Here $K_\mu$ is the modified Bessel function of the second kind. We use the half-integer formulas and integer-order series in \cite[Sections 10.31 and 10.49]{NIST10}, together with the recurrence
	\begin{equation}\label{eq:kernel-recurrence}
		G_{n+2}(\kappa,r)=-\frac1{2\pi r}\partial_rG_n(\kappa,r).
	\end{equation}
	This follows from $K_\mu'(w)=\mu w^{-1}K_\mu(w)-K_{\mu+1}(w)$; see \cite[Eq.~(10.29.2)]{NIST10}. For bounded open sets $U,V\subset\R^n$, localization means that the source is extended by zero from $U$ and the output is restricted to $V$.

	\begin{lemma}[Holomorphic extension in odd dimensions]\label{lem:odd-compact}
		Let $n\ge3$ be odd, and let $U,V\subset\R^n$ be bounded open sets. The localized operators $\mathcal R_0(\kappa)$ and $\nabla\mathcal R_0(\kappa)$ extend holomorphically near zero with values in $\mathcal L(L^2(U),L^2(V))$ and $\mathcal L(L^2(U),L^2(V;\C^n))$, respectively. The operator $K(\kappa)$ extends holomorphically near zero in $\mathcal L(\mathcal V)$, with value $K(0)$. The inverse $T(\kappa)^{-1}$ is holomorphic on a sufficiently small disk centered at zero.
	\end{lemma}
	
	\begin{proof}
		Write $n=2q+3$ with $q\in \Nzero$. Substituting the half integer formula \cite[Eqs.~(10.47.9), (10.49.1), and (10.49.12)]{NIST10} into \eqref{eq:free-kernel}, we obtain
		\begin{equation}\label{eq:odd-kernel}
			G_n(\kappa,r)=c_n r^{2-n}e^{-\kappa r}P_q(\kappa r),
		\end{equation}
		where $c_n=1/(2^{q+2}\pi^{q+1})$ and $P_q(w)=\sum_{m=0}^{q}\frac{(q+m)!}{2^m m!(q-m)!}w^{q-m}$. In particular, $G_3(\kappa,r)=e^{-\kappa r}/(4\pi r)$.
		
		Choose $R>\sup_{x\in V,y\in U}|x-y|$. Formula~\eqref{eq:odd-kernel} and its derivative give, for every $\rho>0$,
		\begin{equation}\label{eq:odd-local-kernel-bounds}
			|G_n(\kappa,|x-y|)|\le C_\rho|x-y|^{2-n},\quad
			|\nabla_xG_n(\kappa,|x-y|)|\le C_\rho|x-y|^{1-n},
		\end{equation}
		when $x\in V$, $y\in U$, $x\ne y$, and $|\kappa|\le\rho$.
		
		For a scalar or vector kernel $H(x,y)$, let $T_Hf(x)=\int_U H(x,y)f(y)\,dy$. Schur's estimate states that
		\begin{equation}\label{eq:schur-estimate}
			\|T_H\|\le
			\Big(\sup_{x\in V}\int_U|H(x,y)|\,dy\Big)^{1/2}
			\Big(\sup_{y\in U}\int_V|H(x,y)|\,dx\Big)^{1/2}.
		\end{equation}
		Indeed, writing the two suprema as $a,b<\infty$, Cauchy--Schwarz gives $|T_Hf(x)|^2\le a\int_U|H(x,y)||f(y)|^2\,dy$. Integration in $x$ and Tonelli's theorem then give $\|T_Hf\|_{L^2(V)}^2\le ab\|f\|_{L^2(U)}^2$. For the kernels in \eqref{eq:odd-local-kernel-bounds}, both suprema are uniformly bounded for $|\kappa|\le\rho$, since the spatial majorants are integrable over $B_R(0)$.
		
		For either kernel $H(\kappa,x,y)$, its dependence on $\kappa$ is entire whenever $x\ne y$, and Cauchy's integral formula gives its Taylor coefficients at zero:
		\[
		H_m(x,y)=\frac{1}{2\pi \mathsf{i}}\int_{|\zeta|=\rho}\frac{H(\zeta,x,y)}{\zeta^{m+1}}\,d\zeta.
		\]
		The integration circle $|\zeta|=\rho$ has length $2\pi\rho$, and $|\zeta^{-(m+1)}|=\rho^{-(m+1)}$ on that circle. Combining these facts with \eqref{eq:odd-local-kernel-bounds} gives $|H_m(x,y)|\le C_\rho\rho^{-m}|x-y|^{2-n}$ for the scalar kernel and $|H_m(x,y)|\le C_\rho\rho^{-m}|x-y|^{1-n}$ for the gradient kernel. By \eqref{eq:schur-estimate}, the corresponding coefficient operators satisfy $\|T_{H_m}\|\le C_\rho\rho^{-m}$.
		
		For any $0<\rho'<\rho$, the pointwise Taylor remainder on $|\kappa|\le\rho'$ is bounded by $C_\rho\sum_{m=N+1}^\infty(\rho'/\rho)^m$ times $|x-y|^{2-n}$ or $|x-y|^{1-n}$, respectively. Applying \eqref{eq:schur-estimate} to this remainder kernel gives
		\[
		\sup_{|\kappa|\le\rho'}\Big\|T_H(\kappa)-\sum_{m=0}^{N}\kappa^mT_{H_m}\Big\|\le C_\rho\frac{(\rho'/\rho)^{N+1}}{1-\rho'/\rho}.
		\]
		The right-hand side tends to zero as $N\to\infty$, so the operator power series converges to $T_H(\kappa)$ in operator norm, locally uniformly for $|\kappa|<\rho$. Since $\rho>0$ is arbitrary, both operator families are entire, proving the first two assertions.
		
		To treat $K(\kappa)$, choose a ball $\Omega$ with $D\Subset\Omega$. For $f\in\mathcal V$, extended by zero outside $D$, define
		\[
		v_\kappa(f)(x)=\int_D\nabla_xG_n(\kappa,|x-y|)\cdot f(y)\,dy,\quad x\in\Omega.
		\]
		The preceding argument shows that $v_\kappa$ is holomorphic with values in $\mathcal L(\mathcal V,L^2(\Omega))$. For $\Rea\kappa>0$, integration by parts identifies $v_\kappa(f)$ with $\mathcal R_0(\kappa)\Div f$, while $v_0(f)$ is defined by the Newton kernel. Set $w_\kappa(f)=v_\kappa(f)-v_0(f)$. Subtracting their equations gives
		\[
		-\Delta w_\kappa(f)=-\kappa^2v_\kappa(f)\quad\text{in }\Omega,
		\]
		initially for $\Rea\kappa>0$, and then for all $\kappa$ by holomorphy.
		
		For $w,\Delta w\in L^2(\Omega)$, interior elliptic regularity gives
		\begin{equation}\label{eq:localized-interior-estimate}
			\norm{w}_{H^2(D)}\le C(\norm{w}_{L^2(\Omega)}+\norm{\Delta w}_{L^2(\Omega)}).
		\end{equation}
		Both $w_\kappa$ and $\Delta w_\kappa=\kappa^2v_\kappa$ are holomorphic with values in $\mathcal L(\mathcal V,L^2(\Omega))$. Applying \eqref{eq:localized-interior-estimate} to differences of their difference quotients shows that $w_\kappa$ is holomorphic with values in $\mathcal L(\mathcal V,H^2(D))$. Since the Newton kernel gives $\nabla v_0(f)|_D=K(0)f$,
		\[
		K(\kappa)f=K(0)f+\nabla w_\kappa(f)|_D
		\]
		provides the required holomorphic extension, with value $K(0)$ at zero because $w_0=0$.
		
		Finally, Lemma~\ref{lem:zero-inverse} gives a bounded inverse for $T(0)$. Since $K(\kappa)\to K(0)$ in operator norm, the Neumann series shows that
		\[
		T(\kappa)^{-1}=(I-T(0)^{-1}(K(\kappa)-K(0))B)^{-1}T(0)^{-1}
		\]
		is defined and holomorphic on a sufficiently small disk.
	\end{proof}
	
	For even $n=2\nu\ge2$ with $\nu\in \N$, write $z=\kappa^2$ and define the scalar function
	\[
	\mathcal B_n(z,r)=\sum_{m=0}^{\infty}\frac{(-1)^\nu}{(4\pi)^\nu4^m m!(m+\nu-1)!}z^mr^{2m},
	\quad z\in\C,\quad r>0.
	\]
	In particular, $\mathcal B_n(0,r)=(-1)^\nu/((4\pi)^\nu(\nu-1)!)$ is a constant independent of $r$. For $z\ne0$ with $-\pi<\arg z<\pi$, define $\mathcal A_n(z,r)$ by
	\begin{equation}\label{eq:even-kernel}
		G_n(\sqrt z,r)=\mathcal A_n(z,r)+z^{\nu-1}\log z\,\mathcal B_n(z,r),
	\end{equation}
	using the principal branches of $\sqrt z$ and $\log z$.
	
	For bounded open sets $U,V\subset\R^n$ and either $\mathcal C_n=\mathcal A_n$ or $\mathcal C_n=\mathcal B_n$, the corresponding localized scalar and gradient operators are
	\[
	\begin{aligned}
		f\longmapsto\int_U\mathcal C_n(z,|x-y|)f(y)\,dy,\quad f\longmapsto\int_U\nabla_x\mathcal C_n(z,|x-y|)f(y)\,dy,
	\end{aligned}
	\]
	with $x\in V$, initially defined for $f\in C_c^\infty(U)$.
	
	\begin{lemma}[Logarithmic expansion in even dimensions]\label{lem:even-compact}
		There is $\rho>0$ such that, for every fixed $r>0$, $\mathcal A_n(z,r)$ extends holomorphically to $|z|<\rho$. For any bounded open sets $U,V$, the localized scalar and gradient operators defined above extend to holomorphic families on this disk with values in $\mathcal L(L^2(U),L^2(V))$ and $\mathcal L(L^2(U),L^2(V;\C^n))$, respectively. There are holomorphic operator families $K_1,K_2$ with values in $\mathcal L(\mathcal V)$ such that
		\begin{equation}\label{eq:even-K}
			K(\sqrt z)=K(0)+zK_1(z)+z^\nu\log z\,K_2(z)
		\end{equation}
		for $0<|z|<\rho$ and $-\pi<\arg z<\pi$.
	\end{lemma}
	
	\begin{proof}
		In dimension two, $G_2(\kappa,r)=(2\pi)^{-1}K_0(\kappa r)$. The convergent series \cite[Eq.~(10.31.2)]{NIST10} gives
		\begin{equation}\label{eq:K0-series}
			\begin{aligned}
				K_0(w)=-(\log(w/2)+\gamma_E)I_0(w)+\sum_{m=1}^{\infty}H_m\frac{(w^2/4)^m}{(m!)^2},\quad  I_0(w)=\sum_{m=0}^{\infty}\frac{(w^2/4)^m}{(m!)^2},
			\end{aligned}
		\end{equation}
		where $\gamma_E$ is Euler's constant, $H_0=0$, and $H_m=\sum_{\ell=1}^{m}\ell^{-1}$ for $m\ge1$. Splitting $\log(\sqrt z r/2)=\frac12\log z+\log(r/2)$ proves \eqref{eq:even-kernel} for $\nu=1$, with
		\[
		\begin{aligned}
			\mathcal A_2(z,r)=\frac1{2\pi}\sum_{m=0}^{\infty}(H_m-\log(r/2)-\gamma_E)\frac{(zr^2/4)^m}{(m!)^2},\quad
			\mathcal B_2(z,r)=-\frac1{4\pi}I_0(\sqrt z r).
		\end{aligned}
		\]
		Both expressions are entire in $z$ for each $r>0$, as follows directly from their series.
		
		Set $\mathcal D_r=-(2\pi r)^{-1}\partial_r$. The recurrence \eqref{eq:kernel-recurrence} gives $\mathcal A_n=\mathcal D_r^{\nu-1}\mathcal A_2$. To compute the logarithmic coefficient, note that $\mathcal D_r^{\nu-1}$ kills $r^{2k}$ when $0\le k<\nu-1$, while
		\[
		\mathcal D_r^{\nu-1}r^{2(m+\nu-1)}=\frac{(-1)^{\nu-1}(m+\nu-1)!}{\pi^{\nu-1}m!}r^{2m}.
		\]
		Applying this identity to the series for $\mathcal B_2$ yields $\mathcal D_r^{\nu-1}\mathcal B_2=z^{\nu-1}\mathcal B_n$, where
		\begin{equation}\label{eq:B-series}
			\mathcal B_n(z,r)=\sum_{m=0}^{\infty}\mathfrak c_mz^mr^{2m},
			\quad \mathfrak c_m=\frac{(-1)^\nu}{(4\pi)^\nu4^m m!(m+\nu-1)!}.
		\end{equation}
		This proves \eqref{eq:even-kernel} in every even dimension and shows that $\mathcal B_n(0,r)=\mathfrak c_0$ is constant.
		
		We next verify the localized operator bounds. In the series for $\mathcal A_2$, applying $\mathcal D_r^{\nu-1}$ to the terms with $m<\nu-1$ gives finitely many singular powers of $r$. The remaining terms give convergent series with coefficients containing $\log r$. The factorial denominators ensure convergence after the required spatial derivatives. For fixed $R>0$, $0<\rho'<\rho$, $0<r\le R$, and $|z|\le\rho'$, we have
		\[
		\begin{aligned}
			|\mathcal A_2(z,r)|&\le C(1+|\log r|),
			&|\partial_r\mathcal A_2(z,r)|&\le Cr^{-1},\\
			|\mathcal A_n(z,r)|&\le Cr^{2-n},
			&|\partial_r\mathcal A_n(z,r)|&\le Cr^{1-n}\quad(n\ge4).
		\end{aligned}
		\]
		The series \eqref{eq:B-series} also gives $|\mathcal B_n(z,r)|\le C$ and $|\partial_r\mathcal B_n(z,r)|\le C|z|r$. All the displayed spatial bounds are locally integrable. Applying Cauchy's integral formula in $z$ and Schur's estimate \eqref{eq:schur-estimate} gives the stated operator holomorphy by the coefficient argument in the proof of Lemma~\ref{lem:odd-compact}.
		
		To prove \eqref{eq:even-K}, we use the equations satisfied by $\mathcal A_n$ and $\mathcal B_n$. Since $\Delta r^{2m}=4m(m+\nu-1)r^{2m-2}$ and $4(m+1)(m+\nu)\mathfrak c_{m+1}=\mathfrak c_m$, equation~\eqref{eq:B-series} gives $(-\Delta+z)\mathcal B_n(z,|\cdot|)=0$. The function $\mathcal B_n(z,|\cdot|)$ is smooth at the origin, so this equation also holds distributionally. Combining it with the fundamental solution identity for $G_n$, we obtain
		\[
		(-\Delta+z)\mathcal A_n(z,|\cdot|)=\delta_0.
		\]
		This identity initially holds in the slit disk and extends throughout $|z|<\rho$ by holomorphy. At $z=0$, $\mathcal A_n(0,r)$ is the Newton kernel for $n\ge4$, while $\mathcal A_2(0,r)=-(2\pi)^{-1}(\log(r/2)+\gamma_E)$ differs from the two-dimensional Newton kernel by a constant.
		
		Choose a ball $\Omega$ with $D\Subset\Omega$. For $f\in\mathcal V$, extended by zero outside $D$, define
		\[
		\begin{aligned}
			(v_A(z)f)(x)=\int_D\nabla_x\mathcal A_n(z,|x-y|)\cdot f(y)\,dy,\quad 	(v_B(z)f)(x)=\int_D\nabla_x\mathcal B_n(z,|x-y|)\cdot f(y)\,dy,
		\end{aligned}
		\]
		for $x\in\Omega$. These are holomorphic operators from $\mathcal V$ into $L^2(\Omega)$ and satisfy
		\[
		(-\Delta+z)v_A(z)f=\Div f,
		\quad (-\Delta+z)v_B(z)f=0.
		\]
		Since $\mathcal B_n(0,r)$ is constant, $v_B(0)=0$. The operator quotients
		\[
		Q_1(z)=\frac{v_A(z)-v_A(0)}z,
		\quad Q_2(z)=\frac{v_B(z)-v_B(0)}z=\frac{v_B(z)}z
		\]
		have removable singularities at zero and extend holomorphically with values in $\mathcal L(\mathcal V,L^2(\Omega))$. Subtracting the equations for $v_A(z)f$ and $v_A(0)f$, and using the equation for $v_B(z)f$, gives
		\[
		\Delta(Q_1(z)f)=v_A(z)f,
		\quad \Delta(Q_2(z)f)=zQ_2(z)f.
		\]
		Both $Q_1,Q_2$ and their Laplacians are holomorphic with values in $\mathcal L(\mathcal V,L^2(\Omega))$. Applying \eqref{eq:localized-interior-estimate} to their parameter difference quotients proves that $Q_1,Q_2$ are holomorphic with values in $\mathcal L(\mathcal V,H^2(D))$.
		
		Define $K_1(z)f=\nabla(Q_1(z)f)|_D$ and $K_2(z)f=\nabla(Q_2(z)f)|_D$. The Newton kernel gives $\nabla(v_A(0)f)|_D=K(0)f$. Taking the gradient in the decomposition \eqref{eq:even-kernel}, we obtain
		\[
		K(\sqrt z)f=K(0)f+zK_1(z)f+z^\nu\log z\,K_2(z)f.
		\]
		This proves \eqref{eq:even-K}.
	\end{proof}
	
	On a smaller parameter disk, the holomorphic operators $K_1,K_2$ are bounded, so \eqref{eq:even-K} gives
	\[
	\norm{K(\sqrt z)-K(0)}_{\mathcal V\to\mathcal V}\le C (|z|+|z|^\nu(|\log|z||+\pi))\to0
	\]
	as $z\to0$ in the slit disk. For $\kappa=\mathsf{i}t\ne0$, the kernels are defined by limits from $\Rea\kappa>0$: $z=\kappa^2$ approaches the negative real axis from above when $t>0$ and from below when $t<0$, with $\log z=\log(t^2)+\mathsf{i}\pi\operatorname{sgn}(t)$. The same norm bound holds for these boundary values. Lemma~\ref{lem:zero-inverse} and a Neumann series about $T(0)$ show that $T(\kappa)^{-1}$ exists for all sufficiently small nonzero $\kappa$ in the closed right half-plane and converges in operator norm to $T(0)^{-1}$ as $\kappa\to0$.
	
	\subsection{Root-exponential source estimates}
	
	For $a,\gamma>0$, write $L^2_{a,\gamma}=\{u:e^{a|x|^\gamma}u\in L^2(\R^n)\}$. We first estimate the polynomial moments of a source in $L^2_{a,1/2}$. These estimates will control the parameter derivatives in odd dimensions and the coefficients and remainders of the logarithmic expansion in even dimensions.
	
	\begin{lemma}[Weighted moment estimates]\label{lem:weighted-moments}
		Let $a>0$ and $d\ge0$. There exist constants $C_0,C_1>0$, depending only on $n,a,d$, such that
		\begin{equation}\label{eq:weighted-moments}
			\int_{\R^n}(1+|y|)^{m+d}|u(y)|\,dy\le C_0C_1^m(m!)^2\norm{e^{a\sqrt{|\cdot|}}u}_2
		\end{equation}
		for every $m\in\Nzero$ and $u\in L^2_{a,1/2}$.
	\end{lemma}
	
	\begin{proof}
		By the Cauchy--Schwarz inequality, the left-hand side of \eqref{eq:weighted-moments} is bounded by $I_m^{1/2}\norm{e^{a\sqrt{|\cdot|}}u}_2$, where $I_m=\int_{\R^n}(1+|y|)^{2m+2d}e^{-2a\sqrt{|y|}}\,dy$. Splitting the radial integral at $r=1$ and then substituting $r=t^2$ gives
		\[
		I_m\le C^{m+1}\Big(1+\int_0^\infty r^{2m+2d+n-1}e^{-2a\sqrt r}\,dr\Big)=C^{m+1}\Big(1+\frac{2\Gamma(4m+4d+2n)}{(2a)^{4m+4d+2n}}\Big).
		\]
		For fixed $\sigma>0$, Stirling's formula gives $\Gamma(4m+\sigma)^{1/2}\le C_\sigma16^m(m+1)^{\sigma/2}(m!)^2$. Apply this with $\sigma=4d+2n$. The factors depending exponentially on $m$, including $(2a)^{-2m}$, are included in $C_1^m$. The remaining fixed power of $m+1$ can also be included in this constant, since $(m+1)^{\sigma/2}\le C_\sigma2^m$. $I_m^{1/2}\le C_0C_1^m(m!)^2$, which proves the result.
	\end{proof}
	
	We first consider odd dimensions.
	
	\begin{proposition}[Gevrey estimates in odd dimensions]\label{prop:Gevrey}
		Let $n\ge3$ be odd, $a>0$, $u\in L^2_{a,1/2}$, and $\phi\in C_c^\infty(\R^n)$. There are constants $\eps,C_0,C_1>0$ such that $h(\kappa)=\ip{\mathcal R_A(\kappa)u}{\phi}$ is holomorphic in $U_\eps=\{\kappa:\Rea\kappa>0,\ |\kappa|<\eps\}$, each derivative $h^{(m)}$ extends continuously to $\overline{U_\eps}$, and
		\begin{equation}\label{eq:Gevrey}
			\sup_{\overline{U_\eps}}|h^{(m)}(\kappa)|\le C_0C_1^m(m!)^2,\quad m\in\Nzero.
		\end{equation}
		The constants $C_0,C_1$ may depend on $n,a,A,D,u,\phi$, but are independent of $m$ and $\kappa$.
	\end{proposition}
	
	\begin{proof}
		Choose $\eps>0$ so that $T(\kappa)^{-1}$ is holomorphic for $|\kappa|<4\eps$. Choose balls $D_1\Subset D_2$ such that $\overline D\cup\supp\phi\subset D_1$, and write $u=u_1+u_2$, where $u_1=u\mathbf1_{D_2}$ and $u_2=u\mathbf1_{\R^n\setminus D_2}$. Here $\mathbf1_E$ denotes the characteristic function of a set $E$. Throughout the proof, $C_0,C_1$ denote positive constants that may change from one estimate to another but are always independent of $m$ and $\kappa$. We estimate $u_1$ and $u_2$ separately.
		
		For $u_1$, both the source and the observation variables lie in fixed bounded sets. By Lemma~\ref{lem:odd-compact}, the scalar function $\kappa\mapsto\ip{\mathcal R_0(\kappa)u_1}{\phi}$ and the $\mathcal V$-valued function $\kappa\mapsto g_{u_1}(\kappa)$ admit holomorphic extensions to $|\kappa|<4\eps$, where $g_u$ is defined by \eqref{eq:def_X_u_g_u}. We use the same notation for these extensions. Choose $C_0>0$ such that
		\begin{equation}\label{eq:holo_estimate_C_0}
			\sup_{|\eta|\le3\eps}\left(\left|\ip{\mathcal R_0(\eta)u_1}{\phi}\right|+\norm{g_{u_1}(\eta)}_{\mathcal V}\right)\le C_0.
		\end{equation}
		Cauchy's integral formula holds for both scalar-valued and Banach space valued holomorphic functions. For $|\kappa|\le2\eps$, the circle $|\eta-\kappa|=\eps$ lies in $|\eta|\le3\eps$. With counterclockwise orientation, the formula gives
		\[
		\begin{aligned}
			\partial_\kappa^m\ip{\mathcal R_0(\kappa)u_1}{\phi}=\frac{m!}{2\pi\mathsf{i}}\int_{|\eta-\kappa|=\eps}\frac{\ip{\mathcal R_0(\eta)u_1}{\phi}}{(\eta-\kappa)^{m+1}}\,d\eta,\quad 
			\partial_\kappa^m g_{u_1}(\kappa)	=\frac{m!}{2\pi\mathsf{i}}\int_{|\eta-\kappa|=\eps}\frac{g_{u_1}(\eta)}{(\eta-\kappa)^{m+1}}\,d\eta.
		\end{aligned}
		\]
		Taking the absolute value in the first formula and the $\mathcal V$-norm in the second, and using \eqref{eq:holo_estimate_C_0}, we obtain
		\[
		\begin{aligned}
			\left|\partial_\kappa^m\ip{\mathcal R_0(\kappa)u_1}{\phi}\right|+\norm{\partial_\kappa^m g_{u_1}(\kappa)}_{\mathcal V}\le\frac{m!}{2\pi}\frac{C_0}{\eps^{m+1}}(2\pi\eps)=C_0\eps^{-m}m!.
		\end{aligned}
		\]
		For every $C_1\ge\eps^{-1}$,
		\begin{equation}\label{eq:estimate_for_u_1}
			\sup_{|\kappa|\le2\eps}\left(\left|\partial_\kappa^m\ip{\mathcal R_0(\kappa)u_1}{\phi}\right|+\norm{\partial_\kappa^m g_{u_1}(\kappa)}_{\mathcal V}\right)\le C_0C_1^m m!.
		\end{equation}
		Applying Cauchy's integral formula to the holomorphic families $g_{\overline\phi}$ and $T^{-1}$ from Lemma~\ref{lem:odd-compact} also gives
		\begin{equation}\label{eq:odd-internal-derivative-bounds}
			\sup_{|\kappa|\le2\eps}(\norm{\partial_\kappa^m g_{\overline\phi}(\kappa)}_{\mathcal V}+\norm{\partial_\kappa^m(T(\kappa)^{-1})}_{\mathcal L(\mathcal V)})\le C_0C_1^m m!.
		\end{equation}
		
		For $u_2$, the variables $x\in D_1$ and $y\notin D_2$ satisfy $r=|x-y|\ge r_0$ for some fixed $r_0>0$. We estimate the kernel using \eqref{eq:odd-kernel}. Expanding its fixed polynomial reduces the derivatives in $\kappa$ to terms of the form
		\[
		\partial_\kappa^m(\kappa^\ell e^{-\kappa r})=e^{-\kappa r}\sum_{k=0}^{\min(m,\ell)}\binom mk\frac{\ell!}{(\ell-k)!}\kappa^{\ell-k}(-r)^{m-k},
		\]
		where $\ell$ ranges over a fixed finite set, so the number of terms is bounded independently of $m$, and their coefficients grow at most as a fixed power of $m+1$. Since $\nabla_x r=(x-y)/r$, a spatial derivative of \eqref{eq:odd-kernel} produces terms of the same form, with different fixed powers of $r$ and with the degree in $\kappa$ increased by at most one. The lower bound $r\ge r_0$ controls all negative powers of $r$. Using also $r\le C(1+|y|)$, $|\kappa|\le2\eps$, and $|e^{-\kappa r}|\le1$ for $\Rea\kappa\ge0$, we obtain
		\begin{equation}\label{eq:odd-tail-kernel}
			|\partial_\kappa^mG_n(\kappa,r)|+|\nabla_x\partial_\kappa^mG_n(\kappa,r)|\le C_0C_1^m(1+|y|)^{m+d_0}
		\end{equation}
		for some fixed $d_0\ge0$. The constants are independent of $m$, $\kappa$, $x$, and $y$ in the indicated ranges; the fixed polynomial growth in $m$ has been absorbed into $C_1^m$.
		
		Integrating \eqref{eq:odd-tail-kernel} against $|u_2|$ and applying \eqref{eq:weighted-moments} with $d=d_0$, followed by pairing with $\phi$ or taking the $L^2(D)$-norm, gives
		\begin{equation}\label{eq:estimate_for_u_2}
			\sup_{\substack{|\kappa|\le2\eps\\\Rea\kappa\ge0}}\left(\left|\partial_\kappa^m\ip{\mathcal R_0(\kappa)u_2}{\phi}\right|+\norm{\partial_\kappa^m g_{u_2}(\kappa)}_{\mathcal V}\right)\le C_0C_1^m(m!)^2,
		\end{equation}
		where the fixed factors $\norm{e^{a\sqrt{|\cdot|}}u}_2$, $\norm{\phi}_1$, and $|D|^{1/2}$ are included in $C_0$, and the constants are independent of $m$ and $\kappa$.
		
		For each fixed $m$, \eqref{eq:odd-tail-kernel} provides an integrable majorant against $|u_2|$, uniformly in $\kappa$ and $x\in D_1$. This justifies differentiation under the integral in the open right half-plane. By the dominated convergence theorem, the derivatives of the scalar pairing and of $g_{u_2}$ extend continuously to $|\kappa|\le2\eps$, $\Rea\kappa\ge0$, including zero and the imaginary axis. Combining the estimates \eqref{eq:estimate_for_u_1} and \eqref{eq:estimate_for_u_2}, we obtain
		\[
		\sup_{\substack{|\kappa|\le2\eps\\\Rea\kappa\ge0}}\left(\left|\partial_\kappa^m\ip{\mathcal R_0(\kappa)u}{\phi}\right|+\norm{\partial_\kappa^m g_u(\kappa)}_{\mathcal V}\right)\le C_0C_1^m(m!)^2.
		\]
		By \eqref{eq:estimate_for_u_1}, \eqref{eq:estimate_for_u_2}, and \eqref{eq:odd-internal-derivative-bounds}, the derivatives of $\ip{\mathcal R_0(\kappa)u}{\phi}$, $g_u$, $g_{\overline\phi}$, and $T^{-1}$ all satisfy the bound $C_0C_1^m(m!)^2$ on $\overline{U_\eps}$ in their respective norms, after increasing $C_0,C_1$ and using $m!\le(m!)^2$. These functions are holomorphic in $U_\eps$, so \eqref{eq:scalar-reconstruction} also gives the holomorphy of $h$ there.
		
		Since $B$ is independent of $\kappa$, Leibniz' rule gives
		\begin{equation}\label{eq:scalar-reconstruction-derivatives}
			\begin{aligned}
				h^{(m)}(\kappa)=\partial_\kappa^m\ip{\mathcal R_0(\kappa)u}{\phi}-\sum_{\substack{m_1,m_2,m_3\ge0\\m_1+m_2+m_3=m}}\frac{m!}{m_1!m_2!m_3!}
				[B(\partial_\kappa^{m_1}(T(\kappa)^{-1}))\partial_\kappa^{m_2}g_u(\kappa),\partial_\kappa^{m_3}g_{\overline{\phi}}(\kappa)]_D.
			\end{aligned}
		\end{equation}
		The Cauchy--Schwarz inequality gives $|[BX,Y]_D|\le\norm{B}_\infty\norm{X}_{\mathcal V}\norm{Y}_{\mathcal V}$. For each triple with $m_1+m_2+m_3=m$, the factorial factor satisfies $\frac{m!}{m_1!m_2!m_3!}(m_1!)^2(m_2!)^2(m_3!)^2=m!m_1!m_2!m_3!\le(m!)^2$. Each summand is bounded by $\norm{B}_\infty C_0^3C_1^m(m!)^2$. Summing over the $\binom{m+2}{2}\le3^m$ triples gives
		\[
		|h^{(m)}(\kappa)|\le (C_0+\norm{B}_\infty C_0^3)(3C_1)^m(m!)^2,
		\quad \kappa\in U_\eps.
		\]
		Finally, the scalar term and every differentiated factor in the formula for $h^{(m)}$ extend continuously to $\overline{U_\eps}$. This formula defines the continuous extension of $h^{(m)}$, with the same bound. Increasing the constants gives \eqref{eq:Gevrey}.
	\end{proof}
	
	We next consider even dimensions, writing $n=2\nu$ and $d=\nu-1$. For $w\ne0$ with $\Rea w\ge0$, we use the principal logarithm, so that $|\Ima\log w|\le\pi/2$. Recall that $K_d$ is the modified Bessel function of the second kind of order $d$. Its integer-order expansion \cite[Eq.~(10.31.1)]{NIST10} is
	\begin{equation}\label{eq:integer-Bessel-series}
		K_d(w)=S_d(w)+\Big(\frac w2\Big)^d\sum_{m=0}^{\infty}\frac{(w^2/4)^m}{m!(m+d)!}\left(\alpha_d\log(w/2)+\beta_{d,m}\right),
	\end{equation}
	where $S_0=0$ and, for $d\ge1$,
	\begin{equation}\label{eq:Bessel-finite-part}
		S_d(w)=\frac12\Big(\frac w2\Big)^{-d}\sum_{m=0}^{d-1}\frac{(d-m-1)!}{m!}\Big(-\frac{w^2}{4}\Big)^m.
	\end{equation}
	Here $\alpha_d=(-1)^{d+1}$, $\beta_{d,m}=\frac{(-1)^d}{2}(\psi(m+1)+\psi(m+d+1))$, and $\psi=\Gamma'/\Gamma$ is the digamma function. To estimate $\beta_{d,m}$, the recurrence $\psi(z+1)=\psi(z)+1/z$ gives $\psi(k+1)=\psi(1)+H_k$, where $H_0=0$ and $H_k=\sum_{\ell=1}^k1/\ell$ for $k\ge1$. Since $H_k\le1+\log(k+1)$, we obtain
	\begin{equation}\label{eq:bound_beta_m}
		\begin{aligned}
			|\beta_{d,m}|\le|\psi(1)|+\frac{H_m+H_{m+d}}2\le|\psi(1)|+1+\log(m+d+1)\le C_d(1+\log(m+d+2)),
		\end{aligned}
	\end{equation}
	where $C_d$ depends only on $d$ and is independent of $m$.
	
	For each integer $J\ge1$, define the remainder
	\begin{equation}\label{eq:def-of-EJ,d}
		\begin{split}
			E_{J,d}(w)&=K_d(w)-S_d(w)-\Big(\frac w2\Big)^d\sum_{m=0}^{J-1}\frac{(w^2/4)^m}{m!(m+d)!}(\alpha_d\log(w/2)+\beta_{d,m})\\
			&=\Big(\frac w2\Big)^d\sum_{m=J}^{\infty}\frac{(w^2/4)^m}{m!(m+d)!}(\alpha_d\log(w/2)+\beta_{d,m}).
		\end{split}
	\end{equation}
	
	\begin{lemma}[Remainder estimate for the Bessel expansion]\label{lem:Bessel-remainder}
		There is $C_d>1$ such that, for every $J\ge1$ and every $w\ne0$ with $\Rea w\ge0$,
		\begin{equation}\label{eq:Bessel-remainder}
			|E_{J,d}(w)|+\frac{|w|}{J+d+1}|E'_{J,d}(w)|\le\frac{C_d^{J+1}|w|^{2J+d}}{J!(J+d)!}(1+|\log w|+\log(J+d+2)).
		\end{equation}
	\end{lemma}
	
	\begin{proof}
		Let $f_m(w)$ denote the $m$th term in \eqref{eq:integer-Bessel-series}, including the factor $(w/2)^d$. $K_d(w)=S_d(w)+\sum_{m=0}^\infty f_m(w)$ and $E_{J,d}(w)=\sum_{m=J}^\infty f_m(w)$. Set $s=|w|$,
		\[
		A_m(s)=\frac{s^{2m+d}}{4^m m!(m+d)!}\quad\text{and}\quad L_m(w)=1+|\log w|+\log(m+d+2).
		\]
		Using \eqref{eq:bound_beta_m}, $|\alpha_d|=1$, and $|\log(w/2)|\le|\log w|+\log2$, we obtain $|f_m(w)|\le C_dA_m(s)L_m(w)$. The product rule then gives
		\[
		\frac{s}{m+d+1}|f'_m(w)|\le\frac{2m+d}{m+d+1}|f_m(w)|+\frac{2^{-d}A_m(s)}{m+d+1}\le2|f_m(w)|+A_m(s).
		\]
		Since $L_m(w)\ge1$, combining these bounds and increasing $C_d$ yields
		\begin{equation}\label{eq:Bessel-term-bound}
			|f_m(w)|+\frac{s}{m+d+1}|f'_m(w)|\le C_dA_m(s)L_m(w),
		\end{equation}
		where the constant $C_d>0$ depends only on $d$.
		
		First, suppose that $s\le J+1$. For $m\ge J$, $\frac{A_{m+1}(s)}{A_m(s)}=\frac{s^2}{4(m+1)(m+d+1)}\le\frac14$, which implies  $A_{J+\ell}(s)\le4^{-\ell}A_J(s)$. We also have $L_{J+\ell}(w)\le(\ell+1)L_J(w)$ and $J+\ell+d+1\le(J+d+1)(\ell+1)$. Applying \eqref{eq:Bessel-term-bound} together with these inequalities to the series for the remainder and its derivative yields
		\[
		|E_{J,d}(w)|+\frac{s}{J+d+1}|E'_{J,d}(w)|\le C_dA_J(s)L_J(w)\sum_{\ell=0}^{\infty}4^{-\ell}(\ell+1)^2.
		\]
		The sum is finite and independent of $J$ and $w$. This proves \eqref{eq:Bessel-remainder} in this region.
		
		Next, suppose that $s\ge J+1$, and set $Q_{J,d}(s)=s^{2J+d}/(J!(J+d)!)$. We estimate separately the finite sum subtracted from $K_d$, the term $S_d$, and $K_d$ itself. For $0\le m<J$,
		\begin{equation}\label{eq:Bessel-power-comparison}
			\frac{A_m(s)}{Q_{J,d}(s)}=\frac{J!(J+d)!}{4^m m!(m+d)!s^{2(J-m)}}\le4^{-m}\Big(\frac{J+d}{J+1}\Big)^{2(J-m)}\le4^{-m}e^{2\max(d-1,0)}.
		\end{equation}
		The first inequality follows by writing the factorial quotient as a product of $J-m$ factors of the form $k(k+d)$ and bounding each factor by $(J+d)^2$. Since $L_m(w)\le L_J(w)$ and $m+d+1\le J+d+1$, estimate~\eqref{eq:Bessel-term-bound} gives
		\[
		\Big|\sum_{m=0}^{J-1}f_m(w)\Big|+\frac{s}{J+d+1}\Big|\sum_{m=0}^{J-1}f'_m(w)\Big|\le C_dJQ_{J,d}(s)L_J(w).
		\]
		
		The large-argument estimates \cite[Eqs.~(10.40.2) and (10.40.4)]{NIST10} are uniform for $|\arg w|\le\pi/2$. Enlarging the constant to include the bounded part of this sector gives
		\begin{equation}\label{eq:Bessel-large}
			|K_d(w)|+|K'_d(w)|\le C_ds^{-1/2},\quad s\ge1,\quad\Rea w\ge0.
		\end{equation}
		For the factorials, we use
		\begin{equation}\label{eq:Bessel-factorial-bound}
			J!(J+d)!=(J!)^2\prod_{\ell=1}^d(J+\ell)\le d!(J+1)^{2J+d},
		\end{equation}
		where $J!\le(J+1)^J$ and $J+\ell\le\ell(J+1)$; for $d=0$, the product equals one. Combining \eqref{eq:Bessel-factorial-bound} with $s\ge J+1$ gives
		\[
		s^{-1/2}+\frac{s^{1/2}}{J+d+1}\le C_dQ_{J,d}(s),\quad s^d\le C_d(J+1)^dQ_{J,d}(s).
		\]
		Estimate~\eqref{eq:Bessel-large} controls $K_d$ and its derivative with the factor $s/(J+d+1)$ required in \eqref{eq:Bessel-remainder}. Since $S_d$ is a finite sum of powers $w^{2m-d}$ with $0\le m<d$, and $S_0=0$, we have
		\[
		|S_d(w)|+\frac{s}{J+d+1}|S'_d(w)|\le C_ds^d\le C_d(J+1)^dQ_{J,d}(s).
		\]
		Combining the estimates for $K_d$, $S_d$, and the finite sum proves \eqref{eq:Bessel-remainder} for $s\ge J+1$, after increasing $C_d$ so that $1+J+(J+1)^d\le C_d^{J+1}$.
	\end{proof}
	
	For $u\in L^2_{a,1/2}$ and $\phi\in C_c^\infty(\R^n)$, set
	\begin{equation}\label{eq:def-a_u-g_u}
		a_u(\kappa)=\ip{\mathcal R_0(\kappa)u}{\phi}\quad \text{and} \quad g_u(\kappa)=\nabla\mathcal R_0(\kappa)u|_D.
	\end{equation} 
	For nonzero $\kappa$ on the imaginary axis, these quantities are defined by the boundary values of the free kernel from $\Rea\kappa>0$. Write $t=|\kappa|$, $\ell(t)=1+|\log t|$, and $\zeta=\log\kappa$, using the same logarithmic branch as above.
	
	\begin{lemma}[Source expansion in even dimensions]\label{lem:even-source-expansion}
		There are scalar coefficients $a_r,b_r$, vector coefficients $g_r,h_r\in\mathcal V$, and constants $C_0,C_1>0$ such that
		\begin{equation}\label{eq:even-coefficient-bounds}
			|a_r|+|b_r|+\norm{g_r}_{\mathcal V}+\norm{h_r}_{\mathcal V}\le C_0C_1^r(r!)^2,\quad r\in\Nzero.
		\end{equation}
		For all sufficiently small nonzero $\kappa$ in the closed right half-plane and every integer $M\ge\nu+1$,
		\begin{equation}\label{eq:even-source-expansions}
			\begin{aligned}
				a_u(\kappa)=\sum_{r=0}^{M-1}\kappa^{2r}(a_r+b_r\zeta)+E_{a,M}(\kappa),\quad
				g_u(\kappa)=\sum_{r=0}^{M-1}\kappa^{2r}(g_r+h_r\zeta)+E_{g,M}(\kappa),
			\end{aligned}
		\end{equation}
		where
		\begin{equation}\label{eq:even-source-remainders}
			|E_{a,M}(\kappa)|+\norm{E_{g,M}(\kappa)}_{\mathcal V}\le C_0C_1^M(M!)^2t^{2M}\ell(t).
		\end{equation}
		The coefficients are independent of $M$, and the constants are independent of $r$, $M$, and $\kappa$.
	\end{lemma}
	
	\begin{proof}
		We first define the coefficients and estimate their growth. The moment estimate \eqref{eq:weighted-moments}, together with $\log(2+|y|)\le C(1+|y|)$, gives
		\begin{equation}\label{eq:even-weighted-moments}
			\int_{\R^n}(1+|y|)^{k+q}(1+\log(2+|y|))|u(y)|\,dy\le C_0C_1^k(k!)^2,\quad k\in\Nzero,
		\end{equation}
		for every fixed $q\ge0$, where the constants may depend on $q$ and $\norm{e^{a\sqrt{|\cdot|}}u}_2$, but not on $k$.
		
		Recall from \eqref{eq:free-kernel} that $G_n(\kappa,\rho)=c_n\kappa^d\rho^{-d}K_d(\kappa\rho)$, where $c_n=(2\pi)^{-\nu}$. Multiplying \eqref{eq:integer-Bessel-series} by $c_n\kappa^d\rho^{-d}$ and splitting $\log(\kappa\rho/2)=\zeta+\log(\rho/2)$ determines two radial coefficient kernels at each degree $\kappa^{2r}$. We denote them by $\mathfrak a_r$ and $\mathfrak b_r$. For $0\le r<d$, define
		\begin{equation}\label{eq:even-coefficient-kernels-low}
			\mathfrak a_r(x)=\frac{c_n(-1)^r2^{d-2r-1}(d-r-1)!}{r!}|x|^{2r-2d},\quad \mathfrak b_r(x)=0,\quad x\ne0.
		\end{equation}
		For $r=m+d$ with $m\ge0$, define
		\begin{equation}\label{eq:even-coefficient-kernels-high}
			\begin{aligned}
				\mathfrak a_r(x)=\frac{c_n|x|^{2m}}{2^{d+2m}m!r!}\left(\alpha_d\log(|x|/2)+\beta_{d,m}\right),\quad
				\mathfrak b_r(x)=\frac{c_n\alpha_d|x|^{2m}}{2^{d+2m}m!r!}.
			\end{aligned}
		\end{equation}
		We then set
		\begin{equation}\label{eq:even-source-coefficient-definitions}
			\begin{aligned}
				a_r=\ip{\mathfrak a_r*u}{\phi},\quad b_r=\ip{\mathfrak b_r*u}{\phi},\quad
				g_r=\left.\nabla(\mathfrak a_r*u)\right|_D,\quad h_r=\left.\nabla(\mathfrak b_r*u)\right|_D.
			\end{aligned}
		\end{equation}
		We verify below that these convolutions and their first gradients belong to $L^2$ on bounded sets. The definitions also show that $b_r=0$ for $r<d$ and $h_r=0$ for $r\le d$, since $\mathfrak b_d$ is constant.
		
		Choose a ball $\Omega$ containing $\overline D\cup\supp\phi$. For each scalar or gradient convolution, split the integral into $\rho=|x-y|\le1$ and $\rho\ge1$. When $r<d$, there are only finitely many kernels to consider. Their strongest singularities are $\rho^{2-n}$ for the scalar kernels and $\rho^{1-n}$ for their gradients, both locally integrable.
		
		For $r=m+d$, \eqref{eq:even-coefficient-kernels-high} and \eqref{eq:bound_beta_m} give
		\begin{equation}\label{eq:even-coefficient-local-bound}
			\norm{\mathfrak a_r}_{L^1(B_1)}+\norm{\mathfrak b_r}_{L^1(B_1)}+\norm{\nabla\mathfrak a_r}_{L^1(B_1)}+\norm{\nabla\mathfrak b_r}_{L^1(B_1)}\le\frac{C_d(r+1)(1+\log(r+2))}{m!r!}.
		\end{equation}
		Indeed, on $0<\rho\le1$ we have $\rho^{2m}\le1$ and $\rho^{2m-1}\le\rho^{-1}$, while $\int_0^1\rho^{n-2}(1+|\log\rho|)\,d\rho<\infty$. Schur's estimate \eqref{eq:schur-estimate}, together with \eqref{eq:even-coefficient-local-bound}, controls the part with $\rho\le1$ in local $L^2$.
		
		For $\rho=|x|\ge1$, direct differentiation gives
		\begin{equation}\label{eq:even-coefficient-far-bound}
			|\mathfrak a_r(x)|+|\mathfrak b_r(x)|+|\nabla\mathfrak a_r(x)|+|\nabla\mathfrak b_r(x)|\le\frac{C_d(r+1)(1+\log(r+2))}{m!r!}\rho^{2m}(1+\log\rho).
		\end{equation}
		In the convolution integrals, $x\in\Omega$ and $\rho=|x-y|\ge1$, so $\rho\le C_\Omega(1+|y|)$ and $1+\log\rho\le C_\Omega(1+\log(2+|y|))$. Applying \eqref{eq:even-weighted-moments} with $k=2m$ bounds the distant part uniformly in $x\in\Omega$. These estimates prove local $L^2$ convergence of the convolutions and give, for $r\ge d$,
		\[
		|a_r|+|b_r|+\norm{g_r}_{\mathcal V}+\norm{h_r}_{\mathcal V}\le C_0C^{r+1}(r+1)(1+\log(r+2))\frac{((2r)!)^2}{(r-d)!r!}.
		\]
		The factorial quotient satisfies
		\begin{equation}\label{eq:even-factorial-comparison}
			\frac{((2r)!)^2}{(r-d)!r!}=\binom{2r}{r}^2(r!)^2\frac{r!}{(r-d)!}\le16^r(r+1)^d(r!)^2,\quad r\ge d.
		\end{equation}
		The remaining fixed powers of $r+1$ and the logarithmic factor can be included in an exponential constant. The finitely many coefficients with $r<d$ are included by increasing $C_0$. This proves \eqref{eq:even-coefficient-bounds}.
		
		We now estimate the remainder of the finite expansion. Let $M\ge\nu+1=d+2$, and set $J=M-d$, so that $J\ge2$. For $x\ne0$, define
		\begin{equation}\label{eq:even-remainder-kernel-definition}
			\mathcal E_M(\kappa,x)=c_n\kappa^d|x|^{-d}E_{J,d}(\kappa|x|).
		\end{equation}
		Combining \eqref{eq:def-of-EJ,d} with the coefficient definitions \eqref{eq:even-coefficient-kernels-low} and \eqref{eq:even-coefficient-kernels-high} gives the exact finite identity
		\begin{equation}\label{eq:even-finite-kernel-identity}
			G_n(\kappa,|x|)=\sum_{r=0}^{M-1}\kappa^{2r}(\mathfrak a_r(x)+\zeta\mathfrak b_r(x))+\mathcal E_M(\kappa,x).
		\end{equation}
		Integrating this identity against $u$ gives \eqref{eq:even-source-expansions}, with $E_{a,M}(\kappa)=\ip{\mathcal E_M(\kappa,\cdot)*u}{\phi}$ and $E_{g,M}(\kappa)=\nabla(\mathcal E_M(\kappa,\cdot)*u)|_D$, once the remainder integrals have been estimated.
		
		For $\rho=|x|>0$, estimate~\eqref{eq:Bessel-remainder} gives a bound for $\mathcal E_M$. To estimate its gradient, use $\partial_\rho(\rho^{-d}E_{J,d}(\kappa\rho))=\rho^{-d-1}(-dE_{J,d}(\kappa\rho)+\kappa\rho E'_{J,d}(\kappa\rho))$. Since $J+d=M$ and $|\log(\kappa\rho)|\le|\log t|+\pi/2+|\log\rho|$, we obtain
		\begin{equation}\label{eq:even-remainder-kernel-bound}
			|\mathcal E_M(\kappa,x)|+\frac{\rho}{M+1}|\nabla_x\mathcal E_M(\kappa,x)|\le\frac{C^{M+1}t^{2M}\rho^{2J}}{J!M!}\left(\ell(t)+|\log\rho|+\log(M+2)\right).
		\end{equation}
		This estimate holds throughout the closed right half-plane, with $\kappa\ne0$.
		
		Split each remainder integral into $\rho=|x-y|\le1$ and $\rho\ge1$. For the local part, the scalar and gradient majorants in \eqref{eq:even-remainder-kernel-bound} satisfy
		\begin{equation}\label{eq:even-local-remainder-integrals}
			\int_{|x|\le1}|x|^{2J-\sigma}(1+|\log|x||)\,dx=|\mathbb S^{n-1}|\Bigl(\frac1{2J+n-\sigma}+\frac1{(2J+n-\sigma)^2}\Bigr)\le C_n,\quad \sigma=0,1.
		\end{equation}
		The bound is independent of $J\ge2$. Applying Schur's estimate \eqref{eq:schur-estimate} to \eqref{eq:even-remainder-kernel-bound}, using \eqref{eq:even-local-remainder-integrals}, controls the scalar and gradient remainders in local $L^2$.
		
		For $x\in\Omega$ and $\rho=|x-y|\ge1$, we have $\rho\le C_\Omega(1+|y|)$. The gradient factor $\rho^{-1}$ is bounded by one in this region. Applying \eqref{eq:even-weighted-moments} with $k=2J$ controls both remainder integrals uniformly in the observation variable. Since $J\le M$, we may bound $((2J)!)^2$ by $((2M)!)^2$. The factors $M+1$ and $1+\log(M+2)$ can be included in $C^{M+1}$, and $\ell(t)\ge1$. Combining the two regions gives
		\[
		|E_{a,M}(\kappa)|+\norm{E_{g,M}(\kappa)}_{\mathcal V}\le C_0C^{M+1}\frac{((2M)!)^2}{(M-d)!M!}t^{2M}\ell(t).
		\]
		Estimate~\eqref{eq:even-factorial-comparison}, with $r=M$, bounds the quotient by $16^M(M+1)^d(M!)^2$. Including the fixed power $(M+1)^d$ in the exponential constant proves \eqref{eq:even-source-remainders}. All coefficient kernels were defined before choosing $M$, so the coefficients in \eqref{eq:even-source-expansions} are independent of the truncation order.
	\end{proof}
	
	\subsection{Vanishing power--logarithm coefficients}
	
	Continue to assume $n=2\nu$ with $\nu\in\N$, and $u\in L^2_{a,1/2}$. Given $\phi\in C_c^\infty(\R^n)$, we write $t=|\kappa|$, $\ell(t)=1+|\log t|$, and $\zeta=\log\kappa$, with the branch of the logarithm fixed in the previous subsection.
	
	By \eqref{eq:even-K},
	\[
	\begin{aligned}
		T(\kappa)=\mathcal T(\kappa^2,\kappa^{2\nu}\log\kappa),\quad
		\mathcal T(z,w)=T(0)-zK_1(z)B-2wK_2(z)B,
	\end{aligned}
	\]
	where $B=A-I$. The operator-valued function $\mathcal T$ is holomorphic near $(0,0)$. Since $\mathcal T(0,0)=T(0)$ is invertible by Lemma~\ref{lem:zero-inverse}, continuity of $\mathcal T$ in the operator norm allows us to choose $\delta>0$ such that
	\[
	\norm{T(0)^{-1}(\mathcal T(z,w)-T(0))}_{\mathcal L(\mathcal V)}<\frac12,
	\quad\text{for }|z|<\delta,\text{ and } |w|<\delta.
	\]
	Set $E(z,w)=T(0)^{-1}(\mathcal T(z,w)-T(0))$. Since $\mathcal T(z,w)=T(0)(I+E(z,w))$, the Neumann series gives $\mathcal T(z,w)^{-1}=\sum_{m=0}^{\infty}(-E(z,w))^mT(0)^{-1}$. Each summand is holomorphic in $(z,w)$ and has operator norm at most $2^{-m}\norm{T(0)^{-1}}_{\mathcal L(\mathcal V)}$. The series converges locally uniformly in the operator norm, so its sum is holomorphic.
	
	For a bounded operator $F$ on $\mathcal V$, let $F^{\mathrm t}$ denote its transpose with respect to the bilinear pairing, so that $[FX,Y]_D=[X,F^{\mathrm t}Y]_D$ for $X,Y\in\mathcal V$. The map $F\mapsto F^{\mathrm t}$ is bounded and complex linear, and preserves holomorphic dependence. Define
	\begin{equation}\label{eq:internal-S-definition}
		S_\phi(\kappa)=(T(\kappa)^{-1})^{\mathrm t}B^{\mathrm t}g_{\overline\phi}(\kappa).
	\end{equation}
	
	Since $\mathcal B_n(0,r)=\mathfrak c_0$ is independent of $r$, the series \eqref{eq:B-series} gives $\nabla_x\mathcal B_n(z,|x-y|)=z\mathcal H_n(z,x,y)$, where $\mathcal H_n(z,x,y)=2(x-y)\sum_{m=1}^{\infty}m\mathfrak c_mz^{m-1}|x-y|^{2m-2}$. By Lemma~\ref{lem:even-compact}, the localized gradient operator associated with $\mathcal B_n$ is holomorphic and vanishes at zero, so division by $z$ gives a holomorphic operator family after removing the singularity at zero. Taking the gradient in \eqref{eq:even-kernel}, using $z=\kappa^2$ and $\log z=2\log\kappa$, and integrating against $\overline\phi$ gives
	\[
	g_{\overline\phi}(\kappa)=\gamma_0(\kappa^2)+\kappa^{2\nu}\log\kappa\,\gamma_1(\kappa^2),
	\]
	where
	\[
	\begin{aligned}
		\gamma_0(z)(x)=\int_{\R^n}\nabla_x\mathcal A_n(z,|x-y|)\overline{\phi(y)}\,dy,\quad \gamma_1(z)(x)=2\int_{\R^n}\mathcal H_n(z,x,y)\overline{\phi(y)}\,dy
	\end{aligned}
	\]
	for $x\in D$. Both functions are holomorphic near zero with values in $\mathcal V$ by the preceding operator holomorphy.
	
	Define
	\begin{equation}\label{eq:internal-S-holomorphic-representation}
		\begin{aligned}
			\mathscr S_\phi(z,w)=(\mathcal T(z,w)^{-1})^{\mathrm t}B^{\mathrm t}(\gamma_0(z)+w\gamma_1(z)),\quad S_\phi(\kappa)=\mathscr S_\phi(\kappa^2,\kappa^{2\nu}\log\kappa).
		\end{aligned}
	\end{equation}
	The function $\mathscr S_\phi$ is holomorphic near $(0,0)$. Choose $0<R<1$ so that $\mathscr S_\phi$ is holomorphic on a neighborhood of the closed polydisk $\{|z|\le R,\ |w|\le R\}$. Applying Cauchy's integral formula successively in $z$ and $w$ gives the Taylor expansion and coefficient bounds
	\begin{equation}\label{eq:internal-Taylor-coefficients}
		\mathscr S_\phi(z,w)=\sum_{\ell=0}^{\infty}\sum_{m=0}^{\infty}z^\ell w^mS_{\ell,m},\quad
		\norm{S_{\ell,m}}_{\mathcal V}\le CR^{-\ell-m},
	\end{equation}
	where $C$ is independent of $\ell$ and $m$. For sufficiently small nonzero $\kappa$ in the closed right half-plane, the point $(\kappa^2,\kappa^{2\nu}\zeta)$ lies inside this polydisk, and the substituted Taylor series is absolutely convergent. We may group its terms according to $q=\ell+\nu m$ and obtain
	\begin{equation}\label{eq:internal-power-log-series}
		S_\phi(\kappa)=\sum_{q=0}^{\infty}\kappa^{2q}Q_q(\zeta),
		\quad
		Q_q(\zeta)=\sum_{m=0}^{\lfloor q/\nu\rfloor}S_{q-\nu m,m}\zeta^m.
	\end{equation}
	In particular, $\deg Q_q\le\lfloor q/\nu\rfloor$.
	
	If $\ell+\nu m=q$, then $\ell+m\le q$, and \eqref{eq:internal-Taylor-coefficients} gives $\norm{S_{\ell,m}}_{\mathcal V}\le CR^{-q}$ since $R\in(0,1)$. Moreover, $|\zeta|\le C\ell(t)$ in the closed right half-plane. Since $Q_q$ has at most $q+1$ terms, increasing the exponential constant absorbs their number and gives constants $C_0,C_1>0$ such that
	\begin{equation}\label{eq:internal-power-log-bounds}
		\begin{aligned}
			\norm{Q_q(\zeta)}_{\mathcal V}\le C_0C_1^q\ell(t)^q,\quad 
			\Big\|S_\phi(\kappa)-\sum_{q=0}^{M-1}\kappa^{2q}Q_q(\zeta)\Big\|_{\mathcal V}\le C_0(C_1t^2\ell(t))^M.
		\end{aligned}
	\end{equation}
	These estimates hold for all $q\in\Nzero$, all integers $M\ge1$, and all sufficiently small nonzero $\kappa$ in the closed right half-plane, with constants independent of $q,M,\kappa$. Indeed, the second estimate follows by summing the first estimate over $q\ge M$; the resulting geometric series has ratio $C_1t^2\ell(t)$, which is at most $1/2$ after restricting $t$. Taking $M=1$ shows that $S_\phi(\kappa)\to Q_0$ as $\kappa\to0$, and in particular $S_\phi$ is uniformly bounded near zero.
	
	By \eqref{eq:internal-S-definition}, the reconstruction formula \eqref{eq:scalar-reconstruction} reads
	\begin{equation}\label{eq:even-scalar-reconstruction}
		h_u(\kappa):=\ip{\mathcal R_A(\kappa)u}{\phi}=a_u(\kappa)-[g_u(\kappa),S_\phi(\kappa)]_D,
	\end{equation}
	where $a_u(\kappa)$ is given by \eqref{eq:def-a_u-g_u}. Let $a_r,b_r,g_r,h_r$ be the coefficients defined in \eqref{eq:even-source-coefficient-definitions} and appearing in the finite expansions \eqref{eq:even-source-expansions}. For each $r\in\Nzero$, define
	\begin{equation}\label{eq:power-log-polynomial-definition}
		P_r(\zeta)=a_r+b_r\zeta-\sum_{q=0}^{r}[g_{r-q}+h_{r-q}\zeta,Q_q(\zeta)]_D.
	\end{equation}
	Each $P_r$ is a polynomial. Its definition uses only the coefficients fixed in Lemma~\ref{lem:even-source-expansion} and \eqref{eq:internal-power-log-series}, so it is independent of the order of truncation.
	
	For an integer $M\ge\nu+1$, write the finite sums as
	\begin{equation}\label{eq:power-log-partial-sums}
		\begin{aligned}
			A_M(\kappa)=\sum_{r=0}^{M-1}\kappa^{2r}(a_r+b_r\zeta),\quad
			G_M(\kappa)=\sum_{r=0}^{M-1}\kappa^{2r}(g_r+h_r\zeta),\quad
			S_M(\kappa)=\sum_{q=0}^{M-1}\kappa^{2q}Q_q(\zeta).
		\end{aligned}
	\end{equation}
	When these sums are multiplied, the coefficient multiplying $\kappa^{2r}$ in $A_M-[G_M,S_M]_D$ is $P_r(\zeta)$ for $r<M$. Every other product term has total index at least $M$. For fixed $M$, these latter terms are $O_M(t^{2M}\ell(t)^M)$, since $\deg Q_q\le q\le M-1$. Also, $\norm{G_M}_{\mathcal V}\le C_M\ell(t)$ and $S_\phi$ is bounded. The source remainders in \eqref{eq:even-source-remainders} and the remainder in \eqref{eq:internal-power-log-bounds} give
	\begin{equation}\label{eq:h_u}
		h_u(\kappa)=\sum_{r=0}^{M-1}\kappa^{2r}P_r(\zeta)+O_M(t^{2M}\ell(t)^{M+1}).
	\end{equation}
	Here the constant in $O_M$ may depend on $M$, but is independent of $\kappa$. For nonzero $\kappa$ on the imaginary axis, we use the boundary values of $a_u$ and $g_u$ specified before Lemma~\ref{lem:even-source-expansion}, together with \eqref{eq:internal-power-log-series} and \eqref{eq:even-scalar-reconstruction}.
	
	\begin{lemma}\label{lem:flat-power-log}
		Let $h_u$ be defined by \eqref{eq:even-scalar-reconstruction}, and suppose that every polynomial $P_r$ defined by \eqref{eq:power-log-polynomial-definition} vanishes identically. Then there are constants $c,C,\eps>0$ such that
		\begin{equation}\label{eq:flat-power-log-bound}
			|h_u(\kappa)|\le Ce^{-c/|\kappa|},\quad
			0<|\kappa|<\eps,\quad\Rea\kappa\ge0.
		\end{equation}
	\end{lemma}
	
	\begin{proof}
		We use the estimates whose constants are independent of $M$, namely \eqref{eq:even-coefficient-bounds}, \eqref{eq:even-source-remainders}, and \eqref{eq:internal-power-log-bounds}. Choose $C_*\ge1$ at least as large as each of the constants denoted by $C_1$ in those estimates. Throughout the proof, $C$ denotes a constant independent of $\kappa$ and $M$, and we may increase it from line to line.
		
		Choose $\eta>0$ so that $C_*\eta^2\le1/4$, and set $M=\lfloor\eta/t\rfloor$. For sufficiently small $t$, we have $M\ge\nu+1$ and $\eta/(2t)\le M\le\eta/t$. Put $\alpha=C_*t^2M^2$, $\beta=C_*t^2\ell(t)$, then $\alpha\le1/4$, while $\beta/\alpha=\ell(t)/M^2\le (4t^2\ell(t))/\eta^2\to 0$ as $t\to 0$. After restricting $t$ further, we may assume $\beta/\alpha\le1/2$. In particular, $\beta\le1/8$.
		
		For $0\le r<M$, the elementary bound $r!\le M^r$, together with $|\zeta|\le C\ell(t)$, gives $\norm{\kappa^{2r}(g_r+h_r\zeta)}_{\mathcal V}\le C\ell(t)\alpha^r$. Likewise, \eqref{eq:internal-power-log-bounds} gives $\norm{\kappa^{2q}Q_q(\zeta)}_{\mathcal V}\le C\beta^q$. Using $M!\le M^M$ in \eqref{eq:even-source-remainders}, we also obtain
		\[
		|E_{a,M}|+\norm{E_{g,M}}_{\mathcal V}\le C\ell(t)\alpha^M.
		\]
		Summing the coefficient bounds and using \eqref{eq:internal-power-log-bounds}, we obtain
		\[
		\begin{aligned}
			\norm{G_M}_{\mathcal V}\le C\ell(t)\sum_{r=0}^{M-1}\alpha^r\le C\ell(t),\quad 
			\norm{S_\phi-S_M}_{\mathcal V}\le C\beta^M.
		\end{aligned}
		\]
		
		Substituting \eqref{eq:even-source-expansions} into \eqref{eq:even-scalar-reconstruction} yields the exact identity
		\begin{equation}\label{eq:finite-scalar-reconstruction}
			h_u-(A_M-[G_M,S_M]_D)=E_{a,M}-[E_{g,M},S_\phi]_D-[G_M,S_\phi-S_M]_D.
		\end{equation}
		The bounds above, the boundedness of $S_\phi$, and the Cauchy--Schwarz inequality for $[\,\cdot,\cdot\,]_D$ give
		\begin{equation}\label{eq:h_u-estimate_1}
			|h_u-(A_M-[G_M,S_M]_D)|\le C\ell(t)(\alpha^M+\beta^M)\le C\ell(t)\alpha^M.
		\end{equation}
		
		It remains to estimate the finite expression $A_M-[G_M,S_M]_D$. By \eqref{eq:power-log-partial-sums}, for each $0\le s<M$, the terms with total index $s$ in $[G_M,S_M]_D$ correspond exactly to $r=s-q$ with $0\le q\le s$. The coefficient multiplying $\kappa^{2s}$ in $A_M-[G_M,S_M]_D$ is
		\[
		a_s+b_s\zeta-\sum_{q=0}^{s}[g_{s-q}+h_{s-q}\zeta,Q_q(\zeta)]_D=P_s(\zeta)=0,
		\]
		where we utilized \eqref{eq:power-log-polynomial-definition} and the assumption $P_s\equiv0$. The remaining terms come from $-[G_M,S_M]_D$ and satisfy $r+q\ge M$. For each $0\le r<M$, this condition is equivalent to $q\ge M-r$, giving
		\[
		A_M-[G_M,S_M]_D=-\sum_{r=0}^{M-1}\sum_{q=M-r}^{M-1}\kappa^{2(r+q)}[g_r+h_r\zeta,Q_q(\zeta)]_D.
		\]
		The estimates for the individual terms give
		\begin{equation}\label{eq:h_u-estimate_2}
			\begin{aligned}
				|A_M-[G_M,S_M]_D|\le C\ell(t)\sum_{r=0}^{M-1}\sum_{q=M-r}^{M-1}\alpha^r\beta^q\le\frac{C\ell(t)\alpha^M}{1-\beta}\sum_{k=1}^{M}\Big(\frac{\beta}{\alpha}\Big)^k\le C\ell(t)\alpha^M.
			\end{aligned}
		\end{equation}
		In the second inequality, we extend the sum in $q$ to infinity and set $k=M-r$. The last inequality uses $\beta\le1/8$ and $\beta/\alpha\le1/2$, so its constant is independent of $M$.
		
		Combining the two estimates \eqref{eq:h_u-estimate_1}, \eqref{eq:h_u-estimate_2}, and using $\alpha\le1/4$ and $M\ge\eta/(2t)$, we conclude that
		\[
		|h_u(\kappa)|\le C\ell(t)\alpha^M\le C\ell(t)4^{-M}\le C\ell(t)e^{-(\eta\log4)/(2t)}.
		\]
		Since $\log\ell(t)=o(t^{-1})$ as $t\to0$, the factor $\ell(t)$ is absorbed by decreasing the constant in the exponent. This proves \eqref{eq:flat-power-log-bound}, with $c=\eta\log4/4$, after decreasing $\eps$ if necessary.
	\end{proof}

	\section{Arbitrary observation}\label{sec:arbitrary}
	
	We prove Theorem~\ref{thm:arbitrary} by separating the fractional orders in \eqref{eq:odd-separation}. In odd dimensions, the order assumption forces every Taylor coefficient of each resolvent pairing to vanish; the bound \eqref{eq:Gevrey} then permits the use of Lemma~\ref{lem:boundary-uniqueness}\ref{item_2_lem_bound_unique}. In even dimensions, every power--logarithm coefficient vanishes, and \eqref{eq:flat-power-log-bound} gives the boundary decay required in Lemma~\ref{lem:boundary-uniqueness}\ref{item_1_lem_bound_unique}.
	
	\subsection{Boundary uniqueness at the root-exponential rate}
	
	\begin{lemma}\label{lem:boundary-uniqueness}
		Let $U_\eps=\{\kappa\in\C:\Rea\kappa>0,\ |\kappa|<\eps\}$, and let $h$ be holomorphic in $U_\eps$ and continuous on $\overline{U_\eps}$.
		\begin{enumerate}[\rm(i)]
			\item\label{item_1_lem_bound_unique} Suppose that
			\begin{equation}\label{eq:boundary-exponential}
				|h(\mathsf{i}t)|\le Ce^{-c/|t|},\quad 0<|t|<\eps,
			\end{equation}
			for some $c,C>0$. Then $h=0$ in $U_\eps$.
			
			\item\label{item_2_lem_bound_unique} Suppose that every derivative extends continuously to $\overline{U_\eps}$, $h^{(m)}(0)=0$ for every $m\in\Nzero$, and
			\begin{equation}\label{eq:boundary-Gevrey}
				\sup_{\overline{U_\eps}}|h^{(m)}|\le C_0C_1^m(m!)^2,\quad m\in\Nzero,
			\end{equation}
			for some $C_0,C_1>0$ independent of $m$. Then $h=0$ in $U_\eps$.
		\end{enumerate}
	\end{lemma}
	
	\begin{proof}
		For \ref{item_1_lem_bound_unique}, suppose that $h(\kappa_*)\ne0$ at some point $\kappa_*\in U_\eps$, and choose $|\kappa_*|<\eps'<\eps$. Let $\omega_{\kappa_*}$ denote harmonic measure on $\partial U_{\eps'}$ evaluated at $\kappa_*$. It gives the mean value formula $H(\kappa_*)=\int_{\partial U_{\eps'}}H\,d\omega_{\kappa_*}$ for every harmonic function $H$ in $U_{\eps'}$ and continuous on $\overline{U_{\eps'}}$. On the segment $\{\mathsf{i}t:|t|\le\eps'/2\}$, its density with respect to $dt$ is continuous and bounded below by a positive constant.
		
		For $\eta>0$, the function $\log(|h|+\eta)$ is continuous on $\overline{U_{\eps'}}$ and subharmonic in $U_{\eps'}$. The maximum principle for subharmonic functions bounds it by the harmonic extension of its boundary values, giving
		\begin{equation}\label{eq:harmonic-measure}
			\log(|h(\kappa_*)|+\eta)\le\int_{\partial U_{\eps'}}\log(|h|+\eta)\,d\omega_{\kappa_*}.
		\end{equation}
		Put $K=\log(1+\sup_{\overline{U_{\eps'}}}|h|)$. We have $\log(|h(\kappa)|+\eta)\le K$ for every $\kappa\in\partial U_{\eps'}$ and $0<\eta\le1$. On the imaginary diameter, \eqref{eq:boundary-exponential} gives $\log|h(\mathsf{i}t)|\le\log C-c/|t|$. Since $|t|^{-1}$ is not integrable near $0$, the lower bound for the harmonic measure density implies that $\int_{\partial U_{\eps'}}\log|h|\,d\omega_{\kappa_*}=-\infty$. To pass to the limit in \eqref{eq:harmonic-measure}, note that the nonnegative functions $K-\log(|h|+\eta)$ increase to $K-\log|h|$ as $\eta\to 0$, and the latter has infinite integral. By the monotone convergence theorem, the right-hand side of \eqref{eq:harmonic-measure} tends to $-\infty$. The left-hand side tends to the finite number $\log|h(\kappa_*)|$, which is a contradiction.
		
		For \ref{item_2_lem_bound_unique}, we show that the derivative hypotheses imply \eqref{eq:boundary-exponential} on a smaller half-disk. The continuous extensions of the derivatives allow Taylor's formula with integral remainder along the segment from $0$ to $\mathsf{i}t$. Since $h^{(k)}(0)=0$ for every $k\in\Nzero$, we obtain, for every integer $m\ge1$ and $0<|t|<\eps$,
		\[
		\begin{aligned}
			|h(\mathsf{i}t)|=\Big|\frac{(\mathsf{i}t)^m}{(m-1)!}\int_0^1(1-\theta)^{m-1}h^{(m)}(\mathsf{i}\theta t)\,d\theta\Big|\le\frac{|t|^m}{m!}\sup_{\overline{U_\eps}}|h^{(m)}|\le C_0(C_1|t|)^m m!,
		\end{aligned}
		\]
		where we used \eqref{eq:boundary-Gevrey} in the last inequality. Choose $0<\delta<\min\{\eps,(4C_1)^{-1}\}$. For $0<|t|<\delta$, choose $m=\lfloor(2C_1|t|)^{-1}\rfloor$. Then $m\ge(4C_1|t|)^{-1}$ and $C_1|t|m\le1/2$. Using $m!\le m^m$, we obtain
		\[
		|h(\mathsf{i}t)|\le C_0(C_1|t|m)^m\le C_0 2^{-m}\le C_0e^{-\log2/(4C_1|t|)},
		\]
		which shows \eqref{eq:boundary-exponential} with $\eps$ replaced by $\delta$. Applying \ref{item_1_lem_bound_unique} on $U_\delta$ gives $h=0$ there, and the identity theorem gives $h=0$ throughout $U_\eps$.
	\end{proof}
	
	\subsection{Separation of the orders}
	
	\begin{proof}[Proof of Theorem~\ref{thm:arbitrary}]
		First, assume that there exists $a>0$ such that $u_k\in L^2_{a,1/2}$ for every $k=1,\ldots,N$. Fix $\phi\in C_c^\infty(\cO)$ and set
		\[
		c_k=b_k\sin(\pi s_k)\ne0,\quad h_k(\kappa)=\ip{\mathcal R_A(\kappa)u_k}{\phi}.
		\]
		Proposition~\ref{prop:resolvent}, with the resolvent parameter written as $\lambda=\kappa^2$, gives
		\begin{equation}\label{eq:odd-separation}
			\sum_{k=1}^N c_k\kappa^{2s_k}h_k(\kappa)=0,\quad \kappa>0.
		\end{equation}
		
		Suppose first that $n$ is odd. By Proposition~\ref{prop:Gevrey}, all derivatives of $h_1,\ldots,h_N$ extend continuously to $\overline{U_\eps}$ and satisfy \eqref{eq:Gevrey}, with the same $\eps,C_0,C_1>0$ for every $k=1,\ldots,N$. We claim that
		\begin{equation}\label{eq:h_k^m=0}
			h_k^{(m)}(0)=0,\quad k=1,\ldots,N,\quad m\in\Nzero.
		\end{equation}
		Suppose that \eqref{eq:h_k^m=0} fails. For each $k$ whose derivatives at zero do not all vanish, let $m_k$ be the least integer such that $h_k^{(m_k)}(0)\ne0$. Taylor's formula along the positive real axis gives
		\[
		h_k(\kappa)=a_k\kappa^{m_k}+o(\kappa^{m_k}),\quad a_k=\frac{h_k^{(m_k)}(0)}{m_k!}\ne0,\quad \kappa\to0^+.
		\]
		For a function whose derivatives at zero all vanish, Taylor's formula instead gives $h_k(\kappa)=o(\kappa^K)$ for every fixed $K\ge0$. We call such a function flat at zero.
		
		Among the functions that are not flat, the exponents $2s_k+m_k$ are pairwise distinct, because $s_k-s_\ell\notin\frac12\Z$. Let $k_*$ be the index with the smallest exponent. Divide \eqref{eq:odd-separation} by $\kappa^{2s_{k_*}+m_{k_*}}$ and let $\kappa\to0^+$. Every other term tends to zero: the other nonflat functions have larger exponents, and the flat functions are smaller than every fixed power. We obtain $c_{k_*}a_{k_*}=0$, contradicting the choice of $k_*$. Hence, \eqref{eq:h_k^m=0} holds, and every $h_k$ is flat at zero. Applying part~\ref{item_2_lem_bound_unique} of Lemma~\ref{lem:boundary-uniqueness} with the derivative bound \eqref{eq:Gevrey}, we obtain $h_k=0$ in $U_\eps$ for $k=1,\ldots,N$.
		
		Suppose next that $n=2\nu$ is even. Denote the power--logarithm coefficients of $h_k$ in \eqref{eq:h_u} by $P_{k,r}$. For every integer $M\ge\nu+1$, that expansion gives
		\[
		h_k(\kappa)=\sum_{r=0}^{M-1}\kappa^{2r}P_{k,r}(\log\kappa)+O_M(\kappa^{2M}(1+|\log\kappa|)^{M+1}),\quad \kappa\to0^+.
		\]
		The polynomials $P_{k,r}$ are independent of $M$, and the remainder constant is independent of $\kappa$.
		
		If some coefficient of $h_k$ is nonzero, let $m_k$ be its least index and put $P_k=P_{k,m_k}$. Choosing a finite expansion beyond this index gives
		\begin{equation}\label{eq:h_k_leading}
			h_k(\kappa)=\kappa^{2m_k}P_k(\log\kappa)+O(\kappa^{2m_k+2}(1+|\log\kappa|)^{d_k})
		\end{equation}
		for a fixed integer $d_k\ge0$. Let $\ell_k=\deg P_k$ and let $\alpha_k\ne0$ be its leading coefficient. Then $P_k(\log\kappa)/(\log\kappa)^{\ell_k}\to\alpha_k$ as $\kappa\to0$. In particular, $|P_k(\log\kappa)|\ge\frac{|\alpha_k|}{2}|\log\kappa|^{\ell_k}$ for all sufficiently small $\kappa>0$. Using \eqref{eq:h_k_leading}, we obtain 
		\[
		\Big|\frac{h_k(\kappa)-\kappa^{2m_k}P_k(\log\kappa)}{\kappa^{2m_k}P_k(\log\kappa)}\Big|
		\le C\kappa^2\frac{(1+|\log\kappa|)^{d_k}}{|\log\kappa|^{\ell_k}}
		\longrightarrow0,\quad\kappa\to0.
		\]
		If all coefficients of $h_k$ vanish, the finite expansions give $h_k(\kappa)=o(\kappa^K)$ for every fixed $K\ge0$.
		
		Among the functions with a nonzero coefficient, the exponents $2(s_k+m_k)$ are pairwise distinct, because $s_k-s_\ell\notin\Z$. Let $k_*$ have the smallest exponent. Divide \eqref{eq:odd-separation} by $\kappa^{2(s_{k_*}+m_{k_*})}P_{k_*}(\log\kappa)$ and let $\kappa\to0^+$. By \eqref{eq:h_k_leading}, the $k_*$ term tends to $c_{k_*}$. Every other function with a nonzero coefficient has a larger exponent, and a positive power of $\kappa$ tends to zero faster than any fixed power of $|\log\kappa|$. Functions whose coefficients all vanish are $o(\kappa^K)$ for every fixed $K$, so their terms also tend to zero after division. This contradicts $c_{k_*}\ne0$. We have proved
		\begin{equation}\label{eq:even-separated-coefficients}
			P_{k,r}\equiv0,\quad k=1,\ldots,N,\quad r\in\Nzero.
		\end{equation}
		
		Applying Lemma~\ref{lem:flat-power-log} to $h_1,\ldots,h_N$, we obtain constants $\eps,c,C>0$, independent of $k$, such that
		\begin{equation}\label{eq:estimate_h_k_kappa}
			|h_k(\kappa)|\le Ce^{-c/|\kappa|},\quad k=1,\ldots,N,\quad 0<|\kappa|<\eps,\quad \Rea\kappa\ge0.
		\end{equation}
		After decreasing $\eps$ if necessary, the kernel formula \eqref{eq:free-kernel}, the estimates in Lemma~\ref{lem:even-compact} and \eqref{eq:Bessel-large} give continuous extensions of $a_{u_k}$ and $g_{u_k}$ to $\overline{U_\eps}\setminus\{0\}$ for $k=1,\ldots,N$. Since $S_\phi$ is also continuous there by \eqref{eq:internal-power-log-series}, the reconstruction formula \eqref{eq:even-scalar-reconstruction} gives a continuous extension of each $h_k$ to $\overline{U_\eps}\setminus\{0\}$. Estimate~\eqref{eq:estimate_h_k_kappa} extends $h_k$ continuously to zero with $h_k(0)=0$. Restricting the estimate \eqref{eq:estimate_h_k_kappa} to $\kappa=\mathsf{i}t$ gives \eqref{eq:boundary-exponential}. Since $h_k$ is holomorphic in $U_\eps$, part~\ref{item_1_lem_bound_unique} of Lemma~\ref{lem:boundary-uniqueness} gives $h_k=0$ in $U_\eps$ for every $k=1,\ldots,N$.
		
		For the fixed $\phi\in C_c^\infty(\cO)$, we have proved that $h_k(\kappa)=\ip{\mathcal R_A(\kappa)u_k}{\phi}=0$ for $\kappa\in U_\eps$ and $k=1,\ldots,N$. Since $L$ is nonnegative and self-adjoint, $\sigma(L)\subset[0,\infty)$. For $\Rea\kappa>0$, we have $-\kappa^2\notin[0,\infty)$. The holomorphy of the resolvent outside $\sigma(L)$ implies that $\kappa\mapsto\mathcal R_A(\kappa)$ is holomorphic in the right half-plane as a bounded operator on $L^2(\R^n)$. Since $u_k,\phi\in L^2(\R^n)$, each $h_k$ is holomorphic there. The right half-plane is connected, so the identity theorem extends the vanishing of $h_k$ from $U_\eps$ to the entire right half-plane. For every $\lambda>0$, taking $\kappa=\sqrt\lambda$ and using $\mathcal R_A(\sqrt\lambda)=R_\lambda$ gives
		\begin{equation}\label{eq:continued-resolvent-pairings}
			\ip{R_\lambda u_k}{\phi}=0,\quad \lambda>0,\quad k=1,\ldots,N.
		\end{equation}
		
		The argument applies to every $\phi\in C_c^\infty(\cO)$, so $R_\lambda u_k=0$ on $\cO$ in the distributional sense, and almost everywhere since $R_\lambda u_k\in L^2(\R^n)$. This holds for every $\lambda>0$. Since $u_k\in L^2_{a,1/2}\subset L^2(\R^n)$, Lemma~\ref{lem:individual} gives $u_k=0$ on $\R^n$ for $k=1,\ldots,N$.
		
		Finally, let $u_k\in\mathscr D_{1/2}$ for $k=1,\ldots,N$. Choose $r\ge0$ such that $u_1,\ldots,u_N\in H^{-r}(\R^n)$. Choose a nonempty ball $\cO_0\Subset\cO$ and apply the spectral smoothing operators $S_\delta$ from \eqref{eq:spectral-smoothing}. For every sufficiently small $\delta>0$, equation~\eqref{eq:smoothed-entanglement}, with $\eps=\delta$, shows that $S_\delta u_1,\ldots,S_\delta u_N$ satisfy \eqref{eq:entanglement} on $\cO_0$. For each such $\delta$, Lemma~\ref{lem:preserve-decay} gives $a_\delta>0$, independent of $k$, such that $S_\delta u_k\in L^2_{a_\delta,1/2}$ for $k=1,\ldots,N$. Applying the weighted $L^2$ case proved above to $S_\delta u_1,\ldots,S_\delta u_N$ with observation set $\cO_0$, we obtain $S_\delta u_k=0$ on $\R^n$ for $k=1,\ldots,N$. By Lemma~\ref{lem:smoothing}, $S_\delta u_k\to u_k$ in $H^{-r}(\R^n)$ as $\delta\to0$ for every $k$. This proves $u_k=0$ for $k=1,\ldots,N$ and completes the proof.
	\end{proof}
	
	\subsection{A logarithmic improvement in odd dimensions}
	
	We now prove Proposition~\ref{prop:odd-logarithmic}. Let $w_p$ be given by \eqref{eq:def_w_p_r}, with $p\ge0$, and define
	\begin{equation}\label{eq:def_M_m}
		M_0=1,\quad M_m=m!\prod_{\ell=1}^m\log^{2p}(e+\ell),\quad m\in \N.
	\end{equation}
	
	\begin{lemma}\label{lem:logarithmic-moments}
		For $a>0$ and fixed $d,p\ge0$, there are constants $C_0,C_1>0$, depending only on $n,a,d,p$, such that
		\begin{equation}\label{eq:logarithmic-moments}
			\int_{\R^n}(1+|y|)^{m+d}|u(y)|\,dy\le C_0C_1^m m!M_m\norm{e^{aw_p(|\cdot|)}u}_2
		\end{equation}
		for every $m\in\Nzero$ and every $u$ for which the weighted norm is finite.
	\end{lemma}
	
	\begin{proof}
		We first estimate the polynomial factor against the weight. For large $r$, put $q=w_p(r)$. Since $p$ is fixed and $\log\log(e+r)/\log r\to0$ as $r\to\infty$, the identity $\log q=\frac12\log r-p\log\log(e+r)$ gives $\log q/\log r\to1/2$. In particular, $q\to\infty$ and $\log(e+q)/\log(e+r)\to1/2$ as $r\to\infty$, so these two logarithms are comparable for sufficiently large $r$. Using $1+r\le2r$ for $r\ge1$ and the comparison of the logarithms above, we obtain
		\[
		1+r\le2r=2q^2\log^{2p}(e+r)\le Cq^2\log^{2p}(e+q)
		\]
		for all sufficiently large $r$, where $C>0$ is independent of $r$. The bounded range of $r$ contributes only a constant raised to $m$ in the estimates below.
		
		For $m\ge1$, write $q=m\tau$. The inequality $\log(e+m\tau)\le\log(e+m)(1+\log_+\tau)$, where $\log_+\tau=\max(0,\log\tau)$, gives
		\[
		q^{2m}\log^{2pm}(e+q)e^{-aq/2}\le m^{2m}\log^{2pm}(e+m)[\tau^2(1+\log_+\tau)^{2p}e^{-a\tau/2}]^m.
		\]
		The expression in square brackets has a finite supremum for $\tau>0$, depending only on $a,p$. Combining this with the preceding estimate and $m^m\le e^m m!$, we obtain
		\begin{equation}\label{eq:logarithmic-supremum}
			\sup_{r\ge0}(1+r)^m e^{-aw_p(r)/2}\le C_0C_1^m(m!)^2\log^{2pm}(e+m).
		\end{equation}
		
		To express the last factor in terms of $M_m$, note that $x\mapsto\log(e+x)/x$ is decreasing for $x>0$. For $1\le\ell\le m$, we have $\log(e+\ell)\ge\frac{\ell}{m}\log(e+m)$, which gives
		\begin{equation}\label{eq:logarithmic-product}
			\prod_{\ell=1}^m\log(e+\ell)\ge\frac{m!}{m^m}\log^m(e+m)\ge e^{-m}\log^m(e+m),\quad m\ge1.
		\end{equation}
		Raising \eqref{eq:logarithmic-product} to the power $2p$ and using \eqref{eq:def_M_m} to bound the logarithmic factor in \eqref{eq:logarithmic-supremum}, we obtain
		\begin{equation}\label{eq:logarithmic-supremum-sequence}
			\sup_{r\ge0}(1+r)^m e^{-aw_p(r)/2}\le C_0C_1^m m!M_m,
		\end{equation}
		after increasing $C_1$.
		
		Finally, we can estimate the left-hand side of \eqref{eq:logarithmic-moments} via the Cauchy--Schwarz inequality
		\[
		\int_{\R^n}(1+|y|)^{m+d}|u(y)|\,dy \le \norm{e^{aw_p(|\cdot|)}u}_2\Big(\int_{\R^n}(1+|y|)^{2m+2d}e^{-2aw_p(|y|)}\,dy\Big)^{1/2}.
		\]
		Bounding $(1+|y|)^m e^{-aw_p(|y|)/2}$ by \eqref{eq:logarithmic-supremum-sequence} leaves the integral $\int_{\R^n}(1+|y|)^{2d}e^{-aw_p(|y|)}\,dy$, which is finite and independent of $m$, since $w_p(r)/\log r\to\infty$ as $r\to\infty$. This proves \eqref{eq:logarithmic-moments}, with $m=0$ included by increasing $C_0$.
	\end{proof}
	
	We now recall the quasianalyticity criterion for holomorphic functions with bounds on all derivatives. Let $(M_m)$ be a logarithmically convex sequence with $M_0=1$ and $M_{m+1}/M_m\to\infty$. On a bounded sector of opening $\pi$, the class defined by
	\[
	\sup|h^{(m)}|\le CA^m m!M_m,\quad m\in\Nzero,
	\]
	is quasianalytic precisely when
	\begin{equation}\label{eq:sector-quasianalyticity}
		\sum_{m=0}^\infty\Big(\frac{M_m}{(m+1)M_{m+1}}\Big)^{1/2}=\infty.
	\end{equation}
	Here quasianalyticity means that $h^{(m)}(0)=0$ for every $m$ implies $h=0$. We use \cite[Theorem~2.15 and Remark~2.18]{JGS20}; the remark gives the criterion for bounded sectors of opening at most $\pi$, including the right half-disk $U_\eps$.
	
	\begin{proof}[Proof of Proposition~\ref{prop:odd-logarithmic}]
		First, assume that there exists $a>0$ such that $e^{aw_p(|x|)}u_k\in L^2(\R^n)$ for every $k=1,\ldots,N$. Fix $\phi\in C_c^\infty(\cO)$ and set $h_k(\kappa)=\ip{\mathcal R_A(\kappa)u_k}{\phi}$. As in the proof of Proposition~\ref{prop:Gevrey}, choose balls $D_1\Subset D_2$ such that $\overline D\cup\supp\phi\subset D_1$, and write $u_k=u_{k,1}+u_{k,2}$, where $u_{k,1}=u_k\mathbf1_{D_2}$ and $u_{k,2}=u_k\mathbf1_{\R^n\setminus D_2}$. Lemma~\ref{lem:odd-compact} gives holomorphic extensions of $\ip{\mathcal R_0(\kappa)u_{k,1}}{\phi}$ and $g_{u_{k,1}}(\kappa)$ to $|\kappa|<4\eps$ for sufficiently small $\eps>0$. Applying the Cauchy estimate \eqref{eq:estimate_for_u_1} with $u_1=u_{k,1}$ bounds their derivatives by $C_0C_1^m m!$ in the scalar and $\mathcal V$ norms, respectively, for $|\kappa|\le2\eps$.
		
		For $u_{k,2}$, the variables $x\in D_1$ and $y\notin D_2$ satisfy $|x-y|\ge r_0>0$. Integrating \eqref{eq:odd-tail-kernel} against $|u_{k,2}(y)|$ reduces the estimate to the moment $\int_{\R^n}(1+|y|)^{m+d_0}|u_{k,2}(y)|\,dy$. By \eqref{eq:logarithmic-moments} with $d=d_0$, this moment is bounded by $C_0C_1^m m!M_m\norm{e^{aw_p(|\cdot|)}u_k}_2$. Pairing the scalar kernel integral with $\phi$ and taking the $L^2(D)$-norm of the gradient kernel integral introduce only the fixed factors $\norm{\phi}_1$ and $|D|^{1/2}$. Combining these bounds with the estimate for $u_{k,1}$ and using $M_m\ge1$, we obtain
		\[
		\sup_{\overline{U_\eps}}\left(\left|\partial_\kappa^m\ip{\mathcal R_0(\kappa)u_k}{\phi}\right|+\norm{\partial_\kappa^m g_{u_k}(\kappa)}_{\mathcal V}\right)\le C_0C_1^m m!M_m
		\]
		for every $m\in\Nzero$ and $k=1,\ldots,N$. Here $C_0$ includes the maximum of the weighted norms of $u_1,\ldots,u_N$, and $C_0,C_1$ are independent of $m,\kappa,k$. For each fixed $m$, the kernel bound is integrable against $|u_{k,2}|$, uniformly for $x\in D_1$ and $\kappa\in\overline{U_\eps}$. It justifies differentiation under the integral in $U_\eps$ and, by dominated convergence, the continuous extension of the derivatives to $\overline{U_\eps}$. The derivatives associated with $u_{k,1}$ are already continuous there by holomorphy.
		
		The factors $g_{\overline\phi}$ and $T^{-1}$ satisfy \eqref{eq:odd-internal-derivative-bounds}. Since $M_m\ge1$, their derivatives are also bounded by $C_0C_1^m m!M_m$ in the respective norms. By the definition \eqref{eq:def_M_m} of $M_m$, we have
		\begin{equation}\label{eq:M_m+1-divide-M_m}
			\frac{M_{m+1}}{M_m}=(m+1)\log^{2p}(e+m+1).
		\end{equation} 
		These quotients are increasing and tend to infinity, so the sequence is increasing and logarithmically convex. Logarithmic convexity and $M_0=1$ give $M_rM_\ell\le M_{r+\ell}$. For $m_1+m_2+m_3=m$, the factors in \eqref{eq:scalar-reconstruction-derivatives} satisfy $\frac{m!}{m_1!m_2!m_3!}(m_1!M_{m_1})(m_2!M_{m_2})(m_3!M_{m_3})=m!M_{m_1}M_{m_2}M_{m_3}\le m!M_m$. There are $\binom{m+2}{2}\le3^m$ such triples. Applying \eqref{eq:scalar-reconstruction-derivatives} with $u=u_k$ and including this number in $C_1^m$ gives
		\begin{equation}\label{eq:logarithmic-resolvent-derivatives}
			\sup_{\overline{U_\eps}}|h_k^{(m)}(\kappa)|\le C_0C_1^m m!M_m,\quad m\in\Nzero.
		\end{equation}
		
		Proposition~\ref{prop:resolvent} gives \eqref{eq:odd-separation}. The minimum-exponent argument proving \eqref{eq:h_k^m=0} uses only this identity, the existence of all derivatives at zero, and $s_k-s_\ell\notin\frac12\Z$. It applies to $h_1,\ldots,h_N$ and gives \eqref{eq:h_k^m=0} for these pairings.
		
		By \eqref{eq:M_m+1-divide-M_m},
		\begin{equation}\label{eq:logarithmic-quasianalytic-series}
			\sum_{m=0}^{\infty}\Big(\frac{M_m}{(m+1)M_{m+1}}\Big)^{1/2}=\sum_{m=0}^{\infty}\frac{1}{(m+1)\log^p(e+m+1)}=\infty,\quad 0\le p\le1.
		\end{equation}
		The criterion \eqref{eq:sector-quasianalyticity}, applied to \eqref{eq:logarithmic-resolvent-derivatives} and \eqref{eq:h_k^m=0}, gives $h_k=0$ in $U_\eps$. The analytic continuation argument leading to \eqref{eq:continued-resolvent-pairings} applies to these $L^2$ inputs. The arbitrariness of $\phi$ gives $R_\lambda u_k=0$ on $\cO$ for every $\lambda>0$, and Lemma~\ref{lem:individual} gives $u_k=0$.
		
		Finally, let $u_k\in\mathscr D_{\log,p}$ for $k=1,\ldots,N$. Since $N$ is finite, the definition of $\mathscr D_{\log,p}$ gives $r\ge0$ and $a,C>0$ such that $u_k\in H^{-r}(\R^n)$ and $|\ip{u_k}{\psi}|\le Ce^{-aw_p(R)}\norm{\psi}_{H^r}$ for every $k=1,\ldots,N$, $R\ge1$, and $\psi\in C_c^\infty(\R^n\setminus\overline{B_R})$.
		
		Fix $\delta>0$, and let $k_\delta(x,y)$ denote the smooth kernel of $S_\delta$ from \eqref{eq:spectral-smoothing}. By Lemma~\ref{lem:smoothing}, the function $k_\delta(x,\cdot)$ is supported in $\overline{B_{c_L\delta}(x)}$. Put $c_\delta=c_L\delta+1$ and $R_x=|x|-c_\delta$. For sufficiently large $|x|$, we have $R_x\ge1$, and every point $y$ in this support satisfies $|y|\ge R_x+1$. We may use $k_\delta(x,\cdot)$ as a test function in the decay estimate above. Taking $\alpha=0$ in \eqref{eq:kernel-row} bounds its $H^r$-norm independently of $x$, and we obtain
		\[
		|S_\delta u_k(x)|=|\ip{u_k}{k_\delta(x,\cdot)}|\le Ce^{-aw_p(R_x)}\norm{k_\delta(x,\cdot)}_{H^r}\le C_\delta e^{-aw_p(R_x)}
		\]
		for all sufficiently large $|x|$ and every $k=1,\ldots,N$, where $C_\delta$ is independent of $x$ and $k$.
		
		For sufficiently large $|x|$, we have $R_x\ge|x|/2$. Since $w_p$ is eventually increasing and $w_p(r/2)\ge w_p(r)/\sqrt2$, the bound for $S_\delta u_k$ gives, for $0<a'<a/\sqrt2$,
		\begin{equation}\label{eq:logarithmic-smoothed-weight}
			e^{a'w_p(|x|)}|S_\delta u_k(x)|\le C_\delta e^{-(a/\sqrt2-a')w_p(|x|)}
		\end{equation}
		for all sufficiently large $|x|$. The right-hand side is square integrable outside a ball because $w_p(r)/\log r\to\infty$ as $r\to\infty$. On bounded sets, integrability follows from the smoothness of $S_\delta u_k$. We obtain $e^{a'w_p(|x|)}S_\delta u_k\in L^2(\R^n)$ for every $k=1,\ldots,N$.
		
		Choose a nonempty ball $\cO_0\Subset\cO$. By \eqref{eq:smoothed-entanglement}, with $\eps=\delta$, the functions $S_\delta u_1,\ldots,S_\delta u_N$ satisfy \eqref{eq:entanglement} on $\cO_0$ for all sufficiently small $\delta>0$. The weighted $L^2$ case proved above, applied with observation set $\cO_0$, gives $S_\delta u_k=0$ on $\R^n$ for $k=1,\ldots,N$. By Lemma~\ref{lem:smoothing} and \eqref{eq:spectral-sobolev-equivalence}, $S_\delta u_k\to u_k$ in $H^{-r}(\R^n)$ as $\delta\to0$, so $u_k=0$ for every $k$.
	\end{proof}

	\subsection{A logarithmic improvement in even dimensions}
	
	We prove Proposition~\ref{prop:even-logarithmic} for the full range $0\le p\le1$. We write $n=2\nu$ and $d=\nu-1$. We first strengthen the Bessel remainder estimate to include higher derivatives.
	
	\begin{lemma}\label{lem:Bessel-remainder-derivatives}
		Let $E_{J,d}$ be the remainder defined by \eqref{eq:def-of-EJ,d}. For $J\ge2$, put $M=J+d$ and $F_M(w)=w^dE_{J,d}(w)$. There is a constant $C_d>1$, independent of $J,k,w$, such that
		\begin{equation}\label{eq:Bessel-remainder-derivatives}
			|F_M^{(k)}(w)|\le\frac{C_d^{M+k+1}k!|w|^{2M-k}}{J!M!}(1+|\log w|+\log(M+2))
		\end{equation}
		for $0\le k\le2J-2$ and every nonzero $w$ with $\Rea w\ge0$.
	\end{lemma}
	
	\begin{proof}
		Write $s=|w|$ and $L_M(w)=1+|\log w|+\log(M+2)$. Throughout the proof, $C_d>1$ depends only on $d$ and may be increased. All logarithms use the principal branch. Since $\Rea w\ge0$, the closed disk $|z-w|\le s/4$ avoids $(-\infty,0]$. On its boundary, $3s/4\le|z|\le5s/4$ and $|\log z-\log w|\le\log(4/3)$. Set $T_r(z)=\frac{z^{2r+2d}}{2^d4^r r!(r+d)!}$, then \eqref{eq:def-of-EJ,d} gives $F_M(z)=\sum_{r=J}^{\infty}T_r(z)(\alpha_d\log(z/2)+\beta_{d,r})$.
		First, suppose that $s\le J+1$. On $|z-w|=s/4$, the factors $T_r$ satisfy
		\[
		\frac{|T_{r+1}(z)|}{|T_r(z)|}
		=\frac{|z|^2}{4(r+1)(r+d+1)}
		\le\frac{25}{64}=:\rho,\quad r\ge J.
		\]
		Thus, $|T_{J+\ell}(z)|\le\rho^\ell|T_J(z)|$ for $\ell\ge0$. Moreover, \eqref{eq:bound_beta_m} and $\log(M+\ell+2)\le\log(M+2)+\log(\ell+1)$ give
		\[
		|\alpha_d\log(z/2)+\beta_{d,J+\ell}|\le C_d(L_M(w)+\log(\ell+1))\le C_d(\ell+1)L_M(w).
		\]
		Since $|T_J(z)|=|z|^{2M}/(2^d4^J J!M!)$, summing the tail yields
		\[
		\begin{aligned}
			|F_M(z)|\le C_d|T_J(z)|L_M(w)\sum_{\ell=0}^{\infty}(\ell+1)\rho^\ell\le\frac{C_d^{M+1}s^{2M}}{J!M!}L_M(w),
		\end{aligned}
		\]
		where the series equals $(1-\rho)^{-2}$, and the factor $(5/4)^{2M}$ from $|z|^{2M}$ is included in $C_d^{M+1}$. Moreover, the Cauchy integral formula on this circle gives
		\[
		|F_M^{(k)}(w)|
		\le k!(4/s)^k\sup_{|z-w|=s/4}|F_M(z)|
		\le\frac{C_d^{M+k+1}k!s^{2M-k}}{J!M!}L_M(w),
		\]
		which proves \eqref{eq:Bessel-remainder-derivatives} in this region.
		
		Next, suppose that $s\ge J+1$, and put $Q_M(s)=s^{2M}/(J!M!)$. We use the decomposition
		\begin{equation}\label{eq:F_M-decompo}
			F_M(z)=z^dK_d(z)-z^dS_d(z)-\sum_{r=0}^{J-1}T_r(z)(\alpha_d\log(z/2)+\beta_{d,r}).
		\end{equation}
		For $0\le r<J$, estimate~\eqref{eq:Bessel-power-comparison}, with $m=r$ and $M=J+d$, gives $\frac{J!M!}{4^r r!(r+d)!s^{2(J-r)}}\le C_d4^{-r}$. On $|z-w|=s/4$, the definition of $T_r$ gives
		\[
		\frac{|T_r(z)|}{Q_M(s)}=\frac1{2^d}\Big(\frac{|z|}{s}\Big)^{2(r+d)}\frac{J!M!}{4^r r!(r+d)!s^{2(J-r)}}\le C_d(5/4)^{2(r+d)}4^{-r}\le C_d\rho^r.
		\]
		For $r<J$, \eqref{eq:bound_beta_m} gives $|\alpha_d\log(z/2)+\beta_{d,r}|\le C_dL_M(w)$. Using $\sum_{r=0}^{J-1}\rho^r\le(1-\rho)^{-1}$, we obtain
		\[
		\sup_{|z-w|=s/4}\Big|\sum_{r=0}^{J-1}T_r(z)(\alpha_d\log(z/2)+\beta_{d,r})\Big|\le C_dQ_M(s)L_M(w).
		\]
		
		For the second term on the right-hand side of \eqref{eq:F_M-decompo}, definition~\eqref{eq:Bessel-finite-part} shows that $z^dS_d(z)$ is a fixed polynomial of degree at most $2d-2$ when $d\ge1$. Taking $m=0$ in \eqref{eq:Bessel-power-comparison} gives $J!M!/(d!s^{2J})\le C_d$, or $s^{2d}\le C_dQ_M(s)$. Since $s\ge J+1\ge3$, we have $\sup_{|z-w|=s/4}|z^dS_d(z)|\le C_ds^{2d-2}\le C_dQ_M(s)$. For $d=0$, this term vanishes because $S_0=0$. Cauchy's integral formula with radius $s/4$ gives the required derivative bounds for this polynomial and the finite sum.
		
		For the remaining term $w^dK_d(w)$ in the right-hand side of \eqref{eq:F_M-decompo}, use the circle $|z-w|=J/4$. Since $\Rea w\ge0$ and $s=|w|\ge J+1\ge3$, every point on this circle satisfies $3s/4\le|z|\le5s/4$, $\Rea z\ge-J/4$, and $|\arg z|\le\pi/2+\arcsin(1/4)<3\pi/4$.
		
		We first derive a bound for $K_d$ on the fixed sector $|\arg z|\le3\pi/4$. Taking the Bessel order to be $d$ and using $a_0(d)=1$ in \cite[Eq.~(10.40.2)]{NIST10}, we have
		\[
		K_d(z)=\Big(\frac{\pi}{2z}\Big)^{1/2}e^{-z}(1+\eta_d(z)),\quad |\eta_d(z)|\le\frac{A_d}{|z|}
		\]
		for all sufficiently large $|z|$, uniformly in this sector. The uniform remainder estimate can also be read from \cite[Eqs.~(10.40.10)--(10.40.12)]{NIST10} with $\ell=1$. Choose $R_0\ge\max\{1,A_d\}$ sufficiently large that this estimate holds for $|z|\ge R_0$. Then $|\eta_d(z)|\le1$, and taking absolute values gives
		\[
		|z^dK_d(z)|\le2(\pi/2)^{1/2}|z|^{d-1/2}e^{-\Rea z},\quad |z|\ge R_0,\quad |\arg z|\le3\pi/4.
		\]
		To cover $1\le|z|\le R_0$, consider the compact set $E_d=\{z:1\le|z|\le R_0,\ |\arg z|\le3\pi/4\}$. Since $|z|\ge1$ and $|\arg z|<\pi$ on $E_d$, it is disjoint from the principal branch cut $(-\infty,0]$. The principal branches of $z^{1/2}$ and $K_d(z)$ are holomorphic on $\C\setminus(-\infty,0]$, so the continuous real-valued function $|z^{1/2}e^zK_d(z)|$ attains a finite maximum $B_d$ on $E_d$. The identity $|z^dK_d(z)|=|z|^{d-1/2}e^{-\Rea z}|z^{1/2}e^zK_d(z)|$ shows that choosing $C_d\ge\max\{2\sqrt{\pi/2},B_d\}$ gives $|z^dK_d(z)|\le C_d|z|^{d-1/2}e^{-\Rea z}$ throughout the sector with $|z|\ge1$.
		
		On the circle $|z-w|=J/4$, we have $|z|\ge3s/4\ge9/4>1$ and $e^{-\Rea z}\le e^{J/4}$. For $d\ge1$, the upper bound $|z|\le5s/4$ gives $|z|^{d-1/2}\le(5/4)^{d-1/2}s^{d-1/2}$; for $d=0$, the lower bound gives $|z|^{-1/2}\le(4/3)^{1/2}s^{-1/2}$. Enlarging $C_d$ to include these factors, we obtain
		\[
		|z^dK_d(z)|\le C_d|z|^{d-1/2}e^{-\Rea z}\le C_de^{J/4}s^{d-1/2},\quad |z-w|=J/4.
		\]
		Cauchy's integral formula gives $|(w^dK_d(w))^{(k)}|\le C_dk!(4/J)^ke^{J/4}s^{d-1/2}$. Since $M=J+d$ and $0\le k\le2J-2$, we have $k+d-1/2-2M\le-d-5/2<0$. Estimate~\eqref{eq:Bessel-factorial-bound} and $s\ge J+1$ give
		\[
		J!M!J^{-k}s^{k+d-1/2-2M}\le d!J^{-k}(J+1)^{k-1/2}=d!(1+1/J)^k(J+1)^{-1/2}\le d!e^2,
		\]
		where $(1+1/J)^k\le e^{k/J}\le e^2$. Substitution into the Cauchy estimate gives
		\[
		|(w^dK_d(w))^{(k)}|\le\frac{C_d4^ke^{J/4}k!s^{2M-k}}{J!M!}\le\frac{C_d^{M+k+1}k!s^{2M-k}}{J!M!}.
		\]
		Combining the bounds for the three terms in \eqref{eq:F_M-decompo} and using $L_M(w)\ge1$ proves \eqref{eq:Bessel-remainder-derivatives}.
	\end{proof}
	
	\begin{lemma}\label{lem:even-logarithmic-remainder}
		Let $n=2\nu\ge2$, $d=\nu-1$, $p\ge0$, and $a>0$, and assume \eqref{eq:ellipticity} and \eqref{eq:compact-perturbation}. Suppose that $e^{aw_p(|x|)}u\in L^2(\R^n)$, and fix $\phi\in C_c^\infty(\R^n)$. Put $a_u(\kappa)=\ip{\mathcal R_0(\kappa)u}{\phi}$ and $g_u(\kappa)=\nabla\mathcal R_0(\kappa)u|_D$, and let $h_u$ be defined by \eqref{eq:even-scalar-reconstruction}. The coefficient kernels in Lemma~\ref{lem:even-source-expansion} define scalar coefficients $a_r,b_r$ and vector coefficients $g_r,h_r\in\mathcal V$, independent of the truncation order, such that
		\begin{equation}\label{eq:even-logarithmic-coefficients}
			|a_r|+|b_r|+\norm{g_r}_{\mathcal V}+\norm{h_r}_{\mathcal V}\le C_0C_1^r(r!)^2\log^{4pr}(e+r),\quad r\in\Nzero.
		\end{equation}
		There is $0<\eps<1$ such that the finite expansions \eqref{eq:even-source-expansions} hold for every integer $M\ge d+2$ and every $\kappa$ with $0<|\kappa|<\eps$ and $\Rea\kappa\ge0$. Their remainders satisfy
		\begin{equation}\label{eq:even-logarithmic-remainder}
			|\partial_\kappa^kE_{a,M}(\kappa)|+\norm{\partial_\kappa^kE_{g,M}(\kappa)}_{\mathcal V}\le C_0C_1^{M+k}k!(M!)^2\log^{4pM}(e+M)t^{2M-k}\ell(t),
		\end{equation}
		where $t=|\kappa|$, $\ell(t)=\log(e/t)$, and $0\le k\le2(M-d)-3$. At nonzero points of the imaginary axis, the functions and derivatives in these statements are understood as continuous limits from the right half-plane.
		
		Let $P_r$ be defined by \eqref{eq:power-log-polynomial-definition}. The finite expansion \eqref{eq:h_u} holds for every $M\ge d+2$. If every $P_r$ vanishes identically, every derivative of $h_u$ extends continuously to $\overline{U_\eps}$, satisfies $h_u^{(m)}(0)=0$, and obeys
		\begin{equation}\label{eq:even-logarithmic-derivatives}
			\sup_{\overline{U_\eps}}|h_u^{(m)}|\le C_0C_1^m m!M_m,\quad m\in\Nzero,
		\end{equation}
		after decreasing $\eps$ and increasing the constants if necessary. Here $U_\eps=\{\kappa\in\C:\Rea\kappa>0,\ |\kappa|<\eps\}$ and $(M_m)$ is defined by \eqref{eq:def_M_m}. The constants may depend on $n,p,a,A,D,\phi$ and $\norm{e^{aw_p(|\cdot|)}u}_2$, but are independent of $r,M,k,m,\kappa$.
	\end{lemma}
	
	\begin{proof}
		Write $L_r=\log(e+r)$ and $W_r=(r!)^2L_r^{4pr}$. Applying \eqref{eq:logarithmic-moments} with the fixed polynomial weight $q+1$, and using $M_m\le m!L_m^{2pm}$ and $1+\log(2+|y|)\le C(1+|y|)$, gives
		\begin{equation}\label{eq:even-logarithmic-moment-bound}
			\int_{\R^n}(1+|y|)^{m+q}(1+\log(2+|y|))|u(y)|\,dy\le C_0C_1^m(m!)^2L_m^{2pm}
		\end{equation}
		for every fixed $q\ge0$, with constants independent of $m$.
		
		Use the kernels \eqref{eq:even-coefficient-kernels-low}--\eqref{eq:even-coefficient-kernels-high} and the coefficient definitions \eqref{eq:even-source-coefficient-definitions}. The local convolution bounds follow from \eqref{eq:even-coefficient-local-bound} and Schur's estimate \eqref{eq:schur-estimate}, using only $u\in L^2$. For $r\ge d$, integrate the distant kernel bound \eqref{eq:even-coefficient-far-bound} using \eqref{eq:even-logarithmic-moment-bound} with $m=2(r-d)$. The inequality $L_{2(r-d)}^{4p(r-d)}\le2^{4pr}L_r^{4pr}$ gives
		\[
		|a_r|+|b_r|+\norm{g_r}_{\mathcal V}+\norm{h_r}_{\mathcal V}
		\le C_0C_1^r\frac{((2r)!)^2}{(r-d)!r!}L_r^{4pr}.
		\]
		Estimate~\eqref{eq:even-factorial-comparison} bounds the factorial quotient by $16^r(r+1)^d(r!)^2$. Including $(r+1)^d$ in the exponential constant proves \eqref{eq:even-logarithmic-coefficients}; increasing $C_0$ includes $r<d$. The bounds \eqref{eq:even-coefficient-local-bound}, \eqref{eq:even-coefficient-far-bound}, and \eqref{eq:even-logarithmic-moment-bound} also justify the coefficient convolutions in local $L^2$.
		
		For the remainders, put $J=M-d\ge2$. By \eqref{eq:even-remainder-kernel-definition} and the definition of $F_M$ in Lemma~\ref{lem:Bessel-remainder-derivatives}, $\mathcal E_M(\kappa,x)=c_n\rho^{-2d}F_M(\kappa\rho)$, where $\rho=|x|>0$. Its derivatives satisfy
		\[
		\begin{aligned}
			\partial_\kappa^k\mathcal E_M(\kappa,x)&=c_n\rho^{k-2d}F_M^{(k)}(\kappa\rho),\\
			\partial_\rho\partial_\kappa^k\mathcal E_M(\kappa,x)&=c_n\rho^{k-2d-1}((k-2d)F_M^{(k)}(\kappa\rho)+\kappa\rho F_M^{(k+1)}(\kappa\rho)).
		\end{aligned}
		\]
		For $0\le k\le2J-3$, estimate~\eqref{eq:Bessel-remainder-derivatives} applies to both $F_M^{(k)}$ and $F_M^{(k+1)}$. Using $(k+1)!=(k+1)k!$ and including the factors linear in $k$ in the exponential constant gives
		\begin{equation}\label{eq:even-logarithmic-kernel-derivatives}
			\begin{aligned}
				|\partial_\kappa^k\mathcal E_M(\kappa,x)|+\rho|\nabla_x\partial_\kappa^k\mathcal E_M(\kappa,x)|
				\le\frac{C^{M+k+1}k!t^{2M-k}\rho^{2J}}{J!M!}(\ell(t)+|\log\rho|+\log(M+2)).
			\end{aligned}
		\end{equation}
		The spatial power is independent of $k$, since $\rho^{k-2d}(t\rho)^{2M-k}=t^{2M-k}\rho^{2J}$.
		
		Integrate \eqref{eq:even-logarithmic-kernel-derivatives}, splitting at $\rho=1$. On $\rho\le1$, combine \eqref{eq:even-local-remainder-integrals} with Schur's estimate \eqref{eq:schur-estimate} for the scalar and gradient kernels. On $\rho\ge1$, use \eqref{eq:even-logarithmic-moment-bound} with $m=2J$. Including $1+\log(M+2)$ in the exponential constant and using $\ell(t)\ge1$ gives
		\[
		|\partial_\kappa^kE_{a,M}(\kappa)|+\norm{\partial_\kappa^kE_{g,M}(\kappa)}_{\mathcal V}\le C_0C_1^{M+k}k!t^{2M-k}\ell(t)
		\frac{((2J)!)^2L_{2J}^{4pJ}}{J!M!}.
		\]
		Here $((2J)!)^2\le((2M)!)^2$ and $L_{2J}^{4pJ}\le2^{4pM}L_M^{4pM}$. Estimate~\eqref{eq:even-factorial-comparison} with $r=M$ proves \eqref{eq:even-logarithmic-remainder}. The majorants in \eqref{eq:even-logarithmic-kernel-derivatives}, together with \eqref{eq:even-local-remainder-integrals} and \eqref{eq:even-logarithmic-moment-bound}, justify differentiation in local $L^2$ and continuous boundary limits at nonzero points of the imaginary axis by the dominated convergence theorem. Integrating \eqref{eq:even-finite-kernel-identity} proves \eqref{eq:even-source-expansions}, with the coefficients \eqref{eq:even-source-coefficient-definitions} independent of $M$.
		
		The bounds \eqref{eq:internal-power-log-bounds} for $S_\phi$ depend only on $A,D,\phi$. The identity \eqref{eq:finite-scalar-reconstruction} applies to the finite sums \eqref{eq:power-log-partial-sums}. In the derivation of \eqref{eq:h_u}, replace \eqref{eq:even-source-remainders} by \eqref{eq:even-logarithmic-remainder} with $k=0$; the polynomial identity \eqref{eq:power-log-polynomial-definition} and the tail bound \eqref{eq:internal-power-log-bounds} apply unchanged. This proves \eqref{eq:h_u} for each fixed $M$.
		
		Suppose now that every $P_r$ vanishes. The representation \eqref{eq:internal-S-holomorphic-representation} extends $S_\phi(\kappa)=\mathscr S_\phi(\kappa^2,\kappa^{2\nu}\log\kappa)$ as a bounded holomorphic function to
		\[
		\Sigma_{2\eps}=\{\kappa\in\C:0<|\kappa|<2\eps,\ |\arg\kappa|<3\pi/4\},
		\]
		after decreasing $\eps$. Both arguments of $\mathscr S_\phi$ tend uniformly to zero in this sector. The coefficient estimate \eqref{eq:internal-Taylor-coefficients} and the definition of $Q_q$ in \eqref{eq:internal-power-log-series} give $\norm{Q_q(\log\xi)}_{\mathcal V}\le C_0C_1^q(1+|\log\xi|)^q$ there. The series converges absolutely when $C_1|\xi|^2(1+|\log\xi|)<1$, which holds throughout the sector after decreasing $\eps$. For $0<t=|\kappa|<\eps$ and $\Rea\kappa\ge0$, the disk $|\xi-\kappa|\le t/4$ lies in $\Sigma_{2\eps}$. Cauchy's integral formula gives
		\begin{equation}\label{eq:even-logarithmic-S-derivatives}
			\norm{S_\phi^{(k)}(\kappa)}_{\mathcal V}\le k!(4/t)^k\sup_{|\xi-\kappa|=t/4}\norm{S_\phi(\xi)}_{\mathcal V}\le C_0C_1^k k!t^{-k}.
		\end{equation}
		
		Let $A_M,G_M$ be defined by \eqref{eq:power-log-partial-sums}, and put $V_M=A_M-[G_M,S_\phi]_D$. Equations~\eqref{eq:even-scalar-reconstruction} and \eqref{eq:even-source-expansions} give $h_u=V_M+E_{a,M}-[E_{g,M},S_\phi]_D$. By \eqref{eq:power-log-polynomial-definition} and $P_r\equiv0$, all terms of total index less than $M$ vanish. Multiplication by the convergent series \eqref{eq:internal-power-log-series} gives
		\begin{equation}\label{eq:even-logarithmic-V-series}
			V_M(\kappa)=-\sum_{r=0}^{M-1}\sum_{q=M-r}^{\infty}\kappa^{2(r+q)}[g_r+h_r\log\kappa,Q_q(\log\kappa)]_D.
		\end{equation}
		This identity holds on $\Sigma_{2\eps}$, where $V_M$ is holomorphic.
		
		On the circle $|\xi-\kappa|=t/4$, we have $|\xi|\le5t/4$ and $1+|\log\xi|\le C\ell(t)$. Applying \eqref{eq:even-logarithmic-coefficients} and the polynomial bound obtained from \eqref{eq:internal-Taylor-coefficients} to \eqref{eq:even-logarithmic-V-series} gives
		\[
		\sup_{|\xi-\kappa|=t/4}|V_M(\xi)|
		\le C_0\ell(t)\sum_{r=0}^{M-1}(Ct^2)^rW_r\sum_{q=M-r}^{\infty}(Ct^2\ell(t))^q.
		\]
		Choose $\eps$ with $\theta(t)=Ct^2\ell(t)\le1/2$ for $0<t<\eps$. The inner sum equals $\theta(t)^{M-r}/(1-\theta(t))\le2\theta(t)^{M-r}$, which gives
		\[
		\sup_{|\xi-\kappa|=t/4}|V_M(\xi)|
		\le C_0C_1^Mt^{2M}\sum_{r=0}^{M-1}W_r\ell(t)^{M-r+1}.
		\]
		Another application of Cauchy's integral formula yields
		\begin{equation}\label{eq:even-logarithmic-V-derivatives}
			|V_M^{(m)}(\kappa)|
			\le C_0C_1^{M+m}m!t^{2M-m}\sum_{r=0}^{M-1}W_r\ell(t)^{M-r+1}.
		\end{equation}
		
		For $m\le2(M-d)-3$, apply Leibniz' rule to $[E_{g,M},S_\phi]_D$ using \eqref{eq:even-logarithmic-remainder} and \eqref{eq:even-logarithmic-S-derivatives}. The term containing $\partial_\kappa^kE_{g,M}$ and $S_\phi^{(m-k)}$ has factorial factor $\binom mkk!(m-k)!=m!$ and power $t^{2M-k}t^{-(m-k)}=t^{2M-m}$. Summing the $m+1\le2^m$ terms and combining the scalar remainder with \eqref{eq:even-logarithmic-V-derivatives} gives
		\begin{equation}\label{eq:even-logarithmic-flat-derivative-bound}
			|h_u^{(m)}(\kappa)|
			\le C_0C_1^{M+m}m!t^{2M-m}
			\Bigl(W_M\ell(t)+\sum_{r=0}^{M-1}W_r\ell(t)^{M-r+1}\Bigr).
		\end{equation}
		For each fixed $m$, choose an admissible $M$. The right-hand side is $O_{M,m}(t^{2M-m}\ell(t)^{M+1})$ and tends to zero uniformly in the closed right half-plane. Every derivative extends continuously to zero with value zero. Continuity at the other boundary points follows from the continuous boundary limits of the source remainders established using \eqref{eq:even-logarithmic-kernel-derivatives}, together with the holomorphy of $V_M$ and $S_\phi$ on $\Sigma_{2\eps}$.
		
		For a bound uniform in $m$, choose $M=\lceil m/2\rceil+d+2$. Then $m\le2(M-d)-3$ and $2M-m\ge1$. For $k\ge0$, writing $v=\ell(t)$ and using $e^v\ge v^{k+1}/(k+1)!$ gives
		\[
		\sup_{0<t<1}t^{2M-m}\ell(t)^{k+1}
		\le\sup_{0<t<1}t\ell(t)^{k+1}\le e(k+1)!.
		\]
		For $r<M$, set $k=M-r\ge1$. Since $M!/r!=\prod_{\ell=1}^k(r+\ell)\ge k!$ and $L_r^{4pr}\le L_M^{4pM}$,
		\[
		W_r(k+1)!\le W_M\frac{(k+1)!}{(k!)^2}
		=W_M\frac{k+1}{k!}\le2W_M.
		\]
		Using these bounds in \eqref{eq:even-logarithmic-flat-derivative-bound}, with $k=0$ for the term $W_M\ell(t)$, bounds the sum by $e(2M+1)W_M$. Including $2M+1\le3^M$ in the exponential constant gives
		\[
		\sup_{\overline{U_\eps}}|h_u^{(m)}|\le C_0C_1^{M+m}m!W_M.
		\]
		
		To compare $W_M$ with the required weight, put $b=\lfloor m/2\rfloor$ and $c=d+3$. Then $M-b\le c$, $M\le m+c$, and $4pM\le2pm+4pc$. The fixed factorial shift and $(b!)^2\le m!$ give
		\[
		(M!)^2\le C_d(m+1)^{2c}(b!)^2\le C_d(m+1)^{2c}m!.
		\]
		Also, $L_M\le C_dL_m$, so $L_M^{4pM}\le C_0C_1^mL_m^{2pm+4pc}$. For $m\ge1$, the remaining factors satisfy $(m+1)^{2c}L_m^{4pc}\le2^{(2c+4pc)m}$, because $L_m\le m+1\le2^m$. Increasing $C_0$ includes $m=0$. We obtain $W_M\le C_0C_1^m m!L_m^{2pm}$, and the preceding derivative bound becomes
		\[
		\sup_{\overline{U_\eps}}|h_u^{(m)}|\le C_0C_1^m(m!)^2L_m^{2pm}.
		\]
		Raising \eqref{eq:logarithmic-product} to the power $2p$ gives $L_m^{2pm}\le e^{2pm}\prod_{\ell=1}^mL_\ell^{2p}$. Using \eqref{eq:def_M_m} in the preceding derivative bound and replacing $C_1$ by $e^{2p}C_1$ proves \eqref{eq:even-logarithmic-derivatives}.
	\end{proof}
	
	\begin{proof}[Proof of Proposition~\ref{prop:even-logarithmic}]
		First, assume that there exists $a>0$ such that $e^{aw_p(|x|)}u_k\in L^2(\R^n)$ for every $k=1,\ldots,N$. Fix $\phi\in C_c^\infty(\cO)$ and put $h_k(\kappa)=\ip{\mathcal R_A(\kappa)u_k}{\phi}$. Proposition~\ref{prop:resolvent} gives \eqref{eq:odd-separation}. By Lemma~\ref{lem:even-logarithmic-remainder}, each $h_k$ has the finite expansion \eqref{eq:h_u}, with polynomials independent of the truncation order.
		
		The minimum-exponent argument based on \eqref{eq:h_k_leading} and leading to \eqref{eq:even-separated-coefficients} applies to these expansions because $s_k-s_\ell\notin\Z$. It gives $P_{k,r}\equiv0$ for every $k,r$. Lemma~\ref{lem:even-logarithmic-remainder} gives continuous extensions of all derivatives of $h_k$ to $\overline{U_\eps}$, with $h_k^{(m)}(0)=0$ and the bound \eqref{eq:even-logarithmic-derivatives}. Since the number of inputs is finite, $\eps,C_0,C_1$ can be chosen independently of $k$.
		
		Equation~\eqref{eq:logarithmic-quasianalytic-series} verifies \eqref{eq:sector-quasianalyticity} for $0\le p\le1$. Applied with \eqref{eq:even-logarithmic-derivatives} and the vanishing derivatives at zero, the quasianalyticity criterion gives $h_k=0$ in $U_\eps$. The analytic continuation argument leading to \eqref{eq:continued-resolvent-pairings} gives $\ip{R_\lambda u_k}{\phi}=0$ for every $\lambda>0$. Since $\phi\in C_c^\infty(\cO)$ is arbitrary, Lemma~\ref{lem:individual} gives $u_k=0$ on $\R^n$.
		
		For $u_k\in\mathscr D_{\log,p}$, estimate~\eqref{eq:logarithmic-smoothed-weight}, whose proof is independent of the parity of $n$, gives $e^{a'w_p(|x|)}S_\delta u_k\in L^2(\R^n)$ for some $a'>0$ and every $k$. Equation~\eqref{eq:smoothed-entanglement} gives the required relation on any fixed nonempty ball $\cO_0\Subset\cO$ when $\delta$ is sufficiently small. The weighted $L^2$ case just proved gives $S_\delta u_k=0$. Lemma~\ref{lem:smoothing} and \eqref{eq:spectral-sobolev-equivalence} give $S_\delta u_k\to u_k$ in $H^{-r}(\R^n)$ for a common $r\ge0$, proving $u_k=0$ for every $k$.
	\end{proof}

	\section{Counterexamples and sharpness}\label{sec:sharpness}
	
	We prove Theorem~\ref{thm:sharpness} and Proposition~\ref{prop:logarithmic-sharpness} by constructing pairs whose weighted resolvents cancel on $B_R$. We prove Theorem~\ref{thm:sharpness}\ref{item_2_thm_sharp} by splitting the proof into odd- and even-dimensional cases. The converse part of Proposition~\ref{prop:resolvent} then gives the required fractional relation. In dimension three, a radial resolvent formula reduces the construction to a relation between Laplace transforms. Integer derivatives of compactly supported profiles give examples for half-integer differences of the orders. For arbitrary differences, we use profiles with all integer moments zero and prove that their causal fractional derivatives preserve the required decay. An Abel integral gives the corresponding construction in dimension two. We then pass to higher dimensions by a radial differential identity and transfer the examples to every compact coefficient perturbation. We prove the logarithmic counterexamples in Subsection~\ref{subsec:logarithmic-sharpness}.
	
	\subsection{Radial resolvents and outgoing waves}
	
	Let $F\in\cS(\R)$ vanish for $r\le R$, where $R>0$. Define $v_F(x)=F(|x|)/|x|$ on $\R^3$, with value zero near the origin, and denote the Laplace transform of $F$ by $\mathcal LF(\kappa)=\int_0^\infty e^{-\kappa t}F(t)\,dt$ for $\kappa>0$.
	
	\begin{lemma}\label{lem:radial-resolvent}
		For $\kappa>0$ and $0\le r<R$,
		\begin{equation}\label{eq:radial-resolvent}
			(\mathcal R_0(\kappa)v_F)(r)=\frac{\sinh(\kappa r)}{\kappa r}\mathcal LF(\kappa),
		\end{equation}
		where the quotient is defined by continuity at $r=0$.
	\end{lemma}
	
	\begin{proof}
		The resolvent solution $w(r)=(\mathcal R_0(\kappa)v_F)(r)$ is radial. Multiplication of its radial equation by $r$ shows that $y(r)=rw(r)$ solves $-y''+\kappa^2y=F$ on $(0,\infty)$, with $y(0)=0$ and $y\in L^2(0,\infty)$. The Green kernel for this problem is $(2\kappa)^{-1}(e^{-\kappa|r-t|}-e^{-\kappa(r+t)})$. Since $r<R$ and $F(t)=0$ for $t\le R$, the part of the kernel contributing to the integral is $\kappa^{-1}\sinh(\kappa r)e^{-\kappa t}$. Integrating against $F$ gives $y(r)=\kappa^{-1}\sinh(\kappa r)\mathcal LF(\kappa)$. Division by $r$ gives \eqref{eq:radial-resolvent}, including its value at zero by continuity.
	\end{proof}
	
	For nonintegral $\beta>0$, put $m=\lceil\beta\rceil$ and $\delta=m-\beta\in(0,1)$. For a smooth profile $F$ vanishing near $(-\infty,0]$, define its causal fractional derivative by
	\begin{equation}\label{eq:causal-derivative}
		D_+^\beta F(r)=\frac1{\Gamma(\delta)}\int_0^r(r-t)^{\delta-1}F^{(m)}(t)\,dt,\quad r>0,
	\end{equation}
	and set it equal to zero for $r\le0$. For integral $\beta=m$, set $D_+^\beta F=F^{(m)}$. If $F\in\cS(\R)$, integration by parts gives $\mathcal L(F^{(m)})(\kappa)=\kappa^m\mathcal LF(\kappa)$, since all derivatives of $F$ vanish at zero. For nonintegral $\beta$, Fubini's theorem and the Gamma integral $\Gamma(\delta)^{-1}\int_0^\infty e^{-\kappa r}r^{\delta-1}\,dr=\kappa^{-\delta}$ supply the additional factor $\kappa^{-\delta}$ when taking the Laplace transform in \eqref{eq:causal-derivative}. In both cases, we obtain
	\begin{equation}\label{eq:Laplace-derivative}
		\mathcal L(D_+^\beta F)(\kappa)=\kappa^\beta\mathcal LF(\kappa),\quad\kappa>0.
	\end{equation}
	The exponential factor ensures absolute convergence in these calculations.
	
	\begin{proposition}\label{prop:radial-cancellation}
		Let $0<s<t$ be nonintegers and put $\beta=2(t-s)$. Suppose that $F,D_+^\beta F\in\cS(\R)$ and both vanish for $r\le R$. Then the functions
		\begin{equation}\label{eq:radial-inputs}
			u=-\sin(\pi t)\,v_{D_+^\beta F},\quad v=\sin(\pi s)\,v_F
		\end{equation}
		are zero on $B_R$ and satisfy $((-\Delta)^su+(-\Delta)^tv)|_{B_R}=0$ in dimension three.
	\end{proposition}
	
	\begin{proof}
		Equations~\eqref{eq:radial-resolvent} and \eqref{eq:Laplace-derivative} give $\mathcal R_0(\kappa)v_{D_+^\beta F}=\kappa^\beta\mathcal R_0(\kappa)v_F$ on $B_R$ for every $\kappa>0$. Substituting \eqref{eq:radial-inputs} and using $2s+\beta=2t$, we obtain
		\begin{equation}\label{eq:radial-cancellation}
			(\sin(\pi s)\kappa^{2s}\mathcal R_0(\kappa)u+\sin(\pi t)\kappa^{2t}\mathcal R_0(\kappa)v)|_{B_R}=0.
		\end{equation}
		Both inputs vanish on $B_R$, so the converse part of Proposition~\ref{prop:resolvent}, with $\lambda=\kappa^2$, gives the fractional relation.
	\end{proof}
	
	The case $(s,t)=(1/4,3/4)$ has a direct wave interpretation. Consider
	\begin{equation}\label{eq:outgoing-wave}
		w(\tau,x)=\frac{F(|x|-\tau)}{|x|},\quad\tau\ge0,\quad x\in\R^3\setminus\{0\},
	\end{equation}
	and set $w(\tau,0)=0$. Since $F$ vanishes for $r\le R$, this defines a smooth function on $[0,\infty)\times\R^3$. It solves $w_{\tau\tau}-\Delta w=0$ in $(0,\infty)\times\R^3$ and vanishes on $B_R$ for every $\tau\ge0$. The functions in \eqref{eq:radial-inputs} are $u=\sin(\pi/4)w_\tau(0,\cdot)$ and $v=\sin(\pi/4)w(0,\cdot)$. The Laplace transform in $\tau$ satisfies $\mathcal L_\tau w=\mathcal R_0(\kappa)(w_\tau(0,\cdot)+\kappa w(0,\cdot))$. Since $\mathcal L_\tau w$ is zero on $B_R$, multiplication by $\sin^2(\pi/4)\kappa^{1/2}$ gives \eqref{eq:radial-cancellation}.
	
	\begin{lemma}\label{lem:dimension}
		Let $n\ge2$, and write $\Delta_m$ for the Laplacian on $\R^m$. For a radial Schwartz function $f$ on $\R^n$, write $f(x)=F(|x|)$, where $F$ is a smooth even function on $\R$. Define $\mathcal Tf$ on $\R^{n+2}$ by $(\mathcal Tf)(y)=-F'(|y|)/|y|$ for $y\ne0$ and $(\mathcal Tf)(0)=-F''(0)$. Then $\mathcal Tf$ is a radial Schwartz function, and for every $\kappa>0$,
		\begin{equation}\label{eq:dimension-identity}
			\Delta_{n+2}\mathcal Tf=\mathcal T\Delta_nf,\quad (-\Delta_{n+2}+\kappa^2)^{-1}\mathcal Tf=\mathcal T(-\Delta_n+\kappa^2)^{-1}f.
		\end{equation}
		The operator $\mathcal T$ is injective and preserves compact support, vanishing on a ball, and stretched-exponential decay of all derivatives.
	\end{lemma}
	
	\begin{proof}
		For $r=|x|>0$, the radial formula $(\Delta_nf)(x)=F''(r)+(n-1)r^{-1}F'(r)$ gives the first identity by direct differentiation. Since $F$ is even, $F'(0)=0$, and $\frac{F'(r)}r=\int_0^1 F''(\theta r)\,d\theta$. The right-hand side is smooth and even in $r$, with value $F''(0)$ at zero. This gives the smooth radial extension of $\mathcal Tf$ across the origin, where the first identity also holds by continuity. At infinity, differentiation produces derivatives of $F$ multiplied by nonpositive powers of $r$ and bounded angular factors. These formulas show that $\mathcal Tf$ is Schwartz and has the stated decay properties. The support properties follow directly from its definition.
		
		For fixed $\kappa>0$, the free resolvent preserves radial Schwartz functions, as follows from its Fourier multiplier $(|\xi|^2+\kappa^2)^{-1}$. If $w=(-\Delta_n+\kappa^2)^{-1}f$, the first identity gives $(-\Delta_{n+2}+\kappa^2)\mathcal Tw=\mathcal Tf$. Since $\mathcal Tw$ is Schwartz, uniqueness of the $L^2$ resolvent solution proves the second identity. Finally, $\mathcal Tf=0$ implies $F'=0$ on $(0,\infty)$. The decay of $F$ forces this constant profile to be zero, so $f=0$.
	\end{proof}
	
	\begin{lemma}\label{lem:transfer}
		Suppose that $\supp(A-I)\Subset B_R$. Let $\kappa>0$ and $f_\kappa\in\cS(\R^n)$. If $W_\kappa=(-\Delta+\kappa^2)^{-1}f_\kappa$ vanishes on $B_R$, then $W_\kappa=(L+\kappa^2)^{-1}f_\kappa$.
	\end{lemma}
	
	\begin{proof}
		By assumption, $W_\kappa=0$ on $B_R$, so $\nabla W_\kappa=0$ there. The support condition gives $A-I=0$ on $\R^n\setminus B_R$. The product therefore vanishes in both regions: $(A-I)\nabla W_\kappa=0$ in $\R^n$. Using this identity and the defining equation for $W_\kappa$, we obtain 
		\[
		\begin{aligned}
			(L+\kappa^2)W_\kappa=(-\Delta+\kappa^2)W_\kappa-\operatorname{div}((A-I)\nabla W_\kappa)=(-\Delta+\kappa^2)W_\kappa=f_\kappa	\quad\text{in }\R^n.
		\end{aligned}
		\]
		The function $W_\kappa$ is Schwartz and belongs to $\Dom(L)$ under \eqref{eq:compact-perturbation}. Uniqueness of the $L^2(\R^n)$ resolvent solution gives $W_\kappa=(L+\kappa^2)^{-1}f_\kappa$.
	\end{proof}
	
	\begin{proof}[Proof of Theorem~\ref{thm:sharpness}\ref{item_1_thm_sharp}]
		Choose a nonzero real $F\in C_c^\infty((R+1,R+2))$. If $t-s=k/2$ with $k\in\N$, then $D_+^{2(t-s)}F=F^{(k)}$ is smooth and compactly supported. It is nonzero, since a compactly supported function with identically vanishing $ k$-th derivative must be zero. Since $\sin(\pi s)$ and $\sin(\pi t)$ are nonzero, Proposition~\ref{prop:radial-cancellation} gives nonzero radial functions $u_3,v_3\in C_c^\infty(\R^3\setminus\overline{B_R})$ such that, for every $\kappa>0$,
		\[
		W_{3,\kappa}:=(-\Delta_3+\kappa^2)^{-1}(\sin(\pi s)\kappa^{2s}u_3+\sin(\pi t)\kappa^{2t}v_3)=0
		\quad\text{on }B_R\subset\R^3.
		\]
		
		Let $n=3+2q$ with $q\in \Nzero$, and set $u=\mathcal T^q u_3$ and $v=\mathcal T^q v_3$. Each application of $\mathcal T$ increases the dimension by two, and $\mathcal T^0$ denotes the identity. Lemma~\ref{lem:dimension} gives $u,v\in C_c^\infty(\R^n\setminus\overline{B_R})$, and its injectivity ensures that both functions are nonzero. For $\kappa>0$, define $f_\kappa=\sin(\pi s)\kappa^{2s}u+\sin(\pi t)\kappa^{2t}v$. Linearity of $\mathcal T$ and repeated use of \eqref{eq:dimension-identity} give
		\[
		(-\Delta_n+\kappa^2)^{-1}f_\kappa
		=\mathcal T^q W_{3,\kappa}=0
		\quad\text{on }B_R\subset\R^n.
		\]
		The last equality uses the fact that $\mathcal T$ preserves vanishing on a ball.
		
		Since $\supp(A-I)\Subset B_R$, Lemma~\ref{lem:transfer}, applied in $\R^n$, gives $(L+\kappa^2)^{-1}f_\kappa=0$ on $B_R$ for every $\kappa>0$. Using $u=v=0$ on $B_R$, the converse part of Proposition~\ref{prop:resolvent} proves \eqref{eq:counterexample}.
	\end{proof}
	
	\subsection{Zero-moment profiles}
	
	For nonintegral $2(t-s)$, causal fractional differentiation of a compactly supported profile generally produces an algebraic tail. To preserve stretched-exponential decay, we use profiles whose integer moments all vanish. The starting construction is classical in the Stieltjes moment problem; see \cite[Section 4]{Kleiber14} and the references there. We give the moment calculation and then smooth the profile while preserving its moments.
	
	We say that a function $F$ on $\R$ has all integer moments zero if $t^mF\in L^1(0,\infty)$ for every $m\in\Nzero$ and
	\begin{equation}\label{eq:profile-moments}
		\int_0^\infty t^mF(t)\,dt=0\quad\text{for every }m\in\Nzero.
	\end{equation}
	
	\begin{lemma}\label{lem:zero-moments}
		For $0<\gamma<1/2$ and $R>0$, there is a nonzero real $F\in C^\infty(\R)$, zero for $t\le R+1$, such that \eqref{eq:profile-moments} holds and $|F^{(k)}(t)|\le C_ke^{-a t^\gamma}$ for all $t\ge0$ and $k\in\Nzero$, with $a>0$ independent of $k$.
	\end{lemma}
	
	\begin{proof}
		Set $f_\gamma(t)=e^{-t^\gamma}\sin(t^\gamma\tan(\pi\gamma))$ for $t>0$ and $f_\gamma(t)=0$ for $t\le0$. This function is nonzero, locally integrable, and bounded in absolute value by $e^{-t^\gamma}$ on $(0,\infty)$. Put $\theta=\pi\gamma\in(0,\pi/2)$. The substitution $z=t^\gamma$ and the Gamma integral give
		\begin{equation}\label{eq:zero-moment-computation}
			\begin{split}
				\int_0^\infty t^mf_\gamma(t)\,dt
				&=\frac{\Gamma((m+1)/\gamma)}{\gamma}\Ima(1-\mathsf{i}\tan\theta)^{-(m+1)/\gamma}\\
				&=\frac{\Gamma((m+1)/\gamma)}{\gamma}(\cos\theta)^{(m+1)/\gamma}\sin(\pi(m+1))=0,
			\end{split}
		\end{equation}
		where we used $1-\mathsf{i}\tan\theta=e^{-\mathsf{i}\theta}/\cos\theta$ and $\theta(m+1)/\gamma=\pi(m+1)$.
		
		Choose a nonnegative $\rho\in C_c^\infty((R+1,R+2))$ with integral one and set $F=\rho*f_\gamma$ on $\R$. Since $F^{(k)}=\rho^{(k)}*f_\gamma$, the function $F$ is smooth and vanishes for $t\le R+1$. For large $t$ and $s\in\supp\rho$, we have $t-s\ge t/2$, so $|F^{(k)}(t)|\le\norm{\rho^{(k)}}_1e^{-(t/2)^\gamma}$. Increasing $C_k$ to cover bounded intervals gives the required bounds with one rate independent of $k$.
		
		The convolution also preserves the zero moments. Absolute convergence permits Fubini's theorem. Expanding $(s+t)^m$ and using \eqref{eq:zero-moment-computation} gives
		\begin{equation}\label{eq:convolution-moments}
			\int_0^\infty r^mF(r)\,dr=\int_0^\infty\rho(s)\int_0^\infty(s+t)^mf_\gamma(t)\,dt\,ds=0.
		\end{equation}
		
		The convolution identity for Laplace transforms is
		\begin{equation}\label{eq:convolution-Laplace}
			\mathcal LF(\kappa)=(\mathcal L\rho)(\kappa)(\mathcal Lf_\gamma)(\kappa),\quad\kappa>0.
		\end{equation}
		If $F=0$, positivity of $\mathcal L\rho(\kappa)$ and injectivity of the Laplace transform imply $f_\gamma=0$, a contradiction.
	\end{proof}
	
	\subsection{Decay under causal fractional differentiation}
	
	The next lemma shows how the zero moments remove the algebraic tail. The main point is that the resulting exponential rate can be chosen independently of the number of derivatives.
	
	\begin{lemma}\label{lem:fractional-tail}
		Let $0<\gamma<1$ and let $F$ be smooth, vanish on a neighborhood of $(-\infty,0]$, and satisfy \eqref{eq:profile-moments}. Suppose that $|F^{(k)}(t)|\le C_ke^{-a t^\gamma}$ for all $t\ge0$ and $k\in\Nzero$, with $a>0$ independent of $k$. For every fixed $\beta>0$, the function $D_+^\beta F$ defined by \eqref{eq:causal-derivative} is smooth, vanishes on every interval $(-\infty,R_0]$ on which $F$ vanishes, and satisfies
		\begin{equation}\label{eq:fractional-tail}
			|(D_+^\beta F)^{(k)}(r)|\le C'_ke^{-a'r^\gamma},\quad r\ge0,\quad k\in\Nzero,
		\end{equation}
		where $a'>0$ may be chosen to depend only on $a,\gamma,\beta$, and is independent of $k$. The constants $C'_k$ may also depend on $F$ and $k$.
	\end{lemma}
	
	\begin{proof}
		For integral $\beta$, the conclusion follows directly from the assumed derivative bounds. Otherwise, write $\beta=m-\delta$, with $m\in\N$ and $0<\delta<1$, and let $I_+^\delta$ denote convolution with $r_+^{\delta-1}/\Gamma(\delta)$. Since $F$ and all its derivatives vanish near zero, differentiation can be taken on $F$ in this convolution, giving
		\begin{equation}\label{eq:causal-derivative-derivatives}
			(D_+^\beta F)^{(k)}(r)=I_+^\delta F^{(m+k)}(r)=\frac1{\Gamma(\delta)}\int_0^r(r-t)^{\delta-1}F^{(m+k)}(t)\,dt,\quad r\ge0.
		\end{equation}
		For $r<0$, the convolution vanishes because both factors are supported in $[0,\infty)$. The kernel is locally integrable, so this formula gives continuity of every derivative and proves smoothness. The same formula shows that $D_+^\beta F$ vanishes on every interval $(-\infty,R_0]$ on which $F$ vanishes.
		
		Fix $k$ and put $g=F^{(m+k)}$. Its decay estimate has the same rate $a$, while its prefactor may depend on $k$ and the fixed order $\beta$; we denote this prefactor by $C_k$. Repeated integration by parts shows that all integer moments of $g$ vanish:
		\begin{equation}\label{eq:jth-moment}
			\int_0^\infty t^\ell g(t)\,dt=0,\quad \ell\in\Nzero.
		\end{equation}
		All boundary terms vanish because $F$ and its derivatives vanish near zero and decay at infinity. If $\ell<m+k$, differentiating $t^\ell$ repeatedly gives zero. If $\ell\ge m+k$, the remaining integral is a constant multiple of $\int_0^\infty t^{\ell-m-k}F(t)\,dt$, which vanishes by \eqref{eq:profile-moments}. We split
		\begin{equation}\label{eq:tail-split}
			\Gamma(\delta)I_+^\delta g(r)=\int_0^{r/2}(r-t)^{\delta-1}g(t)\,dt+\int_{r/2}^r(r-t)^{\delta-1}g(t)\,dt.
		\end{equation}
		The second integral is bounded by $C_kr^\delta e^{-a(r/2)^\gamma}$, which has the required decay after decreasing the exponential rate.
		
		To estimate the first integral, expand $(1-x)^{\delta-1}$ for $0\le x\le1/2$:
		\[
		(1-x)^{\delta-1}=\sum_{\ell=0}^{M-1}d_\ell x^\ell+E_M(x),\quad d_0=1,\quad d_\ell=\prod_{q=1}^\ell(1-\delta/q)\quad(\ell\ge1),
		\]
		where the integer $M$ will be chosen later. Since $0<d_\ell\le1$, the tail of the series satisfies
		\begin{equation}\label{eq:binomial-tail-bound}
			|E_M(x)|\le\frac{x^M}{1-x}\le2x^M,\quad0\le x\le1/2.
		\end{equation}
		For each $\ell=0,\ldots,M-1$, \eqref{eq:jth-moment} gives $\int_0^{r/2}t^\ell g(t)\,dt=-\int_{r/2}^\infty t^\ell g(t)\,dt$. Substitution into the finite expansion and \eqref{eq:binomial-tail-bound} give
		\begin{equation}\label{eq:tail-decomposition}
			\int_0^{r/2}(r-t)^{\delta-1}g(t)\,dt=-r^{\delta-1}\sum_{\ell=0}^{M-1}d_\ell r^{-\ell}\int_{r/2}^\infty t^\ell g(t)\,dt+\mathcal E_M(r),
		\end{equation}
		and
		\begin{equation}\label{eq:fractional-tail-remainder-general}
			\begin{aligned}
				|\mathcal E_M(r)|\le2r^{\delta-1-M}\int_0^\infty t^M|g(t)|\,dt\le\frac{2C_k}{\gamma}r^{\delta-1-M}a^{-(M+1)/\gamma}\Gamma((M+1)/\gamma).
			\end{aligned}
		\end{equation}
		
		Stirling's inequality gives constants $C,C_0>0$, depending only on $a,\gamma$, such that
		\[
		a^{-(M+1)/\gamma}\Gamma((M+1)/\gamma)\le C(C_0M^{1/\gamma})^M,\quad M\ge1.
		\]
		The fixed powers of $M$ in Stirling's inequality have been included in $C_0^M$. Choose $\eps_0>0$ so that $C_0\eps_0^{1/\gamma}\le1/2$, and take $M=\lfloor\eps_0r^\gamma\rfloor$ for large $r$. Then
		\[
		|\mathcal E_M(r)|\le C_kr^{\delta-1}(C_0M^{1/\gamma}/r)^M\le C_kr^{\delta-1}2^{-M}.
		\]
		For sufficiently large $r$, we have $M\ge\eps_0r^\gamma/2$, so this remainder is bounded by $C_ke^{-c r^\gamma}$ with $c>0$ independent of $k$.
		
		It remains to estimate the finite sum in \eqref{eq:tail-decomposition}. Decrease $\eps_0$ further so that $\eps_0\le ae\gamma/2$. For $t\ge r$ and $0\le\ell<M$, the identity $\sup_{q\ge1}(\log q)/q^\gamma=1/(e\gamma)$ gives $\ell\log(t/r)\le\eps_0r^\gamma\log(t/r)\le\frac a2t^\gamma$. For $r/2\le t\le r$, we instead use $(t/r)^\ell\le1$. In both ranges, we have
		\[
		(t/r)^\ell e^{-a t^\gamma}\le e^{-a t^\gamma/2},\quad t\ge r/2.
		\]
		Since $d_\ell\le1$, the absolute value of the finite sum, including its factor $r^{\delta-1}$, is at most $C_kMr^{\delta-1}\int_{r/2}^\infty e^{-a t^\gamma/2}\,dt$. Writing the integrand as $e^{-a t^\gamma/4}e^{-a t^\gamma/4}$ shows that this integral is bounded by $Ce^{-(a/4)(r/2)^\gamma}$. The factor $Mr^{\delta-1}$ is at most a fixed power of $r$ and can be absorbed by decreasing the exponential rate.
		
		The choices of $\eps_0$ depend only on $a$, $\gamma$, and $\beta$. The preceding estimates give \eqref{eq:fractional-tail} for large $r$ with $a'>0$ depending only on these parameters, independently of $k$. Increasing $C'_k$ on a bounded interval proves the estimate for all $r\ge0$.
	\end{proof}
	
	\begin{proof}[Proof of Theorem~\ref{thm:sharpness}\ref{item_2_thm_sharp} in odd dimensions]
		Fix $0<\gamma<1/2$ and choose a nonzero real $F$ from Lemma~\ref{lem:zero-moments}. Put $\beta=2(t-s)$ and $G=D_+^\beta F$. Lemma~\ref{lem:fractional-tail} gives $|F^{(k)}(r)|+|G^{(k)}(r)|\le C_ke^{-a'r^\gamma}$ for $r\ge0$ and $k\in\Nzero$, where $a'>0$ is independent of $k$. Both profiles belong to $\cS(\R)$ and vanish for $r\le R+1$. Moreover, \eqref{eq:Laplace-derivative} gives $\mathcal LG(\kappa)=\kappa^\beta\mathcal LF(\kappa)$ for $\kappa>0$. If $G=0$, injectivity of the Laplace transform would imply $F=0$, so $G$ is nonzero.
		
		Define the three-dimensional inputs
		\[
		u_3(x)=-\sin(\pi t)\frac{G(|x|)}{|x|},\quad v_3(x)=\sin(\pi s)\frac{F(|x|)}{|x|},\quad x\in\R^3\setminus\{0\},
		\]
		and set both functions equal to zero at the origin. They are nonzero because $F,G$ and the two sine factors are nonzero. Since the profiles vanish for $r\le R+1$, division by $|x|$ causes no singularity, and both functions vanish on $B_{R+1}\subset\R^3$. At infinity, their spatial derivatives are finite sums of profile derivatives multiplied by nonpositive powers of $|x|$ and bounded angular factors. The bounds on $F$ and $G$ show that $u_3,v_3\in\cS(\R^3)$ and satisfy \eqref{eq:counterexample-decay}, with a decay rate independent of the derivative order. Proposition~\ref{prop:radial-cancellation} gives, for every $\kappa>0$,
		\[
		W_{3,\kappa}:=(-\Delta_3+\kappa^2)^{-1}(\sin(\pi s)\kappa^{2s}u_3+\sin(\pi t)\kappa^{2t}v_3)=0
		\quad\text{on }B_R\subset\R^3.
		\]
		
		Let $n=3+2q$ with $q\in\Nzero$, and set $u=\mathcal T^q u_3$ and $v=\mathcal T^q v_3$, using the same convention for $\mathcal T^q$ as in the proof of Theorem~\ref{thm:sharpness}\ref{item_1_thm_sharp}. Lemma~\ref{lem:dimension} ensures that $u,v$ are nonzero radial Schwartz functions on $\R^n$, vanish on $B_R$, and satisfy \eqref{eq:counterexample-decay} with a rate independent of the derivative order. For $\kappa>0$, set $f_\kappa=\sin(\pi s)\kappa^{2s}u+\sin(\pi t)\kappa^{2t}v$. Linearity and repeated use of \eqref{eq:dimension-identity} give
		\[
		(-\Delta_n+\kappa^2)^{-1}f_\kappa
		=\mathcal T^q W_{3,\kappa}=0
		\quad\text{on }B_R\subset\R^n,
		\]
		where the last equality uses preservation of vanishing on a ball.
		
		Since $\supp(A-I)\Subset B_R$, Lemma~\ref{lem:transfer}, applied in $\R^n$, gives $(L+\kappa^2)^{-1}f_\kappa=0$ on $B_R$ for every $\kappa>0$. Together with $u=v=0$ on $B_R$, the converse part of Proposition~\ref{prop:resolvent} proves \eqref{eq:counterexample}. The proof of the even-dimensional case will be given in the next subsection.
	\end{proof}

	\subsection{Abel integrals in even dimensions}
	
	In dimension two, the radial resolvent involves an integral of a smooth source profile $f$ with the following form. Setting $x=\tau^2$ and substituting $y=r^2$ gives
	\[
	\int_0^\tau\frac{r f(r)}{\sqrt{\tau^2-r^2}}\,dr
	=\frac12\int_0^x\frac{f(\sqrt y)}{\sqrt{x-y}}\,dy,\quad \tau>0.
	\]
	The integral on the right has the classical Abel kernel $(x-y)^{-1/2}$ and equals $\sqrt\pi/2$ times the causal fractional integral of order $1/2$ of $y\mapsto f(\sqrt y)$. We will invert this integral to obtain the radial resolvent representation in Lemma~\ref{lem:even-radial-resolvent}. The inversion uses a causal derivative of order $1/2$, so we first verify that causal fractional differentiation preserves the zero moments.
	
	\begin{lemma}\label{lem:fractional-moments-preserved}
		Let $F$ satisfy the hypotheses of Lemma~\ref{lem:fractional-tail}, and let $\beta>0$. Then $G=D_+^\beta F$ has all integer moments zero:
		\begin{equation}\label{eq:G-zero-moments}
			\int_0^\infty t^mG(t)\,dt=0,\quad m\in\Nzero.
		\end{equation}
	\end{lemma}
	
	\begin{proof}
		The decay of $F$ and Lemma~\ref{lem:fractional-tail} give $t^mF,t^mG\in L^1(0,\infty)$ for every $m\in\Nzero$. Differentiation under the integral and dominated convergence show that both Laplace transforms are smooth up to $\kappa=0$ from the right. Since $F$ satisfies \eqref{eq:profile-moments}, we have $(\mathcal LF)^{(m)}(0)=(-1)^m\int_0^\infty t^mF(t)\,dt=0$ for every $m\in\Nzero$. Taylor's formula applied to each derivative gives
		\begin{equation}\label{eq:Laplace-flatness}
			(\mathcal LF)^{(\ell)}(\kappa)=O(\kappa^M),\quad\kappa\to0^+,\quad \ell,M\in\Nzero.
		\end{equation}
		
		By \eqref{eq:Laplace-derivative}, $\mathcal LG(\kappa)=\kappa^\beta\mathcal LF(\kappa)$ for $\kappa>0$. Fix $m\in\Nzero$. Leibniz' rule gives
		\begin{equation}\label{eq:Laplace-fractional-product}
			(\mathcal LG)^{(m)}(\kappa)=\sum_{\ell=0}^m\binom m\ell(\partial_\kappa^{m-\ell}\kappa^\beta)(\mathcal LF)^{(\ell)}(\kappa),\quad\kappa>0.
		\end{equation}
		Each factor $\partial_\kappa^{m-\ell}\kappa^\beta$ is a constant multiple of $\kappa^{\beta-m+\ell}$. Choosing $M>m-\beta$ in \eqref{eq:Laplace-flatness} shows that every term tends to zero as $\kappa\to0^+$. Since $t^mG\in L^1(0,\infty)$, dominated convergence gives
		\begin{equation}\label{eq:Laplace-moment-limit}
			\int_0^\infty t^mG(t)\,dt
			=(-1)^m\lim_{\kappa\to0^+}(\mathcal LG)^{(m)}(\kappa)=0.
		\end{equation}
		This proves \eqref{eq:G-zero-moments} for every $m\in\Nzero$.
	\end{proof}
	
	Let $R_0>0$. For a smooth function $F$ vanishing for $t\le R_0$, define
	\begin{equation}\label{eq:even-Abel-inverse}
		H_F(x)=
		\begin{cases}F(\sqrt{x}),&x>0,\\0,&x\le0,\end{cases}
		\quad q_F(x)=\frac{2}{\sqrt\pi}D_{+,x}^{1/2}H_F(x),\quad f_F(r)=q_F(r^2),\quad r\ge0.
	\end{equation}
	Here $D_{+,x}^{1/2}$ acts in the variable $x$. Since $H_F$ vanishes on $(-\infty,R_0^2]$, the change of variables causes no loss of smoothness at zero.
	
	\begin{lemma}\label{lem:even-Abel-profiles}
		Given $R_0>0$, suppose that $F$ is real, smooth, vanishes for $t\le R_0$, has all integer moments zero, and satisfies $|F^{(k)}(t)|\le C_ke^{-a t^\gamma}$ for all $t\ge0$ and $k\in\Nzero$, with $a>0$ independent of $k$ and $0<\gamma<1$. The functions in \eqref{eq:even-Abel-inverse} have the following properties:
		\begin{enumerate}[\rm(i)]
			\item\label{item_1_lem_even_Abel} $H_F,q_F\in\cS(\R)$, both vanish on $(-\infty,R_0^2]$, and $f_F(|x|)\in\cS(\R^2)$ vanishes on $B_{R_0}$.
			\item\label{item_2_lem_even_Abel} There is $a'>0$, independent of the order of differentiation, such that all spatial derivatives of $f_F(|x|)$ satisfy bounds $C_\alpha e^{-a'|x|^\gamma}$.
			\item\label{item_3_lem_even_Abel} For every $\tau>0$,
			\begin{equation}\label{eq:Abel-profile-identity}
				\int_0^\tau\frac{r f_F(r)}{\sqrt{\tau^2-r^2}}\,dr=F(\tau).
			\end{equation}
			\item\label{item_4_lem_even_Abel} If $F\ne0$, then $f_F\ne0$.
		\end{enumerate}
	\end{lemma}
	
	\begin{proof}
		For $x>0$, put $r=\sqrt x$. Since $\partial_x=(2r)^{-1}\partial_r$, repeated differentiation gives
		\begin{equation}\label{eq:Abel-chain-derivatives}
			H_F^{(k)}(x)=\sum_{\ell=1}^k c_{k,\ell}r^{\ell-2k}F^{(\ell)}(r),\quad r=\sqrt x,\quad k\ge1,
		\end{equation}
		where the constants $c_{k,\ell}$ depend only on $k$ and $\ell$. For $x\ge1$, we have $r\ge1$ and $\ell-2k\le0$, so $|H_F^{(k)}(x)|\le\sum_{\ell=1}^k|c_{k,\ell}|C_\ell e^{-ar^\gamma}=C'_ke^{-a x^{\gamma/2}}$. The case $k=0$ follows directly from the bound on $F$. Since $H_F$ is smooth and vanishes near zero, each derivative is bounded on $[0,1]$. Increasing $C'_k$ to cover this interval gives the same estimate for all $x\ge0$, with the original rate $a$ independent of $k$. In particular, $H_F\in\cS(\R)$. Changing variables $x=t^2$ gives
		\begin{equation}\label{eq:Abel-input-moments}
			\int_0^\infty x^mH_F(x)\,dx=2\int_0^\infty t^{2m+1}F(t)\,dt=0,\quad m\in\Nzero.
		\end{equation}
		Lemma~\ref{lem:fractional-tail}, with decay exponent $\gamma/2$ and derivative order $1/2$, gives $q_F\in\cS(\R)$ and
		\[
		|q_F^{(k)}(x)|\le C'_ke^{-a_1x^{\gamma/2}},\quad x\ge0,
		\]
		with $a_1>0$ independent of $k$. It also vanishes on $(-\infty,R_0^2]$.
		
		The radial function in the conclusion can be written as $f_F(|x|)=q_F(|x|^2)$. Each spatial derivative is a finite sum of derivatives of $q_F$ at $|x|^2$ multiplied by polynomials in $x$. For each derivative order, these polynomial factors are bounded by a constant times $e^{(a_1/2)|x|^\gamma}$. Taking $a'=a_1/2$ proves the bound in \ref{item_2_lem_even_Abel} for every multi-index $\alpha$, and also proves the Schwartz property. Since $q_F$ vanishes on $(-\infty,R_0^2]$, the radial function is zero on $B_{R_0}$. This proves \ref{item_1_lem_even_Abel} and \ref{item_2_lem_even_Abel}.
		
		To prove \ref{item_3_lem_even_Abel}, let $I_{+,x}^{1/2}$ denote the causal fractional integral in $x$. The Gamma integral gives the Laplace transform multiplier $\kappa^{-1/2}$ for $I_{+,x}^{1/2}$, while \eqref{eq:Laplace-derivative} gives the multiplier $\kappa^{1/2}$ for $D_{+,x}^{1/2}$. Since $H_F$ vanishes near zero, there are no initial-value terms, and injectivity of the Laplace transform gives $I_{+,x}^{1/2}D_{+,x}^{1/2}H_F=H_F$. For $x>0$, the definition of $q_F$ gives
		\begin{equation}\label{eq:Abel-causal-inversion}
			\frac12\int_0^x\frac{q_F(y)}{\sqrt{x-y}}\,dy=\frac{\sqrt\pi}{2}I_{+,x}^{1/2}q_F(x)=H_F(x).
		\end{equation}
		Taking $x=\tau^2$ and substituting $y=r^2$ proves \eqref{eq:Abel-profile-identity}. If $f_F=0$, the left-hand side of that identity is zero for every $\tau>0$, so $F=0$. This proves \ref{item_4_lem_even_Abel}.
	\end{proof}
	
	\begin{lemma}\label{lem:even-radial-resolvent}
		Let $F$ satisfy the hypotheses of Lemma~\ref{lem:even-Abel-profiles} for some $R_0>0$, let $f_F$ be defined by \eqref{eq:even-Abel-inverse}, and set $v_F(x)=f_F(|x|)$ on $\R^2$. Let $I_0$ denote the modified Bessel function of the first kind of order zero, with $I_0(0)=1$, and write $\mathcal LF(\kappa)=\int_0^\infty e^{-\kappa t}F(t)\,dt$ for the Laplace transform of $F$. For every $\kappa>0$,
		\begin{equation}\label{eq:even-radial-resolvent}
			((-\Delta+\kappa^2)^{-1}v_F)(x)=I_0(\kappa|x|)\mathcal LF(\kappa),\quad x\in B_{R_0}\subset\R^2.
		\end{equation}
	\end{lemma}
	
	\begin{proof}
		Fix $\kappa>0$ and set $w=(-\Delta+\kappa^2)^{-1}v_F$. By Lemma~\ref{lem:even-Abel-profiles}, $v_F$ is a radial Schwartz function. The free resolvent preserves this class, so we may write $w(x)=W(|x|)$, where $W$ is smooth, bounded near zero, and rapidly decreasing together with its derivatives. The radial resolvent equation is $-W''(r)-r^{-1}W'(r)+\kappa^2W(r)=f_F(r)$ for $r>0$.
		
		Let $K_0$ denote the modified Bessel function of the second kind of order zero. The functions $r\mapsto I_0(\kappa r)$ and $r\mapsto K_0(\kappa r)$ are linearly independent solutions of $-y''-r^{-1}y'+\kappa^2y=0$ on $(0,\infty)$. Since $f_F(r)=0$ for $0<r<R_0$, the profile $W$ satisfies this homogeneous equation on $(0,R_0)$. There are constants $c_\kappa,d_\kappa$ such that $W(r)=c_\kappa I_0(\kappa r)+d_\kappa K_0(\kappa r)$ for $0<r<R_0$. By \eqref{eq:K0-series}, $I_0(\kappa r)\to1$ and $K_0(\kappa r)\to+\infty$ as $r\to0$. Boundedness of $W$ near zero forces $d_\kappa=0$, giving $W(r)=c_\kappa I_0(\kappa r)$ for $0<r<R_0$.
		
		To determine $c_\kappa$, put $Z(r)=K_0(\kappa r)$, which satisfies $-Z''(r)-r^{-1}Z'(r)+\kappa^2Z(r)=0$ for $r>0$. Multiplying the equation for $W$ by $Z$ and subtracting the equation for $Z$ multiplied by $W$ gives
		\[
		\frac{d}{dr}(r(W'(r)Z(r)-W(r)Z'(r)))=-r f_F(r)Z(r).
		\]
		The quantity $r(W'Z-WZ')$ tends to zero as $r\to\infty$, by the rapid decay of $W,W'$ and \eqref{eq:Bessel-large} with $d=0$. For $r>0$, integration from $r$ to infinity gives
		\begin{equation}\label{eq:even-radial-Green-identity}
			r(W'(r)Z(r)-W(r)Z'(r))=\int_r^\infty tK_0(\kappa t)f_F(t)\,dt.
		\end{equation}
		Put $P(r)=I_0(\kappa r)$ and recall that $Z(r)=K_0(\kappa r)$. The homogeneous equations for $P$ and $Z$ give $(r(P'Z-PZ'))'=0$ on $(0,\infty)$, so $r(P'Z-PZ')$ is constant. The series \eqref{eq:K0-series} and their derivatives give $P(r)=1+O(r^2)$, $P'(r)=O(r)$, $Z(r)=O(1+|\log r|)$, and $rZ'(r)\to-1$ as $r\to0$. Taking the limit in $r(P'Z-PZ')$ gives $0-1(-1)=1$, proving the Wronskian identity $P'(r)Z(r)-P(r)Z'(r)=1/r$ for every $r>0$.
		
		For $0<r<R_0$, substituting $W=c_\kappa P$ gives $r(W'Z-WZ')=c_\kappa$. Equation~\eqref{eq:even-radial-Green-identity} yields $c_\kappa=\int_r^\infty tK_0(\kappa t)f_F(t)\,dt$. Since $f_F=0$ on $[0,R_0]$, this gives
		\begin{equation}\label{eq:even-radial-Green-formula}
			W(r)=I_0(\kappa r)\int_0^\infty tK_0(\kappa t)f_F(t)\,dt,\quad0<r<R_0.
		\end{equation}
		
		For $t>0$, the integral representation \cite[Eq.~(10.32.9)]{NIST10}, with order zero, and the substitution $\tau=t\cosh\eta$ give
		\[
		K_0(\kappa t)=\int_0^\infty e^{-\kappa t\cosh\eta}\,d\eta
		=\int_t^\infty\frac{e^{-\kappa\tau}}{\sqrt{\tau^2-t^2}}\,d\tau.
		\]
		The interchange of integrals below is justified by absolute convergence. Indeed, Tonelli's theorem applied with $|f_F|$ identifies the corresponding nonnegative double integral with $\int_0^\infty K_0(\kappa t)|f_F(t)|t\,dt$, which is finite because $f_F$ is Schwartz and vanishes near zero, while $K_0(\kappa t)$ decays exponentially at infinity. Fubini's theorem and \eqref{eq:Abel-profile-identity} give
		\begin{equation}\label{eq:even-radial-Laplace-calculation}
			\begin{aligned}
				c_\kappa=\int_0^\infty e^{-\kappa\tau}\Big(\int_0^\tau\frac{t f_F(t)}{\sqrt{\tau^2-t^2}}\,dt\Big)d\tau=\int_0^\infty e^{-\kappa\tau}F(\tau)\,d\tau=\mathcal LF(\kappa).
			\end{aligned}
		\end{equation}
		This proves the asserted identity for $0<|x|<R_0$. Continuity of $w$ and $I_0$, together with $I_0(0)=1$, gives the identity at $x=0$.
	\end{proof}
	
	\begin{proof}[Completion of the proof of Theorem~\ref{thm:sharpness}\ref{item_2_thm_sharp}]
		It remains to treat the even-dimensional case. Fix $0<\gamma<1/2$, put $\beta=2(t-s)$, and choose $F$ and $G=D_+^\beta F$ as in the odd-dimensional proof. Both profiles are nonzero real Schwartz functions, vanish for $r\le R+1$, and satisfy the derivative bounds required in Lemma~\ref{lem:even-Abel-profiles}. All integer moments of $F$ vanish by construction, and Lemma~\ref{lem:fractional-moments-preserved} gives the same property for $G$.
		
		Apply \eqref{eq:even-Abel-inverse} to $F$ and $G$, with $R_0=R+1$. In dimension two, define
		\begin{equation}\label{eq:even-counterexample-inputs}
			u=-\sin(\pi t)f_G(|x|),\quad v=\sin(\pi s)f_F(|x|).
		\end{equation}
		Since $s,t$ are nonintegers, both scalar coefficients are nonzero. Lemma~\ref{lem:even-Abel-profiles} shows that $u,v$ are nonzero real Schwartz functions, zero on $B_R$, and satisfy \eqref{eq:counterexample-decay} with $c>0$ independent of the multi-index $\alpha$.
		
		Since $G=D_+^\beta F$, \eqref{eq:Laplace-derivative} gives $\mathcal LG(\kappa)=\kappa^\beta\mathcal LF(\kappa)$ for $\kappa>0$. Applying Lemma~\ref{lem:even-radial-resolvent} to $F$ and $G$ gives
		\[
		(-\Delta+\kappa^2)^{-1}f_G(|x|)=\kappa^\beta(-\Delta+\kappa^2)^{-1}f_F(|x|)\quad\text{on }B_R,\quad\kappa>0.
		\]
		Substitution of \eqref{eq:even-counterexample-inputs} and $2s+\beta=2t$ gives
		\begin{equation}\label{eq:even-resolvent-cancellation}
			(\sin(\pi s)\kappa^{2s}(-\Delta+\kappa^2)^{-1}u+\sin(\pi t)\kappa^{2t}(-\Delta+\kappa^2)^{-1}v)|_{B_R}=0.
		\end{equation}
		Let $n=2\nu$ with $\nu\in\N$, and define $u_n=\mathcal T^{\nu-1}u$ and $v_n=\mathcal T^{\nu-1}v$ on $\R^n$, where $u,v$ are the two-dimensional functions in \eqref{eq:even-counterexample-inputs} and $\mathcal T^0$ denotes the identity. Lemma~\ref{lem:dimension} shows that $u_n,v_n$ are nonzero real Schwartz functions, vanish on $B_R\subset\R^n$, and satisfy \eqref{eq:counterexample-decay} with a rate independent of the derivative order.
		
		For $\kappa>0$, set $f_\kappa=\sin(\pi s)\kappa^{2s}u_n+\sin(\pi t)\kappa^{2t}v_n$. Repeated use of \eqref{eq:dimension-identity} in \eqref{eq:even-resolvent-cancellation}, together with preservation of vanishing on a ball, gives $(-\Delta_n+\kappa^2)^{-1}f_\kappa=0$ on $B_R\subset\R^n$. Since $\supp(A-I)\Subset B_R$, Lemma~\ref{lem:transfer} gives $(L+\kappa^2)^{-1}f_\kappa=0$ on $B_R$ for every $\kappa>0$. Together with $u_n=v_n=0$ on $B_R$, the converse part of Proposition~\ref{prop:resolvent}, with $\lambda=\kappa^2$, proves \eqref{eq:counterexample} for the pair $u_n,v_n$.
	\end{proof}

	\subsection{Sharpness of the logarithmic decay condition}\label{subsec:logarithmic-sharpness}
	
	We prove Proposition~\ref{prop:logarithmic-sharpness}. We first construct a nonzero profile with all integer moments zero and the required logarithmic decay.
	
	\begin{lemma}\label{lem:logarithmic-zero-moments}
		For every $p>1$ and $R>0$, there is a nonzero real $F\in C^\infty(\R)$, zero for $t\le R+1$, such that \eqref{eq:profile-moments} holds and $|F^{(k)}(t)|\le C_ke^{-a w_p(t)}$ for $t\ge0$ and $k\in\Nzero$, with $a>0$ independent of $k$.
	\end{lemma}
	
	\begin{proof}
		Put $q=p-1>0$ and choose $E>\exp(2q+2)$. For $\Rea w>0$, define
		\[
		\Psi(w)=\frac{w}{(\log(E+w))^q},
		\]
		using principal branches. Since $\Rea\log(E+w)\ge\log E>0$, this is holomorphic in the right half-plane and extends across every finite point of its boundary. We claim that
		\begin{equation}\label{eq:logarithmic-phase-bound}
			\Rea\Psi(w)\ge c\frac{|w|}{\log^p(E+|w|)}
		\end{equation}
		for $\Rea w\ge0$ and all sufficiently large $|w|$.
		
		Since $\Psi(\overline w)=\overline{\Psi(w)}$, it suffices to consider $w=\rho e^{\mathsf{i}\theta}$ with $\rho=|w|$ and $0\le\theta\le\pi/2$. Write $\log(E+w)=a+\mathsf{i}b$, where $a=\log|E+w|$ and $b=\arg(E+w)$, and put $\vartheta=\arg(a+\mathsf{i}b)=\arctan(b/a)$. Adding the positive real number $E$ to $w$ can only decrease its argument, so $0\le b\le\theta$. Also, $a\ge\log E$, and the inequality $\arctan x\le x$ for $x\ge0$ gives $0\le q\vartheta\le\frac{qb}{a}\le\frac{q\theta}{\log E}\le\frac{\theta}{2}$. In particular, $0\le\theta-q\vartheta\le\theta$.
		
		For $\rho\ge2E$, we have $\rho\le|E+w|\le E+\rho\le3\rho/2$, which gives $\log(E+\rho)-\log(3/2)\le a\le\log(E+\rho)$. Since $0\le b\le\pi/2$ and $a\le|\log(E+w)|=\sqrt{a^2+b^2}\le a+\pi/2$, both $a$ and $|\log(E+w)|$ are comparable to $\log(E+\rho)$, with constants independent of $\rho$ and $\theta$.
		
		If $0\le\theta\le\pi/4$, then $0\le\theta-q\vartheta\le\pi/4$, so $\cos(\theta-q\vartheta)\ge2^{-1/2}$. If $\pi/4\le\theta\le\pi/2$, then $\sin\theta\ge2^{-1/2}$ and $E\le\rho/2$ give $\tan b=\frac{\rho\sin\theta}{E+\rho\cos\theta}\ge\frac{\rho/\sqrt2}{3\rho/2}=\frac{2}{3\sqrt2}$. Hence, $b\ge b_0:=\arctan(2/(3\sqrt2))>0$. Since $b/a\le\pi/(2\log E)<1$, the inequality $\arctan x\ge x/2$ for $0\le x\le1$ gives $\vartheta\ge b/(2a)\ge b_0/(2\log(E+\rho))$. Using $\theta\le\pi/2$, $0\le q\vartheta\le\pi/4$, and $\sin x\ge2x/\pi$ for $0\le x\le\pi/2$, we obtain $\cos(\theta-q\vartheta)\ge\sin(q\vartheta)\ge\frac{2q\vartheta}{\pi}\ge\frac{qb_0}{\pi\log(E+\rho)}$.
		
		Both angular ranges give $\cos(\theta-q\vartheta)\ge c/\log(E+\rho)$. Taking the real part of $\Psi(w)=\rho|\log(E+w)|^{-q}e^{\mathsf{i}(\theta-q\vartheta)}$ and using the preceding comparison yields
		\[
		\Rea\Psi(w)=\rho|\log(E+w)|^{-q}\cos(\theta-q\vartheta)\ge c\frac{\rho}{\log^{q+1}(E+\rho)},
		\]
		where $c>0$ may depend on $q$ and $E$, but is independent of $\rho$ and $\theta$. Since $q+1=p$, this proves \eqref{eq:logarithmic-phase-bound}.
		
		For $z\in\C\setminus[0,\infty)$, set
		\[
		H(z)=\exp(-\Psi(\sqrt{-z})),
		\]
		where $\sqrt{-z}$ has positive real part. This function is holomorphic on $\C\setminus[0,\infty)$ and satisfies $H(\overline z)=\overline{H(z)}$. Recall that $w_p(r)=\sqrt r/\log^p(e+r)$ for $r\ge0$. Since $|\sqrt{-z}|=\sqrt{|z|}$, applying \eqref{eq:logarithmic-phase-bound} with $w=\sqrt{-z}$ and comparing $\log(E+\sqrt{|z|})$ with $\log(e+|z|)$ gives $|H(z)|\le e^{-c w_p(|z|)}$ for sufficiently large $|z|$. The continuity of $\Psi$ on the closed right half-plane also shows that $H$ is bounded when $|z|$ is bounded. Increasing a prefactor $C$ gives
		\begin{equation}\label{eq:logarithmic-slit-decay}
			|H(z)|\le C e^{-c w_p(|z|)},\quad z\in\C\setminus[0,\infty),
		\end{equation}
		where $C,c>0$ are independent of $z$.
		
		For $t>0$, the boundary values from above and below the positive real axis are
		\[
		H(t\pm\mathsf{i}0):=\lim_{\varepsilon\downarrow 0}H(t\pm\mathsf{i}\varepsilon)=\exp(-\Psi(\mp\mathsf{i}\sqrt t)).
		\]
		These limits exist by continuity of $\Psi$ at $\mp\mathsf{i}\sqrt t$. Passing to the limit in \eqref{eq:logarithmic-slit-decay} gives $|H(t\pm\mathsf{i}0)|\le C e^{-c w_p(t)}$ for every $t>0$, with the same constants $C,c$. Define $h(t)=\Ima H(t+\mathsf{i}0)$ for $t>0$, and set $h(t)=0$ for $t\le0$. The expansion at zero is $H(z)=1-(\log E)^{-q}\sqrt{-z}+O(|z|)$. Since $\sqrt{-(t+\mathsf{i}0)}=-\mathsf{i}\sqrt t$, we have $h(t)=(\log E)^{-q}\sqrt t+O(t)$ as $t\to 0$. The function $h$ is nonzero and locally integrable. Estimate~\eqref{eq:logarithmic-slit-decay} gives $|h(t)|\le C e^{-c w_p(t)}$ on $(0,\infty)$.
		
		Fix $\ell\in\Nzero$. For $0<\eta<1<T$ and $0<\eps<\pi/2$, let $\Gamma_{\eta,T,\eps}$ be the positively oriented boundary of $\{re^{\mathsf{i}\theta}:\eta<r<T,\ \eps<\theta<2\pi-\eps\}$. The contour runs from $\eta e^{\mathsf{i}\eps}$ to $Te^{\mathsf{i}\eps}$, counterclockwise along the outer circular arc, from $Te^{-\mathsf{i}\eps}$ to $\eta e^{-\mathsf{i}\eps}$, and clockwise along the inner circular arc. For $r>0$, denote the counterclockwise circular-arc integral by
		\[
		J_\eps(r)=\mathsf{i}r^{\ell+1}\int_\eps^{2\pi-\eps}e^{\mathsf{i}(\ell+1)\theta}H(re^{\mathsf{i}\theta})\,d\theta.
		\]
		Cauchy's integral theorem applied to $z^\ell H(z)$ on the region enclosed by $\Gamma_{\eta,T,\eps}$, followed by parametrization of the two radial segments, gives
		\begin{equation}\label{eq:keyhole-contour-identity}
			\begin{aligned}
				0=\int_{\Gamma_{\eta,T,\eps}}z^\ell H(z)\,dz=e^{\mathsf{i}(\ell+1)\eps}\int_\eta^T t^\ell H(te^{\mathsf{i}\eps})\,dt
				-e^{-\mathsf{i}(\ell+1)\eps}\int_\eta^T t^\ell H(te^{-\mathsf{i}\eps})\,dt
				+J_\eps(T)-J_\eps(\eta).
			\end{aligned}
		\end{equation}
		
		For fixed $\eta,T$, continuity up to each side of the positive real axis and dominated convergence give
		\[
		\begin{aligned}
			\lim_{\eps\to0}e^{\pm\mathsf{i}(\ell+1)\eps}\int_\eta^T t^\ell H(te^{\pm\mathsf{i}\eps})\,dt
			&=\int_\eta^T t^\ell H(t\pm\mathsf{i}0)\,dt,\\
			\lim_{\eps\to0}J_\eps(r)
			&=J(r):=\mathsf{i}r^{\ell+1}\int_0^{2\pi}e^{\mathsf{i}(\ell+1)\theta}H(re^{\mathsf{i}\theta})\,d\theta.
		\end{aligned}
		\]
		Taking $\eps\to0$ in \eqref{eq:keyhole-contour-identity} yields
		\begin{equation}\label{eq:keyhole-integral-1}
			\int_\eta^T t^\ell(H(t+\mathsf{i}0)-H(t-\mathsf{i}0))\,dt=J(\eta)-J(T).
		\end{equation}
		
		The boundedness of $H$ near zero and \eqref{eq:logarithmic-slit-decay} imply $|J(\eta)|\le C\eta^{\ell+1}$ and $|J(T)|\le CT^{\ell+1}e^{-c w_p(T)}$. Since $w_p(T)/\log T\to\infty$, these bounds give $J(\eta)\to0$ as $\eta\to0$ and $J(T)\to0$ as $T\to\infty$. The same bounds on $H$ also ensure that $t^\ell H(t\pm\mathsf{i}0)\in L^1(0,\infty)$. Using $H(t-\mathsf{i}0)=\overline{H(t+\mathsf{i}0)}$, $h(t)=\Ima H(t+\mathsf{i}0)$ and \eqref{eq:keyhole-integral-1}, we obtain
		\[
		\begin{aligned}
			2\mathsf{i}\int_0^\infty t^\ell h(t)\,dt=\int_0^\infty t^\ell(H(t+\mathsf{i}0)-H(t-\mathsf{i}0))\,dt=\lim_{\substack{\eta\to0\\T\to\infty}}(J(\eta)-J(T))=0.
		\end{aligned}
		\]
		All integer moments of $h$ vanish.
		
		Choose a nonnegative $\rho_0\in C_c^\infty((R+1,R+2))$ with integral one and put $F=\rho_0*h$. Then $F$ is smooth, real, and zero for $t\le R+1$. For large $t$ and $y\in\supp\rho_0$, we have $t-y\ge t/2$. Since $w_p$ is eventually increasing and $w_p(t/2)\ge2^{-1/2}w_p(t)$, the bound for $h$ gives $|F^{(k)}(t)|\le C\norm{\rho_0^{(k)}}_{L^1}e^{-(c/\sqrt2)w_p(t)}$. The rate $c/\sqrt2$ is independent of $k$; increasing the constants covers bounded $t\ge0$. Applying \eqref{eq:convolution-moments} with $(\rho,f_\gamma)$ replaced by $(\rho_0,h)$ proves that all integer moments of $F$ vanish. Identity~\eqref{eq:convolution-Laplace}, with the same replacements, and injectivity of the Laplace transform prove $F\ne0$, since $h\ne0$ and $\mathcal L\rho_0(\kappa)>0$ for $\kappa>0$.
	\end{proof}
	
	To estimate the causal fractional derivatives, including those in the Abel inversion, put $\omega_{\sigma,p}(r)=r^\sigma/\log^p(e+r)$ for $0<\sigma<1$ and $p\ge0$. In particular, $\omega_{1/2,p}=w_p$.
	
	\begin{lemma}\label{lem:logarithmic-fractional-tail}
		Let $0<\sigma<1$ and $p\ge0$, and let $F$ be a smooth function that vanishes on a neighborhood of $(-\infty,0]$ and satisfies \eqref{eq:profile-moments}. Suppose that $|F^{(k)}(t)|\le C_ke^{-a\omega_{\sigma,p}(t)}$ for $t\ge0$ and $k\in\Nzero$, with $a>0$ independent of $k$. For each fixed $\beta>0$, the function $G=D_+^\beta F$ is smooth, vanishes on every interval $(-\infty,R_0]$ on which $F$ vanishes, has all integer moments zero, and satisfies
		\begin{equation}\label{eq:logarithmic-fractional-tail}
			|G^{(k)}(t)|\le C'_ke^{-a'\omega_{\sigma,p}(t)},\quad t\ge0,\quad k\in\Nzero,
		\end{equation}
		with $a'>0$ independent of $k$.
	\end{lemma}
	
	\begin{proof}
		Write $\omega=\omega_{\sigma,p}$. We first establish the moment bound
		\begin{equation}\label{eq:general-logarithmic-moment-bound}
			\int_0^\infty t^M e^{-a\omega(t)}\,dt
			\le C (C_1M^{1/\sigma}\log^{p/\sigma}(e+M))^M,\quad M\in\N.
		\end{equation}
		For $y=\omega(t)$, we have $\log y=\sigma\log t-p\log\log(e+t)$, so $\log y/\log t\to\sigma>0$ as $t\to\infty$. In particular, $y\to\infty$ and $\log(e+y)$ and $\log(e+t)$ are comparable. Substituting this comparison into $t=y^{1/\sigma}\log^{p/\sigma}(e+t)$ gives $t\le C y^{1/\sigma}\log^{p/\sigma}(e+y)$ for large $t$. For $y=M\tau$, the inequality $\log(e+M\tau)\le(1+\log_+\tau)\log(e+M)$ yields
		\[
		y^{M/\sigma}\log^{pM/\sigma}(e+y)e^{-ay/2}\le M^{M/\sigma}(\tau^{1/\sigma}(1+\log_+\tau)^{p/\sigma}e^{-a\tau/2})^M\log^{pM/\sigma}(e+M).
		\]
		The expression in parentheses is bounded on $(0,\infty)$, independently of $M$. The bounded range of $t$ contributes at most a constant raised to $M$. These estimates bound $\sup_{t\ge0}t^Me^{-a\omega(t)/2}$ by the right-hand side of \eqref{eq:general-logarithmic-moment-bound}. Multiplying by the finite integral $\int_0^\infty e^{-a\omega(t)/2}\,dt$ proves \eqref{eq:general-logarithmic-moment-bound}.
		
		For integral $\beta$, the derivative bounds and the zero moments follow directly by differentiation and integration by parts. Otherwise, write $\beta=m-\delta$ with $m\in\N$ and $0<\delta<1$. The convolution identity \eqref{eq:causal-derivative-derivatives} gives $G^{(k)}=I_+^\delta F^{(m+k)}$ and proves smoothness and the support assertion. Fix $k$ and put $g=F^{(m+k)}$, so $|g(t)|\le C_ke^{-a\omega(t)}$. The integration-by-parts argument giving \eqref{eq:jth-moment} applies because $F$ vanishes near zero and all its derivatives decrease faster than every power; it gives $\int_0^\infty t^\ell g(t)\,dt=0$ for all $\ell\in\Nzero$.
		
		Split the defining integral for $I_+^\delta g(r)$ at $r/2$, as in \eqref{eq:tail-split}. The integral over $[r/2,r]$ is bounded by $C_kr^\delta e^{-a\omega(r/2)}$ for large $r$, since $\omega$ is eventually increasing. For the integral over $[0,r/2]$, apply \eqref{eq:tail-decomposition} to $g$; that identity uses only the zero moments of $g$ and the binomial expansion of $(1-x)^{\delta-1}$. The first inequality in \eqref{eq:fractional-tail-remainder-general} and \eqref{eq:general-logarithmic-moment-bound} give
		\[
		|\mathcal E_M(r)|\le C_kr^{\delta-1}\Big(\frac{C_1M^{1/\sigma}\log^{p/\sigma}(e+M)}r\Big)^M.
		\]
		Choose $0<\varepsilon_0\le1$ so that $C_1\varepsilon_0^{1/\sigma}\le1/2$, and take $M=\lfloor\varepsilon_0\omega(r)\rfloor$. For large $r$, we have $1\le M\le r$, and $\frac{M^{1/\sigma}\log^{p/\sigma}(e+M)}r\le\varepsilon_0^{1/\sigma}$. Also $M\ge\varepsilon_0\omega(r)/2$ for sufficiently large $r$. The remainder is bounded by $C_ke^{-c\omega(r)}$, with $c>0$ independent of $k$.
		
		For the finite sum, the inequality $\log(e+r\tau)\le(1+\log\tau)\log(e+r)$ for $\tau\ge1$ gives
		\[
		\omega(r)\log\tau\le C_{\sigma,p}\omega(r\tau),\quad
		C_{\sigma,p}=\sup_{\tau\ge1}\tau^{-\sigma}(1+\log\tau)^p\log\tau<\infty.
		\]
		The supremum is finite because the function being maximized is continuous and tends to zero at infinity. Decrease $\varepsilon_0$ so that $\varepsilon_0C_{\sigma,p}\le a/2$. Since $M\le\varepsilon_0\omega(r)$, the preceding inequality gives $\ell\log(t/r)\le a\omega(t)/2$ for $0\le\ell<M$ and $t\ge r$. For $r/2\le t\le r$, we have $(t/r)^\ell\le1$. In both ranges,
		\[
		(t/r)^\ell e^{-a\omega(t)}\le e^{-a\omega(t)/2},\quad 0\le\ell<M.
		\]
		
		For sufficiently large $r$, $\omega(t)\ge\omega(r/2)$ whenever $t\ge r/2$, so $e^{-a\omega(t)/2}\le e^{-a\omega(r/2)/4}e^{-a\omega(t)/4}$. The second factor has a finite integral over $[0,\infty)$ because $\omega(t)/\log t\to\infty$. Using $|g(t)|\le C_ke^{-a\omega(t)}$ and $0<d_\ell\le1$, we obtain the following bound for the finite sum:
		\[
		C_kMr^{\delta-1}\int_{r/2}^\infty e^{-a\omega(t)/2}\,dt
		\le C_kMr^{\delta-1}e^{-a\omega(r/2)/4},
		\]
		where the fixed integral $\int_0^\infty e^{-a\omega(t)/4}\,dt$ has been included in $C_k$.
		
		Since $\omega(r/2)\ge2^{-\sigma}\omega(r)$ and $M\le\varepsilon_0r^\sigma$, the limit $\omega(r)/\log r\to\infty$ allows the factor $Mr^{\delta-1}$ to be absorbed by decreasing the exponential rate. The same argument applies to the factor $r^\delta$ in the earlier estimate for the integral over $[r/2,r]$. Combining these estimates with the remainder bound proves \eqref{eq:logarithmic-fractional-tail} for sufficiently large $r$, with $a'>0$ depending only on $a,\sigma,p,\beta$, independently of $k$. Increasing each $C'_k$ on a bounded interval gives the estimate for all $r\ge0$.
		
		The decay estimates give $t^\ell F,t^\ell G\in L^1(0,\infty)$ for every $\ell\in\Nzero$. The zero moments of $F$ imply \eqref{eq:Laplace-flatness}. Substitution into \eqref{eq:Laplace-fractional-product}, which follows from \eqref{eq:Laplace-derivative}, gives $\lim_{\kappa\to0^+}(\mathcal LG)^{(m)}(\kappa)=0$ for every $m$. Equation~\eqref{eq:Laplace-moment-limit} proves that all integer moments of $G$ vanish.
	\end{proof}
	
	\begin{proof}[Proof of Proposition~\ref{prop:logarithmic-sharpness}]
		Choose $F$ from Lemma~\ref{lem:logarithmic-zero-moments} and set $\beta=2(t-s)$ and $G=D_+^\beta F$. Lemma~\ref{lem:logarithmic-fractional-tail}, with $\sigma=1/2$, gives $|F^{(k)}(r)|+|G^{(k)}(r)|\le C_ke^{-aw_p(r)}$ for $r\ge0$ and $k\in\Nzero$, where $a>0$ is independent of $k$. Both profiles belong to $\cS(\R)$, vanish for $r\le R+1$, and have all integer moments zero. Both are nonzero: if $G=0$, the identity $\mathcal LG(\kappa)=\kappa^\beta\mathcal LF(\kappa)$ and injectivity of the Laplace transform would imply $F=0$.
		
		In dimension three, define $u=-\sin(\pi t)G(|x|)/|x|$ and $v=\sin(\pi s)F(|x|)/|x|$, with both functions zero near the origin. Proposition~\ref{prop:radial-cancellation} gives the free resolvent cancellation \eqref{eq:radial-cancellation}. Their spatial derivatives satisfy \eqref{eq:logarithmic-counterexample-decay}, since differentiation introduces derivatives of $F,G$ and nonpositive powers of $|x|$. For higher odd dimensions, apply $\mathcal T$ from Lemma~\ref{lem:dimension} repeatedly. Its explicit formula shows that it preserves these logarithmic decay bounds, as well as nonvanishing, vanishing on $B_R$, and the resolvent cancellation.
		
		For even dimensions, use the Abel construction \eqref{eq:even-Abel-inverse}. For $J\in\{F,G\}$, let $H_J(x)=J(\sqrt x)$ for $x>0$ and zero otherwise. This function is smooth and vanishes on $(-\infty,(R+1)^2]$. Since $H_J$ vanishes on $[0,1]$, the bounds on $J$, \eqref{eq:Abel-chain-derivatives} with $F$ replaced by $J$, and $\log(e+\sqrt x)\le\log(e+x)$ yield
		\[
		|H_J^{(k)}(x)|\le C_ke^{-a\omega_{1/4,p}(x)},\quad x\ge0,\quad k\in\Nzero.
		\]
		The change of variables in \eqref{eq:Abel-input-moments}, with $F$ replaced by $J$, shows that all integer moments of $H_J$ vanish. Apply Lemma~\ref{lem:logarithmic-fractional-tail} with $\sigma=1/4$ and derivative order $1/2$ to $q_J=(2/\sqrt\pi)D_{+,x}^{1/2}H_J$. After taking the smaller of the two rates for $J=F$ and $J=G$, we have $|q_J^{(k)}(x)|\le C_ke^{-b\omega_{1/4,p}(x)}$ for $x\ge0$, with $b>0$ independent of $k$ and $J$. Since
		\[
		\omega_{1/4,p}(r^2)=\frac{\sqrt r}{\log^p(e+r^2)}\ge2^{-p}w_p(r),
		\]
		the spatial derivatives of $f_J(|x|)=q_J(|x|^2)$ are bounded by polynomials in $x$ times $e^{-2^{-p}bw_p(|x|)}$. Absorbing these polynomial factors gives \eqref{eq:logarithmic-counterexample-decay} with $c=2^{-p-1}b$, independently of the derivative order and $J$. Both functions are Schwartz and zero on $B_{R+1}$.
		
		Identity~\eqref{eq:Abel-causal-inversion} applies to $H_J$, since it is Schwartz and vanishes near zero. It gives \eqref{eq:Abel-profile-identity} with $F$ replaced by $J$, and also proves $f_J\ne0$. Applying \eqref{eq:even-radial-Green-formula} and \eqref{eq:even-radial-Laplace-calculation} to $f_J$ gives
		\[
		(-\Delta+\kappa^2)^{-1}f_J(|x|)=I_0(\kappa|x|)\mathcal LJ(\kappa)\quad\text{on }B_R.
		\]
		In dimension two, set $u=-\sin(\pi t)f_G(|x|)$ and $v=\sin(\pi s)f_F(|x|)$. Since $G=D_+^{2(t-s)}F$, \eqref{eq:Laplace-derivative} gives $\mathcal LG(\kappa)=\kappa^{2(t-s)}\mathcal LF(\kappa)$ for $\kappa>0$. Substituting the definitions of $u,v$ into the preceding resolvent formula and using $\kappa^{2s}\mathcal LG(\kappa)=\kappa^{2t}\mathcal LF(\kappa)$ proves \eqref{eq:even-resolvent-cancellation}. Both inputs are nonzero because $s,t$ are nonintegers. Repeated application of $\mathcal T$ gives the construction in every higher even dimension.
		
		In each dimension, the free resolvent of $f_\kappa=\sin(\pi s)\kappa^{2s}u+\sin(\pi t)\kappa^{2t}v$ vanishes on $B_R$. Lemma~\ref{lem:transfer} gives the same vanishing for $(L+\kappa^2)^{-1}f_\kappa$. The converse implication in Proposition~\ref{prop:resolvent}, with $\lambda=\kappa^2$, proves $(L^su+L^tv)|_{B_R}=0$.
		
		Finally, \eqref{eq:logarithmic-counterexample-decay} and $w_p(r)/\log r\to\infty$ show that $e^{aw_p(|x|)}u,e^{aw_p(|x|)}v\in L^2(\R^n)$ for every $0<a<c$. Both inputs belong to $\mathscr D_{\log,p}$.
	\end{proof}

	\section{Failure of entanglement under small smooth perturbations}\label{sec:general-counterexample}
	
	We prove Theorem~\ref{thm:general-counterexample} by constructing coefficients that converge to the identity in every fixed unweighted $C^k$ norm and admit compactly supported counterexamples in dimension two. The orders $1/4$ and $3/4$ satisfy the even-dimensional restriction in Theorem~\ref{thm:arbitrary}. For each coefficient field, we construct a nonzero outgoing wave that vanishes on $B_1$ for every nonnegative time. Its Laplace transform gives the required resolvent cancellation. We begin with the implication from wave propagation to fractional cancellation.
	
	\begin{lemma}\label{lem:wave-cancellation}
		Let $\cO\subset\R^n$ be a nonempty open set, let $u,v\in C_c^\infty(\R^n)$, and let $w$ solve
		\[
		\begin{cases}
			(\partial_\tau^2 +L)w=0&\text{in }(0,\infty)\times\R^n,\\
			w(0,\cdot)=v&\text{on }\R^n,\\
			w_\tau(0,\cdot)=u&\text{on }\R^n.
		\end{cases}
		\]
		Assume that $w\in C([0,T];H^1(\R^n))\cap C^1([0,T];L^2(\R^n))$ for every $T>0$. If $w(\tau,\cdot)|_{\cO}=0$ for every $\tau\ge0$, then
		\begin{equation}\label{eq:wave-cancellation}
			u|_{\cO}=v|_{\cO}=0,\quad (L^{1/4}u+L^{3/4}v)\big|_{\cO}=0.
		\end{equation}
	\end{lemma}
	
	\begin{proof}
		Since $w(\tau,\cdot)|_{\cO}=0$ for every $\tau\ge0$, evaluating at $\tau=0$ and differentiating there gives $v|_{\cO}=u|_{\cO}=0$. The spectral representation of the wave solution for the nonnegative self-adjoint operator $L$ is
		\[
		w(\tau)=\cos(\tau\sqrt L)v+\frac{\sin(\tau\sqrt L)}{\sqrt L}u;
		\]
		see \cite[Section~5.1, Problem~5.3, formula~(5.10)]{Teschl14}. The second operator is defined by the multiplier $\mu\mapsto\sin(\tau\sqrt\mu)/\sqrt\mu$ for $\mu>0$, with value $\tau$ at $\mu=0$.
		
		The bounds $|\cos(\tau\sqrt\mu)|\le1$ and $|\sin(\tau\sqrt\mu)/\sqrt\mu|\le\tau$ give $\norm{w(\tau)}_2\le\norm{v}_2+\tau\norm{u}_2$. In particular,
		\[
		\int_0^\infty e^{-\kappa\tau}\norm{w(\tau)}_2\,d\tau
		\le\frac{\norm{v}_2}{\kappa}+\frac{\norm{u}_2}{\kappa^2}<\infty,\quad\kappa>0,
		\]
		so the Laplace transform exists as an $L^2(\R^n)$-valued integral. For $\mu\ge0$ and $\kappa>0$, the scalar identities
		\[
		\begin{aligned}
			\int_0^\infty e^{-\kappa\tau}\cos(\tau\sqrt\mu)\,d\tau=\frac{\kappa}{\kappa^2+\mu},\quad 
			\int_0^\infty e^{-\kappa\tau}\frac{\sin(\tau\sqrt\mu)}{\sqrt\mu}\,d\tau=\frac1{\kappa^2+\mu}
		\end{aligned}
		\]
		hold with the same convention at $\mu=0$. Applying these identities through the spectral calculus gives
		\begin{equation}\label{eq:wave-Laplace}
			\int_0^\infty e^{-\kappa\tau}w(\tau)\,d\tau
			=(L+\kappa^2)^{-1}(u+\kappa v),\quad\kappa>0.
		\end{equation}
		The interchange with the spectral integral is justified by the bounds $1$ and $\tau$ for the two scalar multipliers.
		
		Restriction from $L^2(\R^n)$ to $L^2(\cO)$ is bounded and commutes with the integral in \eqref{eq:wave-Laplace}. Since $w(\tau,\cdot)|_{\cO}=0$ for every $\tau\ge0$, the left-hand side of \eqref{eq:wave-Laplace} vanishes on $\cO$. Multiplying \eqref{eq:wave-Laplace} by $\sin(\pi/4)\kappa^{1/2}$ and using $\sin(\pi/4)=\sin(3\pi/4)$ gives
		\[
		\sin(\pi/4)\kappa^{1/2}(L+\kappa^2)^{-1}u
		+\sin(3\pi/4)\kappa^{3/2}(L+\kappa^2)^{-1}v
		=0\quad\text{on }\cO,\quad\kappa>0.
		\]
		With $\lambda=\kappa^2$, this is \eqref{eq:intro-resolvent} for $s_1=1/4$, $s_2=3/4$, and $b_1=b_2=1$. Together with $u|_{\cO}=v|_{\cO}=0$, the converse part of Proposition~\ref{prop:resolvent} proves \eqref{eq:wave-cancellation}.
	\end{proof}
	
	\begin{proof}[Proof of Theorem~\ref{thm:general-counterexample}]
		For $x=r(\cos\theta,\sin\theta)\ne0$, write $e_r=(\cos\theta,\sin\theta)$ and $e_\theta=(-\sin\theta,\cos\theta)$. Choose a fixed function $\chi\in C^\infty([0,\infty))$ with $0\le\chi\le1$, $\chi=0$ on $[0,1/2]$, and $\chi=1$ on $[1,\infty)$. For $\ell\in\N$, define
		\begin{equation}\label{eq:near-identity-coefficient}
			\begin{gathered}
				a_\ell^\infty=\frac{(\ell+1/2)^2}{\ell^2},\quad
				a_\ell(r)=1+(a_\ell^\infty-1)\chi(r),\\
				A_\ell(x)=e_r\otimes e_r+a_\ell(r)e_\theta\otimes e_\theta,\quad x\ne0,
			\end{gathered}
		\end{equation}
		and set $A_\ell(0)=I$. Since $e_r\otimes e_r+e_\theta\otimes e_\theta=I$, the matrix $A_\ell$ equals $I$ on $B_{1/2}$ and is smooth and real symmetric on $\R^2$. Its eigenvalues are $1$ and $a_\ell(r)$, both contained in $[1,9/4]$, independently of $\ell$.
		
		We have $A_\ell-I=(\frac1\ell+\frac1{4\ell^2})\chi(r)e_\theta\otimes e_\theta$. The matrix field $\chi(r)e_\theta\otimes e_\theta$, extended by zero at the origin, has bounded Cartesian derivatives of every order. Indeed, it vanishes near zero, and derivatives of order $k$ of the angular tensor are $O(r^{-k})$ away from zero. It follows that
		\[
		\norm{A_\ell-I}_{C^k(\R^2)}\le\frac{C_k}{\ell},\quad k\in\Nzero,
		\]
		where $C_k$ depends only on $k$ and $\chi$. Since $a_\ell^\infty>1$, we also have $A_\ell(x)\ne I$ for every $|x|\ge1$, so $A_\ell-I$ is not compactly supported.
		
		Set $L_\ell=-\Div(A_\ell\nabla)$ on $\R^2$. The polar coordinate formula for divergence gives
		\[
		L_\ell=-\partial_r^2-\frac1r\partial_r-\frac{a_\ell(r)}{r^2}\partial_\theta^2,\quad r>0.
		\]
		For $r>1$, we have $a_\ell(r)=a_\ell^\infty$. In this region, consider $w(\tau,x)=r^{-1/2}z(\tau,r)\cos(\ell\theta)$. The identities $(-\partial_r^2-\frac1r\partial_r)(r^{-1/2}z)=r^{-1/2}(-z_{rr}-\frac{z}{4r^2})$ and $\partial_\theta^2\cos(\ell\theta)=-\ell^2\cos(\ell\theta)$ show how the factor $r^{-1/2}$ removes the first radial derivative. The choice of $a_\ell^\infty$ gives $a_\ell^\infty\ell^2-1/4=\ell(\ell+1)$, and we obtain
		\begin{equation}\label{eq:angular-reduction}
			w_{\tau\tau}+L_\ell w=r^{-1/2}\Big(z_{\tau\tau}-z_{rr}+\frac{\ell(\ell+1)}{r^2}z\Big)\cos(\ell\theta)
			\quad\text{for }r>1.
		\end{equation}
		
		To construct a solution of the resulting radial wave equation, consider the differential expressions
		\[
		H_m=-\partial_r^2+\frac{m(m+1)}{r^2},\quad m\in\Nzero,
		\]
		on $(0,\infty)$. For $m\ge1$, put $D_m=-\partial_r+m/r$. Its formal adjoint with respect to the measure $dr$ is $D_m^*=\partial_r+m/r$. Multiplication of these differential expressions gives $D_mD_m^*=H_m$ and $D_m^*D_m=H_{m-1}$ on $(0,\infty)$. Composing with $D_m$ gives
		\begin{equation}\label{eq:intertwining}
			H_mD_m=D_mD_m^*D_m=D_mH_{m-1}.
		\end{equation}
		Since $D_m$ acts only on $r$, it commutes with $\partial_\tau^2$. Identity~\eqref{eq:intertwining} gives $(\partial_\tau^2+H_m)D_m=D_m(\partial_\tau^2+H_{m-1})$.
		
		Choose a nonzero real $F\in C_c^\infty((2,3))$, extended by zero to $\R$. For $\tau\ge0$ and $r>0$, define $z_\ell(\tau,r)=D_\ell\cdots D_1[F(\cdot-\tau)](r)$, where each $D_m=-\partial_r+m/r$ acts in the variable $r$. The coefficients of these operators are smooth on $(0,\infty)$, so $z_\ell$ is smooth up to $\tau=0$ for $r>0$. Since $(\partial_\tau^2-\partial_r^2)F(r-\tau)=0$, repeated application of \eqref{eq:intertwining} gives
		\begin{equation}\label{eq:z_ell}
			(\partial_\tau^2+H_\ell)z_\ell=0
			\quad\text{in }(0,\infty)\times(0,\infty).
		\end{equation}
		Each $D_m$ is local on $(0,\infty)$ and preserves support there. For every $\tau\ge0$, we have $\supp z_\ell(\tau,\cdot)\subset\tau+\supp F\Subset(2+\tau,3+\tau)$. In particular, $z_\ell(\tau,r)=0$ whenever $0<r\le2+\tau$.
		
		For $\tau\ge0$, define $w_\ell(\tau,\cdot)$ on $\R^2$ by
		\[
		w_\ell(\tau,x)=
		\begin{cases}
			r^{-1/2}z_\ell(\tau,r)\cos(\ell\theta), & x\ne0,\\
			0, & x=0,
		\end{cases}
		\]
		where $x=r(\cos\theta,\sin\theta)$ when $x\ne0$. The expression away from the origin is smooth, and the support inclusion shows that it vanishes on the fixed neighborhood $B_2$ of the origin for every $\tau\ge0$. This defines a smooth function on $[0,\infty)\times\R^2$, with $w_\ell(\tau,x)=0$ whenever $|x|\le2+\tau$.
		
		For $r>1$, equations~\eqref{eq:angular-reduction} and \eqref{eq:z_ell} give $(\partial_\tau^2+L_\ell)w_\ell=0$. The same equation holds on $(0,\infty)\times B_2$, where $w_\ell$ vanishes identically. These regions cover $(0,\infty)\times\R^2$.
		
		Since each $D_m$ is independent of $\tau$, the initial displacement and velocity are given, for $x\ne0$, by
		\[
		\begin{aligned}
			v_\ell(x)&:=w_\ell(0,x)
			=r^{-1/2}(D_\ell\cdots D_1F)(r)\cos(\ell\theta),\\
			u_\ell(x)&:=\partial_\tau w_\ell(0,x)
			=-r^{-1/2}(D_\ell\cdots D_1F')(r)\cos(\ell\theta),
		\end{aligned}
		\]
		with both functions defined to be zero at the origin. Since $\supp F'\subset\supp F\Subset(2,3)$ and each $D_m$ preserves support, both $u_\ell$ and $v_\ell$ belong to $C_c^\infty(\{x\in\R^2:2<|x|<3\})$. The function $w_\ell$ solves the Cauchy problem
		\[
		\begin{cases}
			(\partial_\tau^2+L_\ell)w_\ell=0
			&\text{in }(0,\infty)\times\R^2,\\
			w_\ell(0,\cdot)=v_\ell
			&\text{on }\R^2,\\
			\partial_\tau w_\ell(0,\cdot)=u_\ell
			&\text{on }\R^2.
		\end{cases}
		\]
		
		To prove that both initial data are nonzero, let $f\in C_c^\infty((0,\infty))$ satisfy $D_mf=0$. The equation $f'=mf/r$ gives $f(r)=Cr^m$ on $(0,\infty)$, and compact support forces $C=0$. Each $D_m$ is injective on $C_c^\infty((0,\infty))$ and maps this space into itself, so their composition is also injective. The chosen function $F$ is nonzero, and $F'$ is nonzero because a constant compactly supported function must vanish. Both $D_\ell\cdots D_1F$ and $D_\ell\cdots D_1F'$ are nonzero. Since $\cos(\ell\theta)$ is not identically zero, neither $v_\ell$ nor $u_\ell$ vanishes identically.
		
		For every $T>0$, the functions $w_\ell(\tau,\cdot)$ and their derivatives have support in $\overline{B_{T+3}}$ for $0\le\tau\le T$. Smoothness gives $w_\ell\in C([0,T];H^1(\R^2))\cap C^1([0,T];L^2(\R^2))$. Uniqueness in this energy class identifies $w_\ell$ with the wave evolution generated by $L_\ell$, with initial displacement $v_\ell$ and initial velocity $u_\ell$. Since $w_\ell(\tau,\cdot)|_{B_1}=0$ for every $\tau\ge0$, Lemma~\ref{lem:wave-cancellation}, applied with $L=L_\ell$ and $\cO=B_1$, proves \eqref{eq:general-counterexample}.
	\end{proof}
	
	In these examples, the radial factor $D_\ell\cdots D_1F$ is nonzero at some $r_\ell\in(2,3)$. The entire circle $\{|x|=r_\ell\}$ lies in $\supp v_\ell$, including the angular zeros of $\cos(\ell\theta)$ by closure. This circle separates $B_1$ from infinity, so $B_1$ lies in a bounded component of $\R^2\setminus(\supp u_\ell\cup\supp v_\ell)$. By \eqref{eq:near-identity-coefficient}, we have $A_\ell-I=(a_\ell^\infty-1)e_\theta\otimes e_\theta$ for $|x|\ge1$, while $\norm{A_\ell-I}_{C^k(\R^2)}\to0$ for every fixed $k\in\Nzero$.

	\section{Applications to inverse problems}\label{sec:inverse-problems}
	
	Throughout this section, the operator $L$ is known and satisfies \eqref{eq:ellipticity}. Let $\Omega$ and $\Omega_e=\R^n\setminus\overline\Omega$ be the bounded Lipschitz domain and its exterior introduced in Section~\ref{subsec:inverse-introduction}. Recovery of a potential and of a local anisotropic operator uses Theorem~\ref{thm:single}, with arbitrary nonempty open observation sets in $\Omega_e$. To recover unknown fractional orders, we first group orders with the same fractional part and then apply Theorem~\ref{thm:exterior}.
	
	\subsection{The exterior Dirichlet problem}
	
	Let $P$ be given by \eqref{eq:inverse-P}. Set $t=s_N$ if $c=0$ and $t=\max\{m,s_N\}$ otherwise, and write $X=\HL^t$ and $X_0(\Omega)=\overline{C_c^\infty(\Omega)}^X$. For $q\in L^\infty(\Omega)$, consider the bounded Hermitian form
	\begin{equation}\label{eq:inverse-energy-form}
		\mathcal B_{P,q}(u,v)=c(L^{m/2}u,L^{m/2}v)_{L^2}+\sum_{k=1}^N b_k(L^{s_k/2}u,L^{s_k/2}v)_{L^2}+\int_\Omega qu\overline v\,dx,\quad u,v\in X,
	\end{equation}
	where the first term is omitted when $c=0$. We assume that the homogeneous problem has only the zero solution:
	\begin{equation}\label{eq:inverse-nonresonance}
		v\in X_0(\Omega),\quad \mathcal B_{P,q}(v,\phi)=0\text{ for every }\phi\in X_0(\Omega)\implies v=0.
	\end{equation}
	For $f\in C_c^\infty(\Omega_e)$, a weak solution of \eqref{eq:inverse-problem} is a function $u_f\in X$ satisfying
	\begin{equation}\label{eq:inverse-weak-solution}
		u_f-f\in X_0(\Omega),\quad \mathcal B_{P,q}(u_f,\phi)=0\quad\text{for every }\phi\in X_0(\Omega).
	\end{equation}
	Every element of $X_0(\Omega)$ vanishes almost everywhere outside $\Omega$, by its defining approximation in $L^2(\R^n)$. In particular, $u_f=f$ on $\Omega_e$ and $\supp u_f\subset\overline\Omega\cup\supp f$. Testing against $C_c^\infty(\Omega)$ gives $(P+q)u_f=0$ in $\Omega$ in the sense of distributions.
	
	We prove unique solvability under \eqref{eq:inverse-nonresonance}. The spectral theorem gives $X=\Dom(L^{t/2})$ with $\norm{v}_X^2\asymp\norm{v}_{L^2}^2+\norm{L^{t/2}v}_{L^2}^2$. The identity $\HL^1=H^1(\R^n)$, complex interpolation for $0<t<1$, and the inclusion $\HL^t\hookrightarrow\HL^1$ for $t\ge1$ imply $X\hookrightarrow H^{\min\{t,1\}}(\R^n)$. Since $t>0$ and the elements of $X_0(\Omega)$ are supported in the fixed bounded set $\overline\Omega$, the Rellich compactness theorem gives the compact embedding $X_0(\Omega)\hookrightarrow L^2(\Omega)$. We also have
	\begin{equation}\label{eq:inverse-Poincare}
		\norm{v}_{L^2}\le C\norm{L^{t/2}v}_{L^2},\quad v\in X_0(\Omega).
	\end{equation}
	Otherwise, a sequence $v_k\in X_0(\Omega)$ with $\norm{v_k}_{L^2}=1$ and $L^{t/2}v_k\to0$ in $L^2(\R^n)$ would, by compactness, have a subsequence converging in $L^2(\R^n)$ to $v\in\ker L^{t/2}=\ker L=\{0\}$ with $\norm{v}_{L^2}=1$, a contradiction.
	
	Let $d\ne0$ be the coefficient of $L^t$ in $P$, and put $\sigma=d/|d|$. The highest power has a single nonzero coefficient, since the fractional orders are distinct nonintegers and the local order $m$ is an integer. For $0<\rho<t$ and $\eps>0$, splitting at $\lambda=\eps^{-1/(t-\rho)}$ gives $\lambda^\rho\le\eps\lambda^t+C_{\eps,\rho,t}$ for all $\lambda\ge0$; the case $\rho=0$ is immediate. Integration against the spectral measure of $v\in X$ gives
	\begin{equation}\label{eq:inverse-lower-power-bound}
		\norm{L^{\rho/2}v}_{L^2}^2\le\eps\norm{L^{t/2}v}_{L^2}^2+C_{\eps,\rho,t}\norm{v}_{L^2}^2.
	\end{equation}
	Let $S$ be the sum of the absolute values of the coefficients of the powers below $t$ in $P$. Apply this estimate to each lower power with $\eps=|d|/(2(1+S))$. Since $\sigma d=|d|$, $\eps S\le|d|/2$, and $q\in L^\infty(\Omega)$, we obtain
	\[
	\sigma\mathcal B_{P,q}(v,v)\ge\frac{|d|}{2}\norm{L^{t/2}v}_{L^2}^2-C\norm{v}_{L^2}^2,\quad v\in X_0(\Omega),
	\]
	where $C\ge0$ is independent of $v$. For $\mu>C$, the form $\widetilde{\mathcal B}(u,v)=\sigma\mathcal B_{P,q}(u,v)+\mu(u,v)_{L^2}$ is bounded and coercive on $X_0(\Omega)$ by the norm equivalence above.
	
	By the Lax--Milgram theorem, there is a bounded linear operator $K:X_0(\Omega)\to X_0(\Omega)$ satisfying
	\begin{equation}\label{eq:inverse-Fredholm-operator}
		\widetilde{\mathcal B}(Kv,\phi)=\mu(v,\phi)_{L^2},\quad v,\phi\in X_0(\Omega).
	\end{equation}
	Taking $\phi=Kv$ and using coercivity gives $\norm{Kv}_X\le C_1\norm{v}_{L^2(\Omega)}$. Applying this estimate to differences of a subsequence converging in $L^2(\Omega)$ proves that $K$ is compact. The identity $\sigma\mathcal B_{P,q}(v,\phi)=\widetilde{\mathcal B}((I-K)v,\phi)$ and \eqref{eq:inverse-nonresonance} imply that $I-K$ is injective. The Fredholm alternative gives a bounded inverse.
	
	For any bounded conjugate linear functional $F$ on $X_0(\Omega)$, let $z\in X_0(\Omega)$ be the unique solution of $\widetilde{\mathcal B}(z,\phi)=\sigma F(\phi)$ for every $\phi\in X_0(\Omega)$. Then $w=(I-K)^{-1}z$ is the unique element of $X_0(\Omega)$ satisfying $\mathcal B_{P,q}(w,\phi)=F(\phi)$. Taking $F(\phi)=-\mathcal B_{P,q}(f,\phi)$ gives the unique solution $u_f=f+w$ of \eqref{eq:inverse-weak-solution}, with $\norm{u_f}_X\le C\norm{f}_X$. Taking $F(\phi)=(h,\phi)_{L^2(\Omega)}$ gives, for every $h\in L^2(\Omega)$, a unique solution of
	\begin{equation}\label{eq:inverse-source-solution}
		w\in X_0(\Omega),\quad \mathcal B_{P,q}(w,\phi)=(h,\phi)_{L^2(\Omega)}\quad\text{for every }\phi\in X_0(\Omega).
	\end{equation}
	If $c\ge0$, all $b_k>0$, and $q\ge0$, the original form is coercive by \eqref{eq:inverse-Poincare}, so \eqref{eq:inverse-nonresonance} holds.
	
	Finally, for $g\in C_c^\infty(\Omega_e)$, the potential term in $\mathcal B_{P,q}(u_f,g)$ vanishes. For each power $\rho$ occurring in $P$, the spectral calculus gives $(L^{\rho/2}u_f,L^{\rho/2}g)_{L^2}=(u_f,L^\rho g)_{L^2}$. The exterior DN map in \eqref{eq:inverse-map} satisfies
	\begin{equation}\label{eq:inverse-weak-DN}
		\ip{\Lambda_{P,q}f}{g}=\mathcal B_{P,q}(u_f,g)=\ip{Pu_f}{g},\quad f,g\in C_c^\infty(\Omega_e).
	\end{equation}
	This formula allows DN maps with different fractional orders to be compared on the same exterior test functions, even when their energy spaces differ.
	
	\subsection{Runge approximation and recovery of a potential}\label{subsec:inverse-potential}
	
	\begin{lemma}\label{lem:inverse-Runge}
		Under \eqref{eq:inverse-nonresonance}, the set $\{u_f|_\Omega:f\in C_c^\infty(W)\}$ is dense in $L^2(\Omega)$ for every nonempty open $W\subset\Omega_e$.
	\end{lemma}
	
	\begin{proof}
		Let $h\in L^2(\Omega)$ be orthogonal to all these restrictions, and let $w$ be the solution of \eqref{eq:inverse-source-solution}. Fix $f\in C_c^\infty(W)$. Since $w\in X_0(\Omega)$, the equation for $u_f$ and Hermitian symmetry give $\mathcal B_{P,q}(w,u_f)=0$. On the other hand, $u_f-f\in X_0(\Omega)$, so the equation for $w$ gives
		\[
		\mathcal B_{P,q}(w,u_f-f)=(h,u_f-f)_{L^2(\Omega)}=(h,u_f)_{L^2(\Omega)}=0,
		\]
		where we used $f=0$ in $\Omega$ and the orthogonality of $h$. Subtraction gives
		\begin{equation}\label{eq:inverse-Runge-exterior-identity}
			\mathcal B_{P,q}(w,f)=0,\quad f\in C_c^\infty(W).
		\end{equation}
		
		For $f\in C_c^\infty(W)$, the identity $\mathcal B_{P,q}(w,f)=0$ is $\ip{Pw}{f}=0$, since the potential term vanishes. We obtain $Pw=0$ on $W$. Also $w=0$ there, since $w\in X_0(\Omega)$. Locality gives $L^mw=0$ on $W$, so Theorem~\ref{thm:single}, applied to the fractional terms of $P$, gives $w=0$ on $\R^n$. The equation for $w$ now gives $(h,\phi)_{L^2(\Omega)}=0$ for every $\phi\in C_c^\infty(\Omega)$, proving $h=0$ and the asserted density.
	\end{proof}
	
	\begin{proof}[Proof of Theorem~\ref{thm:inverse-potential}]
		Fix $f\in C_c^\infty(W_1)$ and let $u_{j,f}$ be the solution associated with $q_j$. Their difference $w=u_{1,f}-u_{2,f}$ is zero on $\Omega_e$, since both solutions have exterior value $f$. Equations~\eqref{eq:inverse-potential-data} and \eqref{eq:inverse-weak-DN} give $Pw=0$ on $W_2$. The local term $cL^mw$ vanishes there by locality, so Theorem~\ref{thm:single} applied to the remaining fractional sum gives $w=0$ on $\R^n$.
		
		Since $u_{1,f}=u_{2,f}$, subtracting their interior equations gives
		\begin{equation}\label{eq:inverse-potential-product}
			(q_1-q_2)u_{1,f}=0\quad\text{almost everywhere in }\Omega.
		\end{equation}
		By Lemma~\ref{lem:inverse-Runge}, there is a sequence $f_\ell\in C_c^\infty(W_1)$ such that $u_{1,f_\ell}|_\Omega\to1$ in $L^2(\Omega)$. Since $q_1-q_2\in L^\infty(\Omega)$, multiplication by this difference is continuous on $L^2(\Omega)$. Taking the limit in \eqref{eq:inverse-potential-product} yields $q_1-q_2=0$ almost everywhere in $\Omega$.
	\end{proof}
	
	If $q_j\in C(\Omega)\cap L^\infty(\Omega)$ for $j=1,2$, one nonzero exterior datum $f\in C_c^\infty(\Omega_e)$ suffices. Suppose its measurements agree on a nonempty open set $W\subset\Omega_e$. Equation~\eqref{eq:inverse-weak-DN}, locality of $L^m$, and Theorem~\ref{thm:single} give $u_{1,f}=u_{2,f}=:u$, with the product identity \eqref{eq:inverse-potential-product} for this fixed datum. If $q_1-q_2$ were nonzero at some point, continuity would give a ball $B\Subset\Omega$ and $\eta>0$ such that $|q_1-q_2|\ge\eta$ on $B$. The product identity gives $u=0$ on $B$, and the interior equation gives $Pu=0$ there. The local term $cL^mu$ vanishes on $B$, so Theorem~\ref{thm:single} applies to the fractional terms and gives $u=0$ on $\R^n$. This contradicts $u=f\ne0$ on $\Omega_e$, proving
	\begin{equation}\label{eq:inverse-single-datum-potential}
		q_1=q_2\quad\text{in }\Omega.
	\end{equation}
	
	\subsection{Recovery of a local anisotropic operator}
	
	Let $\mathcal N=\sum_{k=1}^N b_kL^{s_k}$, where $0<s_1<\cdots<s_N$ are nonintegers and $b_k\in\R\setminus\{0\}$, as in Section~\ref{subsec:inverse-introduction}. For \eqref{eq:inverse-anisotropic-problem}, put $t=\max\{1,s_N\}$, $Y=\HL^t$, and $Y_0(\Omega)=\overline{C_c^\infty(\Omega)}^Y$. Define
	\begin{equation}\label{eq:inverse-anisotropic-form}
		\mathcal B_{C,q}^{\mathcal N}(u,v)=\int_\Omega C\nabla u\cdot\nabla\overline v\,dx+\sum_{k=1}^N b_k(L^{s_k/2}u,L^{s_k/2}v)_{L^2}+\int_\Omega qu\overline v\,dx.
	\end{equation}
	Since $t\ge1$ and $\HL^1=H^1(\R^n)$, this form is bounded on $Y$. For $v\in Y_0(\Omega)$, the gradient vanishes almost everywhere outside $\overline\Omega$. If $s_N<1$, uniform positive definiteness of $C$ and \eqref{eq:ellipticity} give $\int_\Omega C\nabla v\cdot\nabla\overline v\,dx\ge c_0\norm{L^{1/2}v}_{L^2}^2$ for some $c_0>0$, and all fractional terms have lower order. If $s_N>1$, then $b_N\norm{L^{s_N/2}v}_{L^2}^2$ is the highest-order term, while the local term is bounded by a constant times $\norm{L^{1/2}v}_{L^2}^2$. Set $\sigma=1$ in the first case and $\sigma=b_N/|b_N|$ in the second. Applying \eqref{eq:inverse-lower-power-bound} to the terms below order $t$ gives constants $a>0$ and $C_0\ge0$ such that
	\[
	\sigma\mathcal B_{C,q}^{\mathcal N}(v,v)\ge a\norm{L^{t/2}v}_{L^2}^2-C_0\norm{v}_{L^2}^2,\quad v\in Y_0(\Omega).
	\]
	The embedding $Y_0(\Omega)\subset L^2(\Omega)$ is compact. Add a sufficiently large positive multiple of the $L^2$ inner product to $\sigma\mathcal B_{C,q}^{\mathcal N}$, and apply the construction \eqref{eq:inverse-Fredholm-operator} on $Y_0(\Omega)$ with this coercive form. It gives a compact operator $K$ for which the assumed uniqueness of the homogeneous problem implies injectivity of $I-K$. The Fredholm alternative gives existence and uniqueness for the weak problem and for the corresponding source equation.
	
	The weak solution is characterized by $u_f-f\in Y_0(\Omega)$ and $\mathcal B_{C,q}^{\mathcal N}(u_f,v)=0$ for every $v\in Y_0(\Omega)$. If $g\in C_c^\infty(\Omega_e)$, both integrals over $\Omega$ vanish in the pairing with $g$, and
	\begin{equation}\label{eq:inverse-anisotropic-weak-DN}
		\ip{\Lambda_{C,q}^{\mathcal N}f}{g}=\mathcal B_{C,q}^{\mathcal N}(u_f,g)=\ip{\mathcal Nu_f}{g}.
	\end{equation}
	
	For every nonempty open $W\subset\Omega_e$, we have
	\begin{equation}\label{eq:inverse-anisotropic-Runge}
		\overline{\{u_f|_\Omega:f\in C_c^\infty(W)\}}^{L^2(\Omega)}=L^2(\Omega).
	\end{equation}
	Let $h\in L^2(\Omega)$ be orthogonal to these restrictions, and let $w\in Y_0(\Omega)$ solve \eqref{eq:inverse-source-solution} with $X_0(\Omega)$ and $\mathcal B_{P,q}$ replaced by $Y_0(\Omega)$ and $\mathcal B_{C,q}^{\mathcal N}$. The subtraction giving \eqref{eq:inverse-Runge-exterior-identity} applies to this form and gives $\mathcal B_{C,q}^{\mathcal N}(w,f)=0$ for every $f\in C_c^\infty(W)$. The local and potential integrals vanish for these exterior test functions, giving $\mathcal Nw=0$ on $W$. Since $w=0$ there, Theorem~\ref{thm:single} gives $w=0$ globally, and its source equation gives $h=0$.
	
	\begin{proof}[Proof of Theorem~\ref{thm:inverse-anisotropic}]
		Fix $f\in C_c^\infty(W_1)$ and let $u_{j,f}$ be the solution associated with $C_j,q_j$ for $j=1,2$. Their difference is zero on $\Omega_e$, and \eqref{eq:inverse-anisotropic-data} and \eqref{eq:inverse-anisotropic-weak-DN} give $\mathcal N(u_{1,f}-u_{2,f})=0$ on $W_2$. Theorem~\ref{thm:single} gives $u_{1,f}=u_{2,f}$. Write $u_f=u_{1,f}=u_{2,f}$. Subtracting the two interior equations gives
		\begin{equation}\label{eq:inverse-local-difference}
			Qu_f=0\quad\text{in }\Omega,\quad Q=-\Div(D\nabla)+r,\quad D=C_1-C_2,\quad r=q_1-q_2.
		\end{equation}
		
		To pass from these solutions to arbitrary $L^2$ functions, we place the derivatives on the test function. Since $D$ is smooth and symmetric and $r\in L^\infty(\Omega)$, the distribution $Qv$ for $v\in L^2(\Omega)$ is defined by
		\[
		\ip{Qv}{\phi}=(v,Q\phi)_{L^2(\Omega)},\quad\phi\in C_c^\infty(\Omega),
		\]
		and satisfies $|\ip{Qv}{\phi}|\le\norm{v}_{L^2(\Omega)}\norm{Q\phi}_{L^2(\Omega)}$. Convergence in $L^2(\Omega)$ implies distributional convergence after applying $Q$. Equations~\eqref{eq:inverse-anisotropic-Runge} and \eqref{eq:inverse-local-difference} give $Qv=0$ for every $v\in L^2(\Omega)$.
		
		Fix a ball $B\Subset\Omega$ and choose $\eta\in C_c^\infty(\Omega)$ equal to one on a neighborhood of $\overline B$. Applying $Qv=0$ to $v=\eta$ gives $r=0$ almost everywhere on $B$. Applying it to $v=\eta x_a$ then gives
		\[
		\sum_{\ell=1}^n\partial_\ell D_{\ell a}=0\quad\text{on }B,\quad a=1,\ldots,n.
		\]
		Finally, take $v=\eta x_ax_b$. On $B$, the preceding identities yield
		\[
		Q(\eta x_ax_b)=-\sum_{\ell=1}^n\partial_\ell(D_{\ell a}x_b+D_{\ell b}x_a)=-D_{ba}-D_{ab}=-2D_{ab}.
		\]
		The last equality uses symmetry of $D$. Every entry of $D$ vanishes on $B$. Since the ball was arbitrary, $C_1=C_2$ throughout $\Omega$ and $q_1=q_2$ almost everywhere there.
	\end{proof}
	
	\subsection{Fractional orders and coefficients from one exterior datum}
	
	\begin{lemma}\label{lem:inverse-polynomial}
		Let $Q$ be a nonzero polynomial. If $f\in C_c^\infty(\R^n)$ satisfies $Q(L)f=0$, then $f=0$.
	\end{lemma}
	
	\begin{proof}
		Every integer power of $L$ maps $C_c^\infty(\R^n)$ into itself, so $Q(L)f$ is defined in $L^2$. Let $\mu_f$ be the spectral measure of $f$. The spectral theorem gives
		\[
		0=\norm{Q(L)f}_{L^2}^2=\int_{[0,\infty)}|Q(\lambda)|^2\,d\mu_f(\lambda).
		\]
		The measure $\mu_f$ is supported on the nonnegative real roots of $Q$. If there are no such roots, then $\mu_f=0$ and $f=0$.
		
		Otherwise, let $\rho_1,\ldots,\rho_k$ be the distinct nonnegative real roots. Write $f=\sum_{r=1}^k f_r$, where $f_r$ is the spectral projection of $f$ onto $\ker(L-\rho_r)$. The interpolation polynomial that is one at $\rho_r$ and zero at the other roots gives $f_r=\prod_{\ell\ne r}\frac{L-\rho_\ell}{\rho_r-\rho_\ell}\,f$. The empty product is the identity. This formula shows that each $f_r$ is smooth and compactly supported, since the polynomial involves only integer powers of the local operator $L$. Also $(L-\rho_r)f_r=0$ on $\R^n$. Since $f_r$ vanishes outside a large ball, local unique continuation gives $f_r=0$. Summing proves $f=0$.
	\end{proof}
	
	\begin{proof}[Proof of Theorem~\ref{thm:inverse-orders}]
		For $j=1,2$, let $u_j$ be the solution associated with $P_j+q_j$, with $u_j=f$ on $\Omega_e$. Both belong to $L^2(\R^n)$ and are supported in $K=\overline\Omega\cup\supp f$. Let $\mathcal A=\{\alpha_1,\ldots,\alpha_M\}\subset(0,1)$ be the set of distinct fractional parts of all the orders in \eqref{eq:inverse-unknown-P}. For $\alpha\in\mathcal A$, define
		\begin{equation}\label{eq:inverse-order-polynomials}
			Q_{j,\alpha}(z)=\sum_{\substack{1\le k\le N_j\\s_{j,k}-\alpha\in\Nzero}}b_{j,k}z^{s_{j,k}-\alpha},\quad v_\alpha=Q_{1,\alpha}(L)u_1-Q_{2,\alpha}(L)u_2,
		\end{equation}
		with an empty sum equal to zero. Each $Q_{j,\alpha}$ is a polynomial, and locality gives $\supp v_\alpha\subset K$. Each $v_\alpha$ is a compactly supported distribution and vanishes on $W$.
		
		Set $s_*=\max\{s_{1,N_1},s_{2,N_2}\}$. Since $u_j\in L^2(\R^n)$, every term $L^{s_{j,k}}u_j$ belongs to $\HL^{-2s_*}$. The spectral calculus gives $L^\alpha L^m u_j=L^{\alpha+m}u_j$ in this space whenever $\alpha+m=s_{j,k}$. Grouping terms with the same fractional part gives
		\[
		P_1u_1-P_2u_2=\sum_{\alpha\in\mathcal A}L^\alpha v_\alpha.
		\]
		The left-hand side vanishes on $W$ by \eqref{eq:inverse-orders-data} and \eqref{eq:inverse-weak-DN}. Distinct elements of $\mathcal A\subset(0,1)$ have noninteger differences. Since $\bigcup_{\alpha\in\mathcal A}\supp v_\alpha\subset K$, the unbounded component of $\R^n\setminus\bigcup_{\alpha\in\mathcal A}\supp v_\alpha$ contains $F_f$. We may apply Theorem~\ref{thm:exterior} to the distributions $v_\alpha$, with observation set $W\subset F_f$, and obtain $v_\alpha=0$ on $\R^n$ for every $\alpha\in\mathcal A$.
		
		On $\Omega_e$, locality, $u_1=u_2=f$, and $v_\alpha=0$ give
		\[
		(Q_{1,\alpha}-Q_{2,\alpha})(L)f=0.
		\]
		This expression is smooth and has support contained in $\supp f\Subset\Omega_e$, so its vanishing on $\Omega_e$ implies its vanishing on all of $\R^n$. Since $f\ne0$, Lemma~\ref{lem:inverse-polynomial} forces $Q_{1,\alpha}=Q_{2,\alpha}$ for every $\alpha$.
		
		For each $\alpha\in\mathcal A$ and $m\in\Nzero$, equality of the coefficients of $z^m$ identifies the coefficient of the power $L^{\alpha+m}$ in the two operators. Since the orders within each list are distinct and their coefficients are nonzero, $N_1=N_2$ and the ordered lists of orders and coefficients agree.
		
		Set $P=P_1=P_2$. The difference $w=u_1-u_2$ is zero on $\Omega_e$, and equality of the measurements gives $Pw=0$ on $W$. Theorem~\ref{thm:single} gives $w=0$ on $\R^n$.
	\end{proof}
	
	\begin{proof}[Proof of Corollary~\ref{cor:inverse-simultaneous}]
		Theorem~\ref{thm:inverse-orders} gives $P_1=P_2=:P$ and $u_1=u_2$. Apply the single-datum conclusion \eqref{eq:inverse-single-datum-potential} to this operator $P$, the continuous potentials $q_1,q_2$, and the nonzero exterior datum $f$. It gives $q_1=q_2$ in $\Omega$.
	\end{proof}

	\section*{Statements and declarations}
	
	\para{Data availability statement} No datasets were generated or analyzed in this study.
	
	\para{Conflict of interest} The authors declare no conflicts of interest.
	
	\para{Acknowledgments} The first author is partially supported by the National Science and Technology Council (NSTC), Taiwan, under the project 113-2115-M-A49-017-MY3. The first author also acknowledges financial support from the Alexander von Humboldt Foundation through the Henriette Herz Scouting Programme, hosted by Universit\"at Duisburg-Essen, Germany. The second author was supported by the Hong Kong RGC General Research Funds (projects 11311122, 11304224 and 11303125). The main results, arguments, and statements of this paper were developed without the use of AI. AI tools were used for auxiliary mathematical discussions and assistance with parts of the writing. The authors take full responsibility for the final mathematical statements and proofs.

	\raggedbottom
	\bibliographystyle{alpha}
	\bibliography{refs}
\end{document}